\documentclass[10pt, amsart]{article}
\usepackage{amsmath,amsbsy,amsfonts}
\usepackage[frenchb]{babel}
\title{Paquets d'Arthur discrets pour un groupe classique p-adique}
\newcommand{\nl}[0]{\vskip 0.2cm{}\noindent}
\date{}
\author{\sc C. M{\oe}glin\rm\\Institut de Math\'ematiques de
Jussieu\\ CNRS}
\begin{document}
\maketitle

Ici $F$ est un corps p-adique avec $p\neq 2$ pour les
constructions et $p$ grand pour le th\'eor\`eme prouvant la stabilit\'e. On note $W_F$ le groupe de Weyl de $F$ et on
fixe un groupe classique $G$ (cf. \ref{groupe} ci-dessous
pour plus de pr\'ecision) dont on note $^LG$ le groupe dual
au sens de Langlands. Le but de ce travail est de construire
certains paquets de repr\'esentations de $G$ associ\'es
\`a des morphismes, $\psi$, de
$W_F\times SL(2,
{\mathbb C}\times SL(2,{\mathbb C})$ dans le groupe $^LG$. 
On note le centralisateur de $\psi$ dans $^LG^0$ par $Cent_{^LG}\psi$.
Dans tout ce papier, nous allons supposer que $\psi$ est discret avec pour d\'efinition de ce mot le fait que le groupe $Cent_{^LG}(\psi)$ est fini. 

Chaque repr\'esentation \`a l'int\'erieur du paquet associ\'e \`a $\psi$ va \^etre param\'etris\'ee par un caract\`ere du groupe $Cent_{^LG}(\psi)$, mais elle ne sera, en g\'en\'eral pas irr\'eductible et pourra m\^eme dans le cas le plus g\'en\'eral \^etre nulle. Bien s\^ur ces contructions sont dict\'ees par l'endoscopie, endoscopie qui fait intervenir d'une part les groupes endoscopiques de $G$ mais aussi une forme tordue du $GL$ dans lequel $^LG$ se plonge naturellement; ces id\'ees sont d\^ues \`a Arthur. On doit m\^eme pouvoir v\'erifier que les repr\'esentations ainsi associ\'ees sont uniquement d\'etermin\'ees par les relations venant de cette endoscopie. Toutefois, on laisse cette v\'erification pour un papier ult\'erieur qui est loin d'\^etre \'ecrit et qui a de toute fa\c{c}on besoin des travaux r\'ecents de Waldspurger \cite{waldspurger}.

Pour pouvoir faire la construction annonc\'ee, il faut d\'ej\`a
savoir construire les s\'eries discr\`etes, c'est-\`a-dire
les repr\'esentations associ\'ees \`a des morphismes de la forme $\psi$ d\'ecrite ci-dessus mais
dont la restriction \`a la deuxi\`eme copie de
$SL(2,{\mathbb C})$ est triviale. Ceci est fait dans \cite{europe} et \cite{ams} sous l'hypoth\`ese que l'on connait les points de r\'eductibilit\'e de certaines induites de cuspidales. Ici on a vraiment besoin d'avoir une interpr\'etation de ces points de r\'eductibilit\'e en terme de param\'etrisation de Langlands, ce qui n'\'etait pas vraiment n\'ecessaire dans les r\'ef\'erences ci-dessus; les hypoth\`eses pr\'ecises dont on a besoin sont celles de \cite{paquet} 2.1 on les rappelera encore ici. Soyons quand m\^eme optimiste, chaque fois que l'on a une interpr\'etation en terme d'alg\`ebre de Hecke, le probl\`eme est soluble.

En plus, comme on a en vue une construction avec des propri\'et\'es fortes venant de l'endoscopie, il faudra alors aussi savoir que les s\'eries discr\`etes auxquelles on finit par se ramener ont ces propri\'et\'es. Dans l'\'etat actuel un seul cas est \`a peu pr\`es compl\`etement fait, il s'agit du groupe $SO(2n+1)$ et des morphismes $\psi$ qui se factorisent par le Frobenius, cf \cite{inventiones}, \cite{waldspurger}; l'endoscopie n'est pas \'ecrite, elle n\'ecessite un lemme fondamental qui n'est pas encore connue mais on voit comment on peut faire sans rajouter d'id\'ees par rapport \`a ces r\'ef\'erences. Le cas o\`u $\psi$ est trivial sur le groupe de ramification sauvage est lui aussi en vue avec le m\^eme type d'id\'ee (cf. \cite{niveau0}). Les travaux de \cite{debacker} et \cite{kazhdan} qui sont de m\^eme nature supposent que $\psi$ est trivial sur les 2 copies de $SL(2,{\mathbb C})$ et notre papier dans ce cas ne dit rien puisque qu'alors paquet de Langlands et paquet d'Arthur co\"{\i}ncident.

Ainsi la situation consid\'er\'ee ici n'est pas vide et les constructions faites sont \'ecrites dans un cadre
g\'en\'eral; ce sont surtout des constructions de nature combinatoire.

L'id\'ee est alors la suivante; on note $\Delta$ le plongement
diagonal de $SL(2,{\mathbb C})$ dans le produit
$SL(2,{\mathbb C})\times SL(2,{\mathbb C})$. Ainsi pour $\psi$ comme ci-dessus, le compos\'e $\psi\circ \Delta$ est un morphisme de $W_{F} \times SL(2,{\mathbb C})$ dans $^LG$ et doit donc correspondre \`a un paquet de repr\'esentations temp\'er\'ees. Clairement,  $Cent_{^LG}\psi$ contient le groupe analogue pour $\psi \circ \Delta$.

Un cas est plus facile, et nous l'avons trait\'e en \cite{paquet}, il s'agit du cas o\`u l'inclusions naturelle $Cent_{^LG}(\psi)$ dans $Cent_{^LG}(\psi\circ \Delta)$ est un isomorphisme; on appelle ce cas le cas \'el\'ementaire. Dans ce cas, l'application qui a un caract\`ere $\epsilon$ de $Cent_{^LG}(\psi)$ associe une repr\'esentation $\pi(\psi,\epsilon)$ est injective \`a valeurs dans l'ensemble des repr\'esentations irr\'eductibles et $\psi(\psi,\epsilon)$ se calcule en fonction de $\pi(\psi\circ \Delta, \epsilon)$ par une formule du type de celle qui calcule la g\'en\'eralisation de l'involution d'Iwahori-Matsumoto et donn\'ee par Aubert (\cite{aubert}), Schneider-Stuhler (\cite{SS}); les constructions sont rappel\'ees en \ref{rappel} ci-dessous. 

Le deuxi\`eme cas est celui o\`u $\psi\circ \Delta$ est lui aussi discret; ce cas est celui qui occupe la plus grande partie de ce travail; le but est de se ramener au cas \'el\'ementaire, \'eventuellement pour des groupes de rang plus petit.  Les repr\'esentations du paquet sont donc index\'ees par les caract\`eres, $\epsilon$, du centralisateur de $\psi$. Mais \`a la diff\'erence de ce qui se passe (au moins conjecturalement) pour les s\'eries discr\`etes et le cas de \cite{{paquet}} \`a un tel caract\`ere, on associe une repr\'esentation qui n'est pas irr\'eductible en g\'en\'eral.  Ce ph\'enom\`ene avait d\'ej\`a \'et\'e vu dans le cas archim\'edien par Adams et Johnson et r\'ecemment dans le cas du groupe $G_{2}$ par W.T.Gan et N.Gurevich (\cite{gangurevich}).

Dans la premi\`ere partie de cet article (cf \ref{definitioncommeelement}) on donne une formule pour construire ces repr\'esentations, c'est-\`a-dire associer \`a $\psi$ et $\epsilon$ un \'el\'ement $\pi(\psi,\epsilon)$ du groupe de Grothendieck des repr\'esentations de longueur finie de $G$; avec cette formule, il n'est pas difficile de montrer que la somme:
$$
\sum_{\epsilon}\epsilon(\psi(z))\pi(\psi,\epsilon),
$$o\`u $z$ est l'\'el\'ement non trivial de la deuxi\`eme copie de $SL(2,{\mathbb C})$, est un caract\`ere stable si des combinaisons lin\'eaires de s\'eries discr\`etes associ\'ees essentiellement \`a
$\psi\circ \Delta$ est elle-m\^eme stable (ce que l'on sait dans les cas o\`u on sait contruire ces s\'eries discr\`etes cf. \cite{inventiones}). Cette partie est en fait facile. Voici la d\'efinition dans le groupe de Grothendieck: on plonge $^LG$ dans le $GL$ convenable par la repr\'esentation naturelle et on voit donc $\psi$ comme une repr\'esentation de $W_{F}\times SL(2,{\mathbb C}) \times SL(2,{\mathbb C})$ que l'on d\'ecompose en une somme de repr\'esentations irr\'eductibles de la forme $\rho\otimes Rep_{a}\otimes Rep_{b}$ o\`u $\rho$ est une repr\'esentation irr\'eductible de $W_{F}$ (n\'ecessairement autoduale) et $a,b$ sont des entiers avec une condition de parit\'e dont nous n'avons pas besoin ici. Pour un tel triplet $(\rho,a,b)$ on pose $A:=(a+b)/2-1$ et $B:=\vert a-b\vert$ et on note $\zeta$ le signe de $a-b$ en convenant $\zeta=+$ si $a=b$. On appellera l'ensemble des quadruplets $(\rho,A,B,\zeta)$ ainsi obtenus $Jord(\psi)$ et on voit naturellement $\epsilon$ comme une application de $Jord(\psi)$ dans $\pm 1$. Pour avoir une d\'efinition \`a donner qui sort du cas $\psi$ \'el\'ementaire d\'ej\`a \'etudi\'e en \cite{paquet}, on suppose qu'il existe $(\rho,A,B,\zeta)\in Jord(\psi)$ avec $A>B$; on fixe un tel quadruplet. Gr\^ace \`a la correspondance de Langlands locale pour les groupes lin\'eaires, on identifie $\rho$ \`a une repr\'esentation irr\'eductible d'un groupe lin\'eaire convenable. Pour $C\in ]B,A]$, on note $\delta_{C}$ l'unique sous-module irr\'eductible de l'induite, pour un groupe lin\'eaire convenable, $\rho\vert\,\vert^{\zeta B}\times \cdots \times \rho\vert\,\vert^{-\zeta C}$. Ainsi $\delta_{C}$ est une repr\'esentation de Steinberg tordue si $\zeta=+$ et une repr\'esentation de Speh tordue si $\zeta=-$. On pose $\epsilon_{0}$ la valeur de $\epsilon$ sur $(\rho,A,B,\zeta)$ et on note $\psi',\epsilon'$ l'analogue de $\psi,\epsilon$ en enlevant le bloc $(\rho,A,B,\zeta)$. Ci-dessous, on va rajouter \`a $\psi',\epsilon'$ un ou 2 quadruplets cela veut dire que l'on consid\`ere un couple analogue \`a $\psi,\epsilon$ mais dont les blocs de Jordan sont ceux de $\psi',\epsilon'$ auxquels on rajoute ceux \'ecrits. Il faudrait aussi expliquer $Jac_{\cdots}$, c'est fait en \ref{notation}, disons ici simplement que c'est une partie d'un module de Jacquet. On pose alors:
$$
\pi(\psi,\epsilon)=\oplus_{C\in ]B,A]} (-1)^{A-C}\delta_{C}\times Jac_{\zeta (B+2),  \cdots, \zeta C}\pi(\psi',\epsilon', (\rho,A,B+2,\zeta,\epsilon_{0}))
$$
$$
\oplus_{\eta=\pm}(-1)^{[(A-B+1)/2]}\eta^{A-B+1}\epsilon_{0}^{A-B}\pi(\psi',\epsilon',(\rho,A,B+1,\zeta,\eta),(\rho,B,B,\zeta,\eta\epsilon_{0})),
$$
o\`u par convention si $A=B+1$, $\pi(\psi',\epsilon', (\rho,A,B+2,\zeta,\epsilon_{0}))=\pi(\psi',\epsilon')$ si $\epsilon_{0}=+$ et vaut 0 (c'est-\`a-dire n'appara\^{\i}t pas) si $\epsilon_{0}=-$.
C'est une description par induction car les morphismes qui interviennent dans cette d\'efinition soit sont relatif \`a un groupe de rang plus petit soit si on les notes $\psi''$ v\'erifient $\sum_{(\rho'',A'',B'',\zeta'')\in Jord(\psi'')}A''-B'' < \sum_{(\rho',A',B',\zeta')\in Jord(\psi)}A'-B'$; quand cette somme est nulle, on est pr\'ecis\'ement dans le cas des morphismes \'el\'ementaires.
 Il est assez facile de d\'emontrer que cette d\'efinition ne d\'epend pas du choix de $(\rho,A,B,\zeta)$ v\'erifiant $A>B$.

Dans la deuxi\`eme partie de cet article (\ref{caracterisationcommerepresentation}), on montre que
l'\'el\'ement du groupe de Grothendieck est  une somme sans multiplicit\'e de repr\'esentations irr\'eductibles totalement d\'ecrites par exemple dans \ref{explicite}. La description explicite n'est pas \'evidente et cette partie est vraiment technique. Mais j'ai l'impression que l'on va devoir faire la m\^eme chose pour les groupes lin\'eaires tordus et qu'elle est donc indispensable. Le r\'esultat peut s'exprimer ainsi par r\'ecurrence; on utilise la notation $<\delta,X>$ pour signifier le socle de l'induite $\delta\times X$ quand $X$ est une repr\'esentation semi-simple d'un groupe de m\^eme type que $G$ et $\delta$ une repr\'esentation irr\'eductible d'un groupe lin\'eaire. Avec les notations ci-dessus:
$$
\pi(\psi,\epsilon)=<\delta_{A}, \pi(\psi',\epsilon',(\rho,A-1,B+1,\zeta,\epsilon_{0})> $$
$$
\oplus_{\lambda; \lambda^{A-B+1}=\epsilon_{0}\prod_{C\in [B,C]}(-1)^{[C]}}\pi(\psi',\epsilon', \cup_{C\in [B,A]} (\rho,C,C,\zeta,(-1)^{[C]}\lambda)),
$$
avec comme convention que $\pi(\psi',\epsilon',(\rho,A-1,B+1,\zeta,\epsilon_{0})=0$ si $A=B+1$ et $\epsilon_{0}=-$ et vaut $\pi(\psi',\epsilon')$ si $A=B+1$ et $\epsilon_{0}=+$. On remarque que la somme sur $\lambda$ contient exactement un terme si $A-B$ est pair et en contient soit $0$ soit $2$ si $(A-B)$ est impair la diff\'erence d\'ependant de la valeur de $\epsilon_{0}$. Cette description est par induction comme ci-dessus, on peut en faire une description sans induction mais ce n'est gu\`ere plus limpide (cf \ref{explicite}).

A la fin du papier, on traite le cas des morphismes $\psi$ g\'en\'eraux; on le fait assez rapidement, on les obtient comme module de Jacquet de repr\'esentation associ\'es \`a des morphismes $\psi'$ tel que $\psi'\circ \Delta$ est discret. Ici on peut avoir $\pi(\psi,\epsilon)=0$, on en donne un exemple. Il y aurait donc beaucoup \`a dire sur ces repr\'esentations alors que l'on ne dit essentiellement rien dans ce papier. Mais \`a mon avis, avant de se lancer dans une description forc\'ement technique il faut prouver que l'on a fait les bons choix, c'est-\`a-dire v\'erifier qu'ils sont ceux dict\'es par l'endoscopie ce qui n'est pas fait ici.

L'int\'er\^et de ce travail, s'il en a un, est d'obtenir des cons\'equences pour les formes automorphes de carr\'e int\'egrables. On esp\`ere quand m\^eme que si notre couple $\psi,\epsilon$ est la situation locale  d'une situation globale \`a laquelle est attach\'ee un paquet de repr\'esentations automorphes de carr\'e int\'egrable, les repr\'esentations associ\'ees d\'ecrivent les composantes locales. Rien n'est clair mais ce qui est tr\`es encourageant est que je crois bien que les d\'efinitions sont impos\'ees par l'endoscopie et cette unicit\'e forcera les identifications souhait\'ees.

\section{Pr\'ecisions sur les groupes consid\'er\'es et
conventions}
\subsection{D\'efinition de base\label{groupe}}
 Si $G$ est un groupe
symplectique, il n'y a aucun choix, le groupe est
d\'eploy\'e et son
$L$-groupe est
\`a centre connexe. Et on consid\`ere des morphismes $\psi:
W_F\times SL(2,{\mathbb C})\times SL(2,{\mathbb C})$ dans
$S0(2n+1,{\mathbb C})$ ainsi que des morphismes $\epsilon$
du $Cent_{SO(2n+1,{\mathbb C})}\psi$ dans $\{\pm 1\}$. 

Si $G$ est un groupe des automorphismes d'une forme
orthogonale sur un espace de dimension impaire, $G$ est soit
d\'eploy\'e correspondant
\`a une forme orthogonale isotrope et on pose $\sharp G=+1$
ou $G$ correspond \`a une forme orthognale de noyau
anisotrope de dimension 3 et on pose $\sharp_G=-1$. Le
$L$-groupe est
$Sp(2n)$ si $n$ est le rang de $G$ et son centre est
isomorphe \`a $\{\pm 1\}$. Quand on aura un morphisme,
$\epsilon$ d'un sous-groupe de
$Sp(2n)$ contenant le centre dans $\{\pm 1\}$, on dira que 
$\epsilon_Z=\sharp_G$ quand la restriction de $\epsilon$ \`a
ce centre est triviale si $\sharp_G=+1$ et non triviale
sinon . Pour un tel $G$, on consid\`ere dans ce travail les
morphismes
$\psi$ de $W_F\times SL(2,{\mathbb C}) \times SL(2,{\mathbb
C})$ dans $Sp(2n,{\mathbb C})$  des morphismes $\epsilon$ du
centralisateur de $\psi$ dans $\{\pm 1\}$ tel que
$\epsilon_Z=\sharp_G$;  c'est une traduction terre \`a terre d'un isomorphisme d\^u \`a Kottwitz.

Si $G$ est le groupe orthogonal d'une forme
orthogonale sur un espace de dimension paire, la situation
est plus compliqu\'ee. La forme d\'ependant de son
discriminant qui est une classe de carr\'e, not\'ee ici
seulement $\eta$,  et de l'invariant de Hasse. Dans tout ce
qui suit on consid\`ere des morphismes, $\psi$, fix\'es de
$W_F\times SL(2,{\mathbb C})\times SL(2,{\mathbb C})$ dans le
$L$-groupe; ici le $L$-groupe est par d\'efinition
$O(2n,{\mathbb C})$ si $2n$ est la dimension de la forme
orthogonale. On note $\sharp_G$ le discriminant de la forme.
C'est impropre puisque $G$ ne d\'epend de $\sharp_G$ que si
la classe de carr\'e est celle de $1$. 
 Pour un tel groupe orthogonal, on impose
toujours dans cet article que le d\'eterminant de la
restriction de
$\psi$
\`a $W_F$ qui s'identifie \`a un caract\`ere quadratique de
$W_F$ corresponde \`a la classe de carr\'es $\eta$ par 
r\'eciprocit\'e. Et $\eta$ sera toujours implicitement
fix\'e; il n'intervient pas explicitement. On  consid\`ere
aussi des morphismes $\epsilon$ du centralisateur de $\psi$
dans $\{\pm 1\}$ et, comme dans le cas de dimension impaire,
on demande que $\epsilon_Z=\sharp_G$.

Le lecteur doit  \^etre averti que les subtilit\'es ci-dessus sont en fait cach\'ees dans ce papier; elles n'interviennent que pour la construction des s\'eries discr\`etes que l'on va admettre ci-dessous.

On note $\Delta$ l'application diagonale de $SL(2,{\mathbb C})$
dans $SL(2,{\mathbb C})\times SL(2,{\mathbb C})$. Et dans
tout ce travail, on dit  que le compos\'e $\psi\circ
\Delta$ qui est un morphisme de $W_F\times SL(2,{\mathbb
C})$ est un morphisme discret  si  le centralisateur de  son image est un
groupe fini. 
\subsection{Notations\label{notation}}
On appellera segment dans ce travail un intervalle de la
forme $[x,y]$ o\`u $x,y\in 1/2 \mathbb{Z}$ et $x-y \in
\mathbb{Z}$. Le segment sera consid\'er\'e comme croissant
si $x-y$ est n\'egatif et d\'ecroissant sinon. Dans un tel segment, ce qui nous int\'eresse sont les demi-entiers inclus dans ce segment qui sont non entiers exactement quand $x,y$ sont non entiers.

Soit $X$ une repr\'esentation de longueur finie de $G$ et
$x\in
\mathbb{R}$,
$\rho$ une repr\'esentation cuspidale irr\'eductible
autoduale d'un GL convenable, on note $d_\rho$ le rang de ce
groupe. On suppose que $G$ contient un sous-groupe
parabolique dont le Levi est isomorphe \`a $GL(d_\rho)\times
G(n-d_\rho)$, o\`u $G(n-d_\rho)$ est un groupe de m\^eme
type que $G$ mais de rang $n-d_\rho$. On va d\'efinir
$Jac_xX$ comme un \'el\'ement du groupe de Grothendieck
associ\'e aux repr\'esentations de longueur finie de
$G(n-d_\rho)$; pour cela, on calcule le module de Jacquet de
$X$ le long du radical unipotent d'un parabolique de $G$ de
Levi $GL(d_\rho)\times G(n-d_\rho)$, not\'e momentan\'ement
$X_P$. On regarde $X_P$ d'abord comme repr\'esentation de
$GL(d_\rho)$ et on la projette sur le support cuspidal
$\rho\vert\,\vert^x$; on note $X_x$ cette projection. Tout
sous-quotient irr\'eductible de $X_x$ comme repr\'esentation
de $GL(d_\rho)\times G(n-d_\rho)$ est de la forme
$\rho\vert\,\vert^x\otimes Y$ et $Jac_xX$ est par
d\'efinition la somme des $Y$ qui apparaissent comme de tels
sous-quotient compt\'e avec la multiplicit\'e du
sous-quotient irr\'eductible o\`u ils apparaissent.

Soit $x,y\in \mathbb{R}$, on g\'en\'eralise la d\'efinition
en posant $Jac_{x,y}X:= Jac_y Jac_x X$. On remarque que
$Jac_{x,y}X$ peut se calculer en consid\'erant d'abord le
module de Jacquet de $X$ par rapport au radical unipotent
d'un parabolique de Levi $GL(2d_\rho)$ puis en restreignant
encore et en projettant sur le bon support cuspidal. Si
$\rho\vert\,\vert^x\times \rho\vert\,\vert^y$ est
irr\'eductible, on v\'erifie ainsi que
$$Jac_{x,y}X=Jac_{y,x}X \eqno(1)$$
pour tout $X$. La propri\'et\'e (1) est donc vraie pour tout
couple $x,y$ tel que $x- y\neq \pm 1$.

On pourra encore g\'en\'eraliser \`a plus de 2 projections.
Soit $X$ une
repr\'esentation de longueur finie et $[\alpha,\beta]$ un
segment. On a donc d\'efini $Jac_{x\in [\alpha,\beta]}X$
comme \'etant par d\'efinition $Jac_\beta \cdots
Jac_{\beta-i}\cdots Jac_\alpha X$. 

On utilisera \`a plusieurs reprises la propri\'et\'e
suivante: soit $X$ une repr\'esentation de longueur finie et
soit
${\cal E}$ un ensemble totalement ordonn\'e de demi-entiers
relatifs. Soit $y$ le plus grand \'el\'ement de ${\cal E}$
au sens ordinaire. On suppose que
$Jac_{x\in {\cal E}} X\neq 0$. Alors, on va montrer qu'il existe $z$ un
\'el\'ement de
${\cal E}$ qui pr\'ec\`ede
$y$, tel que $[z,y]$ soit un segment et tel que $Jac_{x\in
[y,z]}X\neq 0$. La m\^eme assertion est vraie si $y$ est le
plus petit \'el\'ement de ${\cal E}$.

En effet, on peut supposer que $X$ est irr\'eductible et que
$y$ est le dernier \'el\'ement de ${\cal E}$ pour l'ordre de
${\cal E}$; ainsi il n'y a pas \`a d\'emontrer que $z$
pr\'ec\`ede $y$. Et la non nullit\'e du module de Jacquet
assure par r\'eciprocit\'e de Frobenius qu'il existe une
sous-repr\'esentation irr\'eductible
$\sigma$ de
$GL(d_\rho\vert{\cal E}\vert)$ et une repr\'esentation
irr\'eductible $X'$ convenable d'un groupe de m\^eme type
que $G$ mais de rang plus petit tel que:
$$
X\hookrightarrow \sigma\times X'. \eqno(2)
$$
On \'ecrit $\sigma$ dans la classification de Zelevinski comme
sous-module d'une induite du type $\times_{i\in [1,v]}
\sigma_i$ o\`u les $\sigma_i$ sont des Steinberg tordues
c'est-\`a-dire des repr\'esentations de la forme
$<\rho\vert\,\vert^{x_i}, \cdots, \rho\vert\,\vert^{y_i}>$
avec $v$ un entier, pour tout $i\in [1,v]$, $[x_i,y_i]$ un
segment et $y_i\geq x_i$, ainsi que $x_1 \geq x_2\geq \cdots
\geq x_v$. De plus ${\cal E}$ est l'union ensembliste de ces
segments. On suppose que $y$ est maximal dans ${\cal E}$.
Ainsi il existe
$i$ tel que
$y_i=y$ et on fixe
$i$ minimum avec cette propri\'et\'e. Ainsi les segments
$[x_j,y_j]$ pour $j<i$ ne sont pas li\'es \`a $[x_i,y_i]$ et
il existe donc une repr\'esentation $\sigma''$ convenable
telle que $\sigma\hookrightarrow <\rho\vert\,\vert^{x_i},
\cdots, \rho\vert\,\vert^{y_i}>\times \sigma''$. Ainsi on a
une inclusion analogue \`a (2) en rempla\c cant $\sigma$ par
$<\rho\vert\,\vert^{x_i},
\cdots, \rho\vert\,\vert^{y_i}>$ et l'assertion avec $z=x_i$
s'en d\'eduit. Si $y$ est minimal au lieu de maximal on
raisonne de la m\^eme fa\c con mais en rempla\c cant
Steinberg par Speh.

\subsection{Classification des param\`etres. \label{classification}}
Quand on a un morphisme $\psi$ comme ci-dessus on le
prolonge en un morphisme dans un groupe lin\'eaire
convenable, $GL(2n)$ pour $G$ un groupe orthogonal de rang
$n$ et $GL(2n+1)$ pour $G$ un groupe symplectique de rang
$n$ et on obtient donc une repr\'esentation de $W_F\times
SL(2,{\mathbb C}) \times SL(2,{\mathbb C})$. On d\'ecompose
cette repr\'esentation en sous-repr\'esentations
irr\'eductibles et comme les repr\'esentations
irr\'eductibles de $SL(2,{\mathbb C})$ sont classifi\'ees
par leur dimension uniquement, une sous-repr\'esentation
irr\'eductible est de la forme $\rho\otimes Rep_a\otimes
 Rep_b$ o\`u $a,b\in \mathbb{N}$, $ Rep_a, Rep_b$ sont
les repr\'esentations irr\'eductibles de $SL(2,{\mathbb C})$
de cette dimension et $\rho$ est une repr\'esentation
irr\'eductible de $W_F$. Ainsi \`a $\psi$ on associe
$Jord(\psi)$ qui par d\'efinition est l'ensemble des triples
$(\rho,a,b)$ tel que $\rho\otimes  Rep_a\otimes Rep_b$
intervienne dans la repr\'esentation associ\'ee \`a $\psi$,
ce triple \'etant compt\'e avec la multiplicit\'e de la
sous-repr\'esentation. On peut faire la m\^eme chose pour
$\psi\circ \Delta$ et $Jord(\psi\circ \Delta)$ se constitue
de couples $(\rho,c)$. On montre ais\'ement que 
$$
Jord(\psi\circ \Delta)= \cup_{(\rho,a,b)\in
Jord(\psi)}\{(\rho, c); c\in [\vert a-b\vert+1,a+b-1]_2\},
$$
o\`u $[?,?]_2$ repr\'esente l'ensemble des \'el\'ements du segment de m\^eme parit\'e que les extr\^emit\'es..

On v\'erifie aussi ais\'ement que $\psi$ (resp. $\psi\circ \Delta$) discret
entra\^{\i}ne que l'ensemble $Jord(\psi)$ (resp. $Jord(\psi\circ \Delta)$) est
sans multiplicit\'e. L'absence de multiplicit\'e force que
toutes les repr\'esentations $\rho$ intervenant doivent
\^etre autoduales et on a aussi  une condition de parit\'e sur
$a+b$ pour tout
$(\rho,a,b)\in Jord(\psi)$ parit\'e d\'etermin\'ee par
$\rho$ et le groupe $G$. On n'a pas besoin d'en savoir plus
sauf qu'\'etant donn\'e un ensemble de triple $(\rho,a,b)$
o\`u
$\rho$ est autodual, $a,b$ ont la bonne condition de
parit\'e et que l'on
a l'\'egalit\'e de dimension $\sum_{(\rho,a,b)\in
Jord(\psi)}d_\rho ab=m^*$ o\`u $d_\rho$ est la dimension de
la repr\'esentation $\rho$ et o\`u $m^*$ est la dimension de
la repr\'esentation naturelle du $L$-groupe de $G$. Traduisons aussi le 
cas o\`u $\psi\circ \Delta$ est discret: cela dit exactement
 que $Jord(\psi)$ est sans multiplicit\'e
et pour tout
$\rho$ fix\'ee autoduale comme ci-dessus et pour tous couples
$(a,b)$,
$(a',b')$ distincts tels que $(\rho,a,b)$ et $(\rho,a',b')
\in Jord(\psi)$, 
$$
[\vert a-b\vert+1,a+b-1] \cap [\vert
a'-b'\vert+1,a'+b'-1]=\emptyset.
$$
Soit $\psi$ comme ci-dessus tel que $\psi\circ \Delta$ soit
discret. Alors le centralisateur de l'image de $\psi$ est isomorphe \`a
$\vert Jord(\psi)\vert$ copies de $\{\pm 1\}$ sauf si
$G=Sp(2n)$ o\`u il faut se restreindre aux \'el\'ements de
produit $+1$; en oubliant cette condition, on calcule le
centralisateur dans $O(2n+1)$ au lieu de $SO(2n+1)$ c'est
\`a dire que l'on a ajout\'e le centre. On identifie un
morphisme de ce centralisateur dans
$\{\pm 1\}$
\`a une application de
$Jord(\psi)$ dans
$\{\pm 1\}$ o\`u la restriction au centre du $L$-groupe
donne la condition:
$$
\times_{(\rho,a,b)}\epsilon(\rho,a,b)=\sharp_G,
$$
o\`u pour unifier, on a pos\'e $\sharp_{Sp(2n)}=+1$ ce qui
compense ''l'oubli'' ci-dessus de se limiter aux
\'el\'ements de d\'eterminant $+1$.

Dans tout le papier, au lieu de fixer $G$, on ne fixe que le
type de $G$, cela d\'etermine $G$ si c'est un groupe
symplectique  et on se donne
$\psi,\epsilon$ comme ci-dessus ce qui d\'etermine la forme
orthogonale dont
$G$ est le groupe d'automorphismes gr\^ace au d\'eterminant
de $\psi$ et \`a la restriction de $\epsilon$ au centre du
$L$-groupe.

Jusqu'au chapitre $\ref{casgeneral}$ on suppose que $\psi\circ \Delta$ est discret.

\nl
\bf{Hypoth\`ese g\'en\'erale}: \sl dans tout ce papier on fixe $\psi$ comme ci-dessus mais on fait l'hypoth\`ese: pour tout groupe $H$ de m\^eme type que $G$ mais de rang inf\'erieur ou \'egal et pour tout morphisme $\psi'$ de $W_{F}\times SL(2,{\mathbb C})$ dans $^LH$ et tout caract\`ere $\epsilon'$ de $Cent_{^LH}(\psi')$ avec les propri\'et\'es que l'on va \'ecrire ci-dessous, on sait associer une repr\'esentation cuspidales $\pi(\psi',\epsilon')$ avec les propri\'et\'es de r\'eductibilit\'e \'ecrites aussi ci-dessous.

Propri\'et\'es sur $(\psi',\epsilon')$: la restriction de $\psi'$ \`a $W_{F}$ est une sous-repr\'esentation de la restriction de $\psi$ \`a $W_{F}$, cela s'entend \'evidemment apr\`es avoir envoy\'e les $L$-groupes dans le groupe lin\'eaire convenable par la repr\'esentation naturelle de ces $L$-groupes. Pour tout $\rho$ repr\'esentation irr\'eductible apparaissant dans la restriction de $\psi'$ \`a $W_{F}$, la repr\'esentation de $SL(2,{\mathbb C})$ dans l'espace isotypique correspondant est une somme de repr\'esentations index\'es par leur dimension $a'\in Jord_{\rho}(\psi')$, ensemble qui doit v\'erifier que s'il contient l'entier $a$ il contient aussi tous les entiers de la forme $a-2i$ avec $i\in [0,[a/2]]$. De plus $\epsilon'$ qui d\'efinit une application de $Jord_{\rho}(\psi')$ dans $\pm 1$ v\'erifie pour tout $a\in Jord_{\rho}(\psi')$ (comme dans la phrase pr\'ec\'edente) $\epsilon'(\rho,a)=(-1)^{[a/2]}\eta$ o\`u $\eta$ est un signe ind\'ependant de $a$ qui vaut +1 dans le cas o\`u $Jord_{\rho}(\psi')$ est form\'e d'entier pair et peut valoir $\pm 1$ dans l'autre cas.

Propri\'et\'es sur $\pi(\psi',\epsilon')$: soit $\rho$ une repr\'esentation cuspidale unitaire irr\'eductible d'un groupe lin\'eaire de la forme $GL(d_{\rho},F)$ (ce qui d\'efinit $d_{\rho}$); gr\^ace \`a la correspondance de Langlands locale, on identifie $\rho$ \`a un morphisme irr\'eductible de $W_{F}$ dans $GL(d_{\rho},{\mathbb C})$; on dit que $\rho$ est symplectique si ce morphisme est \`a valeurs dans $Sp(d_{\rho},{\mathbb C})$ et orthogonal sinon. On suppose que $\rho$ est autoduale; on d\'efinit $Jord_{\rho}(\psi')$ comme on l'a fait ci-dessus, en admettant l'ensemble vide et on note $a_{\rho}$ le plus grand \'el\'ement de cet ensemble s'il est non vide. Dans le cas o\`u $Jord_{\rho}(\psi')$ est vide, on pose $a_{\rho}=0$ si $\rho$ est de type oppos\'e \`a $^LG$ (c'est-\`a-dire orthogonal quand $^LG$ est symplectique et vice et versa) et $a_{\rho}=-1$ sinon. La propri\'et\'e demand\'ee \`a $\pi(\psi',\epsilon')$ est alors que l'induite $\rho\vert\,\vert^s\times \pi(\psi',\epsilon')$ pour $s$ un r\'eel est irr\'eductible sauf exactement quand $s=\pm (a_{\rho}+1)/2$.

\rm

\section{D\'efinition de
$\pi(\psi,\epsilon)$ quand $\psi\circ \Delta$ est discret. \label{definitioncommeelement}}

\subsection{D\'efininitions\label{definition}}
On fixe $\psi,\epsilon$ comme en \ref{classification}, 
 on suppose, comme annonc\'e, que $\psi\circ \Delta$ est discret et
on va associer
\`a ce couple un
\'el\'ement du groupe de Grothendieck des repr\'esentations
de $G$.
 Cette d\'efinition
se fait par induction. On a d\'efini en
\cite{paquet} (cf. \ref{rappel}) cet \'el\'ement dans le cas
o\`u pour tout $(\rho,a,b)\in Jord(\psi)$, $inf(a,b)=1$ et
dans ce cas, c'est une repr\'esentation irr\'eductible. Remarquons que la condition que l'on vient juste d'\'ecrire est exactement \'equivalente \`a ce que $Cent_{^LG}(\psi)$ soit naturellement isomorphe \`a $Cent_{^LG}(\psi\circ \Delta)$.
Supposons donc que cette condition n'est plus v\'erifi\'ee
et fixons
$(\rho,a,b)\in Jord(\psi)$ tel que $inf(a,b)\geq 2$.
On note alors $\psi'$ le morphisme qui se d\'eduit de $\psi$
en enlevant le bloc $(\rho,a,b)$ et on note $\epsilon'$ la
restriction de $\epsilon$ \`a $Jord(\psi')$. On pose
$\zeta_{a,b}$ le signe de $a-b$ si $a\neq b$ et $+$ sinon. On
suppose que l'on a d\'efini $\pi(\psi_1,\epsilon_1)$ pour
tout
$(\psi_1,\epsilon_1)$ relatif \`a un groupe de m\^eme type
que $G$ mais de rang plus petit; on suppose que cette
repr\'esentation est aussi d\'efinie si $\psi_1$ est relatif
\`a $G$ mais si $\vert Jord(\psi_1)\vert$ est strictement
plus grand que $\vert Jord(\psi)\vert$. Remarquons que c'est
une bonne hypoth\`ese de r\'ecurrence puisque le cardinal de
$Jord(\psi)$ est n\'ecessairement born\'e sup\'erieurement. 
On d\'efinit alors dans le groupe de Grothendieck:

\nl
si $inf(a,b)=2$ et
$\epsilon(\rho,a, b)=-$, on pose:
$$\pi(\psi,\epsilon):=\oplus_{\eta=\pm 1} \pi(\psi',\epsilon',
(\rho, a+\zeta_{a,b}1,b-\zeta_{a,b}1,\eta)
(\rho,a-1,b-1,-\eta)),$$ 
o\`u la notation signifie que l'on ajoute \`a $Jord(\psi')$
 les 2 blocs
$(\rho,a+\zeta_{a,b}1,b-\zeta_{a,b}1),(\rho,a-1,b-1)$  en
prolongeant
$\epsilon'$ sur ces blocs par 2 valeurs in\'egales (et on
somme sur les possibilit\'es).

\nl
si $inf(a,b)=2$ et
$\epsilon(\rho,a,b)=+$, on pose:
$$\pi(\psi,\epsilon)= $$
$$<\rho\vert\,\vert^{(a-b)/2}, \cdots,
\rho\vert\,\vert^{-\zeta_{a,b}((a+b)/2-1)}>\times
\pi(\psi',\epsilon')\
\ominus_{\eta=\pm}
\pi(\psi',\epsilon',
(\rho,a+\zeta_{a,b}1,b-\zeta_{a,b}1,\eta),(\rho, a-1,
b-1,\eta)).$$

\nl
si $inf(a,b)\geq 3$, on remarque que
si on ajoute \`a $Jord(\psi')$ le bloc
$(\rho,a+\zeta_{a,b}2, b-\zeta_{a,b}2)$ on d\'efinit un
morphisme pour un groupe de rang $\zeta_{a,b}(a-b)+1$ plus
petit que $G$. On sait donc d\'efinir $\pi(\psi',\epsilon',
(\rho, a+\zeta_{a,b}2, b-\zeta_{a,b}2, \epsilon(\rho,a,b)))$
(o\`u l'on prolonge $\epsilon'$ \`a l'aide de
$\epsilon(\rho,a,b)$). On sait aussi, avec l'hypoth\`ese de
r\'ecurrence d\'efinir pour $\eta=\pm$, $\pi(\psi',\epsilon',
(\rho, a+\zeta_{a,b}1,b-\zeta_{a,b}1, \eta) (\rho, sup(0,
a-b) +1, sup (0, b-a)+1, \eta\epsilon(\rho,a,b)))$. On pose
alors:
\nl
$$
\pi(\psi,\epsilon)=\oplus_{j\in
[1,inf(a,b)[}(-1)^{inf(a,b)-1+j}
$$
$$
<\rho\vert\,\vert^{(a-b)/2},\cdots,
\rho\vert\,\vert^{-(a-b)/2 -\zeta_{a,b}j}> \times
Jac_{(a-b)/2+\zeta_{a,b}2,\cdots, (a-b)/2+\zeta_{a,b}j}
\pi(\psi',\epsilon', (\rho,a+\zeta_{a,b}2, b-\zeta_{a,b}2,
\epsilon(\rho,a,b)))
$$
$$
\oplus_{\eta=\pm}(-1)^{[inf(a,b)/2]}\eta^{inf(a,b)}
\epsilon(\rho,a,b)^{inf(a,b)-1}$$
$$
\pi(\psi',\epsilon',(\rho,
a+\zeta_{a,b}1,b-\zeta_{a,b}1,\eta),(\rho,sup(0,a-b)+1,sup(0
,b-a)+1, \eta \epsilon(\rho,a,b))),
$$quand $j=1$ il n'y a pas de $Jac$ par convention.
\nl
Pour unifier les notations, on peut inclure le cas
$inf(a,b)=2$ dans la d\'efinition g\'en\'erale \`a condition
de poser dans ce cas: $$
\pi(\psi',\epsilon',
(\rho,a+\zeta_{a,b}2,b-\zeta_{a,b}2,\eta)=\begin{cases}
=\pi(\psi',\epsilon') \hbox{ si } \eta=+\\
=0 \hbox{ si }\eta=-.
\end{cases}
$$
R\'ecrivons ces formules en supposant que $a\geq b$ de fa\c{c}on \`a ce qu'elles soient plus lisibles:
$$
\pi(\psi,\epsilon)=\sum_{{j\in [1,b[}}(-1)^{b-1+j}$$
$$
<\rho\vert\,\vert^{(a-b)/2},\cdots,\rho\vert\,\vert^{-(a-b)/2-j}> \times Jac_{(a-b)/2+2, \cdots, (a-b)/2+j} \pi(\psi',\epsilon',(\rho,a+2,b-2,\epsilon(\rho,a,b))$$
$$
\oplus_{\eta=\pm} (-1)^{[b/2]}\eta^b\epsilon(\rho,a,b)^{b-1} \pi(\psi',\epsilon',(\rho,a+1,b-1,\eta), (\rho,a-b+1,0,\eta \epsilon(\rho,a,b)).
$$
\subsection{Autre formulation\label{autreformulation}}
Soit $(\psi,\epsilon)$ comme ci-dessus, on a d\'efini ces objets en utilisant $Jord(\psi)$ en tant qu'ensemble de triplet $(\rho,a,b)$ venant naturellement de la d\'ecomposition en repr\'esentations irr\'eductibles de la repr\'esentation de $W_{F}\times SL(2,{\mathbb C})\times SL(2,{\mathbb C})$. On peut remplacer ces triplets par des quadruplets $(\rho,A,B,\zeta)$ en posant $\zeta=\zeta_{a,b}$ avec nos d\'efinitions ant\'erieurs, $A=(a+b)/2-1$ et $B=\vert a-b\vert/2$; ainsi $A,B$ sont des demi-entiers positifs ou nuls avec $A- B\in {\mathbb N}$, si $B=0$, $\zeta=+$ n\'ecessairement. Remarquons aussi que dans cette correspondance, $inf(a,b)=1$ si et seulement si $A=B$.  En particulier $\psi$ est \'el\'ementaire si $A=B$ pour tout quadruplet de $Jord(\psi)$.

Le fait que $\psi\circ \Delta $ est discret se traduit par le fait que pour tout couple de quadruplets $(\rho,A,B,\zeta)$ et $(\rho',A',B',\zeta') \in Jord(\psi)$ soit $\rho\neq \rho'$ soit $[B,A]\cap [B',A']=\emptyset$ ce qui se traduit encore par le fait que soit $B>A'$ soit $B'>A$.

Notons encore $Jord(\psi)$ l'ensemble de ces quadruplets et reformulons les d\'efinitions  avec ces nouvelles notations.
\nl
Supposons qu'il existe $(\rho,A,B)$ dans $Jord(\psi)$ tel que $A>B$; on pose $\epsilon_{0}:=\epsilon(\rho,A,B,\zeta)$. Alors $$
\pi(\psi,\epsilon)=\oplus_{C\in ]B,A]}(-1)^{A-C}<\rho\vert\,\vert^{\zeta B}, \cdots, \rho\vert\,\vert^{-\zeta C}>\times Jac_{\zeta (B+2), \cdots, \zeta(C)}\pi(\psi',\epsilon',(\rho,A, B+2,\epsilon_{0}))
$$
$$
\oplus_{\eta=\pm}(-1)^{[(A-B+1)/2]}\eta^{A-B+1}\epsilon_{0}^{A-B} \pi(\psi',\epsilon',(\rho,A,B+1,\zeta,\eta)),(\rho,B,B,\zeta,\eta\epsilon_{0})).
$$Ici il faut comme convention que les termes contenant $(\rho,A, B+2,\epsilon_{0})$ n'existent pas si $A=B+1$ et $\epsilon_{0}=-$ et qu'ils font seulement intervenir $\psi',\epsilon'$ si $A=B+1$ et $\epsilon_{0}=+$.
\subsection{Propri\'et\'es\label{propriete}}
La premi\`ere propri\'et\'e \`a remarquer est la suivante: soit $D$ un demi-entier tel que pour tout $C$ tel que $D-C\in {\mathbb Z}_{{\geq 0}}$, il existe un signe $\zeta_{C}$ tel que $(\rho,C,C,\zeta_{C})\in Jord(\psi)$ et si $D$ est demi-entier non entier, alors $\epsilon(\rho,1/2,1/2,\zeta_{1/2})=-$. Alors $\pi(\psi,\epsilon)=\pi(\psi',\epsilon')$ si $Jord(\psi)=Jord(\psi')$ et $\epsilon=\epsilon'$ \`a la seule diff\'erence pr\`es que $\zeta_{C}$ peut ne pas \^etre le m\^eme signe pour $\psi$ et pour $\psi'$ (pour le moment on utilise le m\^eme bloc $(\rho,A,B,\zeta)$ pour les d\'efinitions); cela r\'esulte par r\'ecurrence par exemple de \ref{rappel} propri\'et\'e 1.
\subsubsection{ÊInd\'ependance des choix
\label{independance}}
Soit $(\psi,\epsilon), (\rho,a,b)$ comme ci-dessus ou plut\^ot $(\rho,A,B,\zeta)$ avec $A>B$. On
suppose qu'il existe $(\rho',A',B',\zeta')\in Jord(\psi)$ un bloc de Jordan diff\'erent de $(\rho,A,B,\zeta)$ mais v\'erifiant aussi $A'>B'$. On
aurait pu l'utiliser  pour d\'efinir
$\pi(\psi,\epsilon)$ montrons que l'on aurait obtenu le
m\^eme r\'esultat.

On note $\psi'',\epsilon''$ le morphisme qui se d\'eduit de $\psi,\epsilon$ en enlevant les 2 blocs $(\rho,A,B,\zeta)$ et $(\rho',A',B',\zeta')$ et on pose $\epsilon_{0}=\epsilon(\rho,A,B,\zeta)$, $\epsilon'_{0}=(\rho,A',B',\zeta')$. Pour $C\in ]B,A]$ et $C'\in ]B',A']$, on pose:
$$
\delta_{C}:=<\rho\vert\,\vert^{\zeta B}, \cdots, \rho\vert\,\vert^{-\zeta C}>; \delta_{C'}=<\rho'\vert\,\vert^{\zeta' B'}, \cdots, \rho'\vert\,\vert^{-\zeta' C'}.
$$Le seul cas difficile est le cas o\`u $\rho'=\rho$ ce que nous supposerons pour pouvoir employer sans ambiguit\'e la notation $Jac_{x, \cdots}$.

En utilisant $(\rho,A,B,\zeta)$ pour d\'efinir $\pi(\psi,\epsilon)$ puis $(\rho',A',B',\zeta')$ pour d\'efinir les repr\'esentations qui interviennent dans la d\'efinition de $\pi(\psi,\epsilon)$, on obtient:
$$
\pi(\psi,\epsilon)=\oplus_{C\in ]B,A]}(-1)^{A-C} \oplus_{C'\in ]B',A']}(-1)^{A-C}(-1)^{A'-C'}
$$
$$\delta_{C}
\times\biggl( Jac_{\zeta(B+2), \cdots,\zeta C } \bigl(\delta_{C'}\times Jac_{\zeta'(B'+2),\cdots,\zeta' C'} \pi(\psi'',\epsilon'',(\rho,A,B+2,\zeta,\epsilon_{0}),(\rho,A',B'+2,\zeta',\epsilon'_{0}))\bigr)\biggr)\eqno(1)
$$
$$
\oplus_{C\in ]B,A]}\oplus_{\eta'=\pm}(-1)^{A-C}(-1)^{[(A'-B'+1)/2]}(\eta')^{A'-B'+1}(\epsilon'_{0})^{A'-B'}
$$
$$
\delta_{C}\times Jac_{\zeta (B+2), \cdots,\zeta C }\pi(\psi'',\epsilon'',(\rho,A,B+2,\zeta,\epsilon_{0}) (\rho, A',B'+1,\zeta',\eta'),(\rho, B',B',\zeta',\eta'\epsilon'_{0})) \eqno(2)
$$
$$
\oplus_{C'\in ]B',A']}\oplus_{\eta=\pm}(-1)^{A'-C'}(-1)^{[(A-B+1)/2]}(\eta')^{A-B+1}\epsilon_{0}^{A-B}$$
$$
\delta_{C'}\times Jac_{\zeta (B'+2), \cdots,\zeta C' }\pi(\psi'',\epsilon'',(\rho,A',B'+2,\zeta',\epsilon'_{0}) (\rho, A,B+1,\zeta,\eta),(\rho, B,B,\zeta,\eta\epsilon_{0})) \eqno(3)
$$
$$
\oplus_{\eta=\pm, \eta'=\pm} (-1)^{[(A-B+1)/2]}(\eta')^{A-B+1}\epsilon_{0}^{A-B}
(-1)^{[(A'-B'+1)/2]}(\eta')^{A'-B'+1}(\epsilon'_{0})^{A'-B'}$$
$$
\pi(\psi'',\epsilon'',(\rho, A,B+1,\zeta,\eta),(\rho, B,B,\zeta,\eta\epsilon_{0}),
(\rho, A',B'+1,\zeta',\eta'),(\rho, B',B',\zeta',\eta'\epsilon'_{0}).\eqno(4)
$$
Le terme (4) est parfaitement sym\'etrique en $(\rho,A,B,\zeta)$ et $(\rho,A',B',\zeta')$ et il en est de m\^eme des termes (2) et (3) quand on les consid\`ere ensemble. Il n'y a que le terme (1) qui n'est pas sym\'etrique de fa\c{c}on \'evidente. On va donc montrer qu'en fait il l'est.

Fixons $C\in ]B,A]$ et $C'\in ]B',A']$; l'un des segments $[\zeta B,-\zeta C]$, $[\zeta' B',-\zeta' C']$ contient l'autre  d'o\`u $\delta_{C}\times \delta_{C'}$ est une repr\'esentation irr\'eductible du $GL$ convenable et elle est donc isomorphe \`a $\delta_{C'}\times \delta_{C}$. Les segments $[\zeta (B+2),\zeta C]$ et $[\zeta'(B'+2),\zeta' C']$ sont disjoints et non li\'es, on a donc (cf. \ref{notation}) pour tout \'el\'ement du groupe de Grothendieck $X$
$$
Jac_{\zeta(B+2), \cdots, \zeta C}Jac_{\zeta'(B'+2), \cdots \zeta' C'}X \simeq
Jac_{\zeta'(B'+2), \cdots \zeta' C'}Jac_{\zeta(B+2), \cdots, \zeta C} X.
$$
Montrons que $Jac_{\zeta(B+2), \cdots, \zeta C} \delta_{C'}\times X = \delta_{C'}\times Jac_{\zeta(B+2), \cdots, \zeta C} X$, pour tout $X$ comme ci-dessus; c'est compl\`etement clair quand on a remarqu\'e que  ni $\zeta' B'$ ni $\zeta' C'$ ne sont dans le segment $[\zeta (B+2),\zeta C]$; l\`a on utilise clairement le fait que $\psi\circ \Delta$ est discret. Ainsi la contribution au terme (1) des objets relatifs au couple $(C,C')$ n'est autre que 
$$
(-1)^{A+A'-C-C'}\delta_{C}\times \delta_{C'}\times 
Jac_{\zeta(B+2), \cdots, \zeta C}Jac_{\zeta'(B'+2), \cdots \zeta' C'}
\pi(\psi'',\epsilon'', (\rho,A,B+2,\zeta,\epsilon_{0}),(\rho, A',B'+2,\zeta',\epsilon'_{0})).
$$
Et comme on l'a vu ce terme est parfaitement sym\'etrique en $(A,B,C)$ et $(A',B',C')$. Cela termine la preuve.
\subsection{Stabilit\'e} 
\bf Th\'eor\`eme: \sl On fixe un morphisme $\psi$ de
$W_F\times SL(2,\mathbb{C}\times SL(2,\mathbb{C})$ dans
$^LG$ comme en \ref{groupe}; en particulier on suppose que
$\psi\circ \Delta$ vu comme morphisme de $W_F\times
SL(2,\mathbb{C})$ dans $^LG$ est discret. La combinaison
lin\'eaire de caract\`ere de repr\'esentations:
$$
\sum_{\epsilon: Jord(\psi) \rightarrow \{\pm 1\}}
\prod_{(\rho,a,b)\in Jord(\psi)}\epsilon(\rho,a,b)^{b-1}
\pi(\psi,\epsilon)$$
est stable si   toute combinaison lin\'eaire de s\'eries discr\`etes
$\sum_{\epsilon:Jord(\psi')\rightarrow \{\pm 1\}}
\pi(\psi',\epsilon)$ est stable o\`u $\psi'$ est un morphisme
de $W_F\times SL(2,\mathbb{C})$ dans $^LG'$
avec $G'$ de m\^eme type que $G$ v\'erifiant que
$Jord(\psi')\subset Jord(\psi\circ \Delta)$.\rm
\nl
Ici la terminologie employ\'ee en terme de d\'ecomposition de la repr\'esentation d\'efinie par $\psi$  est plus parlante que celle de \ref{autreformulation}.

On a d\'emontr\'e en \cite{paquet} la stabilit\'e sous
l'hypoth\`ese du th\'eor\`eme dans le cas o\`u pour tout
$(\rho,a,b)\in Jord(\psi)$, $inf(a,b)=1$; sous cette
hypoth\`ese, il suffit m\^eme d'avoir la stabilit\'e pour
$\psi'=\psi\circ
\Delta$. On peut donc supposer que $Jord(\psi)$ contient au
moins un triplet $(\rho,a,b)$ avec $inf(a,b)\geq 2$.
 Supposons d'abord que $inf(a,b)=2$ avec les
notations ci-dessus. Les morphismes qui apparaissent sont
$(\psi',\epsilon')$, on a enlev\'e le bloc $(\rho,a,b)$ et
les morphismes  $\psi_-$ qui se d\'eduit de $\psi$ en
rempla\c cant le bloc $(\rho,a,b)$ par les 2 blocs
$ (\rho, a+\zeta_{a,b}1,b-\zeta_{a,b}1),
(\rho,a-1,b-1)$. On a $Jord(\psi'\circ
\Delta)=Jord(\psi\circ \Delta)-\{(\rho,sup(a,b)+1),
(\rho,sup(a,b)-1)\}$. On peut donc l\'egitimement admettre
le r\'esultat de stabilit\'e pour $\psi'$ par r\'ecurrence
(l'hypoth\`ese du th\'eor\`eme pour $\psi$ est plus forte
que pour $\psi'$). De plus
$Jord(\psi_-\circ \Delta)=Jord(\psi\circ \Delta)$. On
peut aussi admettre la stabilit\'e car l'hypoth\`ese du
th\'eor\`eme est la m\^eme pour $\psi$ et $\psi_-$ mais
$\vert Jord(\psi)\vert =\vert Jord(\psi_-)\vert -1< \vert
Jord(\psi_-)\vert$. On d\'ecompose la somme 
$$
\sum_{\epsilon}\prod_{(\rho',a',b')\in Jord(\psi)}
\epsilon(\rho',a',b')^{b'-1} \pi(\psi,\epsilon)=
$$
$$
\sum_{\begin{array}{c}\epsilon, Jord(\psi)\rightarrow
\{\pm 1\};\\
\epsilon(\rho,a,b)=+
\end{array}}
\prod_{(\rho',a',b')\in Jord(\psi')}
\epsilon(\rho',a',b')^{b'-1}<\rho\vert\,\vert^{(a-b)/2},
\cdots,\rho\vert\,\vert^{\zeta_{a,b}((a+b)/2-1)}>\times \pi(
\psi',\epsilon)
$$
$$
-\sum_{\begin{array}{c}
\epsilon:Jord(\psi)\rightarrow
\{\pm 1\}; \\
\eta=\pm
\end{array}}\epsilon(\rho,a,b)
\prod_{(\rho',a',b')\in
Jord(\psi')}\epsilon(\rho',a',b')^{b'-1} 
\epsilon(\rho,a,b)^{
b-1}\eqno(1)$$
$$
\pi(\psi',\epsilon',
(\rho,a+\zeta_{a,b}1,b-\zeta_{a,b}1,\eta) (\rho,a-1,b-1,\eta
\epsilon(\rho,a,b))).
$$
La premi\`ere somme est stable par l'hypoth\`ese de
r\'ecurrence appliqu\'e \`a $\psi'$. Pour montrer la
stabilit\'e de la deuxi\`eme somme, on va appliquer la
r\'ecurrence au morphisme qui se d\'eduit de $\psi'$ en
ajoutant les 2 blocs $(\rho,a+\zeta_{a,b}1,b-\zeta_{a,b}),(
\rho,a-1,b-1)$.  Pour $\epsilon,\eta$ comme dans la somme, on
construit un morphisme $\epsilon''$ de $ Jord(\psi'')$ dans
$\{\pm 1\}$ comme dans l'\'ecriture ci-dessus,
c'est-\`a-dire en prolongeant la restriction de $\epsilon$
\`a $Jord(\psi')$ par $\eta$ sur $
(\rho,a+\zeta_{a,b}1,b-\zeta_{a,b})$ et par
$\eta\epsilon(\rho,a,b)$ sur $(\rho,a-1,b-1,\eta
\epsilon(\rho,a,b))$.
 Et on calcule
$$\prod_{{(\rho'',a'',b'')\in
Jord(\psi'')}}\epsilon''(\rho'',a'',b'')^{b''-1}=\biggl(\prod_{
(\rho',a',b')\in Jord(\psi)-\{(\rho,a,b)\}}
\epsilon'(\rho',a',b')^{b'-1}\biggr)\bigl(
\eta^{b-\zeta_{a,b}1-1}\eta^{b-2}\epsilon(\rho,a,b)^{b-2}
\bigr).$$
On remarque que quelque soit la valeur de $\zeta_{a,b}$,
$\eta^{b-\zeta_{a,b}1-1}\eta^{b-2}=+$. On obtient donc
exactement le signe qui appara\^{\i}t dans la somme (1)
(quand on n'oublie pas le $\epsilon(\rho,a,b)$ qui se trouve
devant le produit). D'o\`u la stabilit\'e.
\nl
On consid\`ere maintenant le cas o\`u $inf(a,b)>2$. Il faut
\'ecrire
$\sum_{\epsilon}\prod_{(\rho,a,b)\in
Jord(\psi)}\epsilon(
\rho,a,b)^{b-1}\pi(\psi,\epsilon)$ en utilisant
\ref{definition}. Cela donne une somme sur $j\in
[1,inf(a,b)[$ et un terme compl\'ementaire. Fixons d'abord
$j\in [1,inf(a,b)[$, et \'etudions 
$$
\sum_{\epsilon}\prod_{(\rho,a,b)\in
Jord(\psi)}\epsilon(\rho,a,b)^{b-1}
$$
$$
<\rho\vert\,\vert^{(a-b)/2}, \cdots,
\rho\vert\,\vert^{-(a-b) /2-\zeta_{a,b}j}>\times Jac_{ (a-b)
/2+\zeta_{a,b}2,\cdots,
\zeta_{a,b}((a+b)/2-1)}\pi(\psi',\epsilon',
(\rho,a+\zeta_{a, b}2,b-\zeta_{a,b}2,\epsilon(\rho,a,b))).
$$Pour $(\psi,\epsilon)$ fix\'e, on note
$\psi'',\epsilon''$ le couple qui se d\'eduit de
$\psi,\epsilon$ en changeant simplement $(\rho,a,b,
\epsilon(\rho,a,b))$ en
$(\rho,a+2\zeta_{a,b},b-2\zeta_{a,b},\epsilon(\rho,a,b))$.
Comme dans ce changement la parit\'e de $b$ est respect\'ee
il est clair que le signe d\'efini avec $\epsilon$ est le
m\^eme que celui avec $\epsilon''$ quand on remplace $\psi$
par
$\psi''$. Donc la stabilit\'e r\'esulte de l'hypoth\`ese de
r\'ecurrence et du fait qu'elle commute \`a la prise de
module de Jacquet et \`a l'induction.

Il reste le terme compl\'ementaire (on enl\`eve le signe
$(-1)^{[inf(a,b)/2]}$ qui est ind\'ependant des termes de la
somme)
$$
\sum_{\begin{array}{c}
\epsilon: Jac(\psi)\rightarrow \{\pm
1\},\\
\eta=\pm
\end{array}}\epsilon(\rho,a,b)^{inf(b,a)-1}\eta^{inf(a,b)}
\epsilon(\rho,a,b)^{b-1}\biggl(\prod_{(\rho',a',b')\in
Jord(\psi')}\epsilon(\rho',a',b')^{b'-1}\biggr)
$$
$$\pi(\psi',\epsilon', 
(\rho,a+\zeta_{a,b},b-\zeta_{a,b},\eta), (\rho,
sup(a-b,1), sup(b-a,1),\eta\epsilon(\rho,a,b).\eqno(2)$$
Pour $\epsilon, \eta$ fix\'e comme dans la somme, on
d\'efinit $\psi'',\epsilon''$ en rempla\c cant
$(\rho,a,b,\epsilon(\rho,a,b))$ par
$(\rho,a+\zeta_{a,b}1,b-\zeta_{a,b}1,\eta)\cup
(\rho,sup(a-b,1),sup(b-a,1),\eta\epsilon(\rho,a,b))$. Et
il faut encore calculer
$$
\prod_{(\rho'',a'',b'')\in
Jord(\psi'')}\epsilon''(\rho'',a'',b'')^{b''-1}=
\eta^{b-\zeta_{a,b}1-1}\eta^{sup(b-a,0)}\epsilon(\rho,a,b)
^{sup(b-a,0)}
\prod_{(\rho',a',b')\in
Jord(\psi')}\epsilon'(\rho',a',b')^{b'-1}
.\eqno(3)
$$Et on v\'erifie que:
$$
\eta^{b-\zeta_{a,b}1-1}\eta^{sup(b-a,0)}\epsilon(\rho,a,b)
^{sup(b-a,0)}=
\eta^{inf(a,b)}\epsilon(\rho,a,b)^{inf(a,b) -1}
\epsilon(\rho,a,b)^{b-1}$$ c'est-\`a-dire que le signe
intervenant dans (2) est le m\^eme que celui de (3).
On obtient alors la stabilit\'e en appliquant le
th\'eor\`eme \`a $\psi''$ et cela termine la preuve.
\section{Propri\'et\'es des modules de
Jacquet\label{moduledejacquet}} Dans tout ce paragraphe, on
fixe
$\psi,\epsilon$, $\rho$ et
$(\rho,A,B,\zeta)\in Jord(\psi)$ (avec les notations de \ref{autreformulation} et on suppose que
$A>B$. On suppose toujours que $\psi \circ \Delta$ est discret. On raisonne par r\'ecurrence d'abord par r\'ecurrence croissante sur le rang du groupe puis par r\'ecurrence croissante sur $s_{\psi}:=\sum_{(\rho,A,B,\zeta)\in Jord(\psi)}A-B$; en effet si $s_{\psi}=0$, $\psi$ est \'el\'ementaire au sens de \ref{rappel} et la situation sera toujours connue dans ce cas.  Le but de tout ce travail est bien de ramen\'e la d\'efinition g\'en\'erale au cas \'el\'ementaire. On aurait sans doute pu se passer dans cette partie de mettre dans la r\'ecurrence que pour les $\psi',\epsilon'$ plus petit que notre $\psi,\epsilon$ au sens ci-dessus, $\pi(\psi',\epsilon')$ est une vraie repr\'esentation, c'est-\`a-dire une combinaison lin\'eaire \`a coefficients positifs de repr\'esentations irr\'eductibles, mais je mets quand m\^eme cette propri\'et\'e dans la r\'ecurrence pour \'eviter des probl\`emes avec l'interpr\'etation de la non nullit\'e des modules de Jacquet que nous allons consid\'erer; sinon il faudrait se m\'efier des simplifications \'eventuelles. Donc cette section ne sera compl\'et\'ee qu'apr\'es la fin de la section \ref{caracterisationcommerepresentation}.

On d\'emontre simultan\'ement les 2
propositions ci-dessous:
\subsection{\label{calculdujac1}}
Supposons ici que
$B\geq 1$ et que
pour tout \'el\'ement $(\rho,A',B',\zeta')$ de $Jord(\psi)$, on ait $A'\neq B-1$. En particulier, en gardant la notation $(\psi',\epsilon')$ de \ref{definition} le morphisme d\'efinit en ajoutant \`a $Jord(\psi')$ le quadruplet $(\rho, A'-1,B'-1,\zeta')$ est encore discret  et de restriction discr\`ete \`a $W_{F}$ fois la diagonale de $Sl(2,{\mathbb C})$. En terme de repr\'esentation, si $(\rho,a,b)$ correspond \`a $(\rho,A,B,\zeta)$, on enl\`eve 2 au sup de $(a,b)$ apr\`es avoir suppos\'e que $sup(a,b)-inf(a,b)\geq 2$. On pose toujours $\epsilon_{0}=\epsilon(\rho,A,B,\zeta)$.
\nl
{\bf Proposition}: {\sl Supposons que $B \geq 1$ mais pas n\'ecessairement $A>B$ comme ci-dessus. Alors
$$
Jac_{\zeta B,\cdots,
\zeta A)}\pi(\psi,\epsilon)=\pi(\psi',\epsilon' 
,(\rho,A-1,B-1,\zeta,\epsilon_{0})).
$$}
\subsection{\label{proprietedujac}}
{\bf Proposition}: {\sl Soit $x\in \mathbb{R}$. On fixe
$\psi,\epsilon$ et
$\rho$  alors
$Jac_{x}\pi(\psi,\epsilon)=0$ sauf \'eventuellement s'il
existe $ (\rho,A',B',\zeta') \in
Jord(\psi)$ et
$x=\zeta' B'$. De plus pour tout $x \in \mathbb{R}$,
$Jac_{x,x}\pi(\psi,\epsilon)=0$.}
\nl
On les d\'emontre par r\'ecurrence, en utilisant la formule
d\'efinissant $\pi(\psi,\epsilon)$ (cf \ref{definition}). Si $A=B$ dans la proposition \ref{calculdujac1}, cela se fait par r\'ecurrence \`a partir du cas o\`u $\psi$ est \'el\'ementaire (cf. \ref{rappel} propri\'et\'e 2), on ne d\'etaille pas, c'est plus simple que le cas o\`u $A>B$ ce que nous supposerons \`a partir de maintenant.

On fixe donc $(\rho,A,B,\zeta,\epsilon_{0}$ dans $Jord(\psi,\epsilon)$ tel que $A>B$ et pour le moment, on n'impose rien \`a $B$. Soit $C\in ]B,A]$. On doit calculer:
$$
Jac_{\zeta B, \cdots, \zeta A}\biggl( <\rho\vert\,\vert^{\zeta B}, \cdots, \rho\vert\,\vert{ -\zeta C}> \times Jac_{\zeta(B+2), \cdots, \zeta C}\pi(\psi',\epsilon',(\rho,A,B+2,\zeta,\epsilon_{0}))\biggr). \eqno(1)
$$
On a $Jac_{\zeta B} Jac_{\zeta( B+2), \cdots, \zeta C}\pi(\psi',\epsilon',(\rho,A,B+2,\zeta,\epsilon_{0}))=$
$Jac_{\zeta( B+2), \cdots, \zeta C}Jac_{\zeta B}\pi(\psi',\epsilon',(\rho,A,B+2,\zeta,\epsilon_{0}))=0$ par la \ref{calculdujac1} appliqu\'ee par r\'ecurrence. Ainsi 
$$
(1)= Jac_{\zeta (B+1), \cdots, \zeta A}\biggl( <\rho\vert\,\vert^{\zeta (B-1)}, \cdots, \rho\vert\,\vert{ -\zeta C}> \times Jac_{\zeta(B+2), \cdots, \zeta C}\pi(\psi',\epsilon',(\rho,A,B+2,\zeta,\epsilon_{0}))\biggr) $$ce qui vaut encore par les calculs standard de modules de Jacquet:
$$
Jac_{\zeta C, \cdots, \zeta A}
\biggl( <\rho\vert\,\vert^{\zeta (B-1)}, \cdots,\rho\vert\,\vert{ -\zeta C}> \times Jac_{\zeta(B+1), \cdots, \zeta(C-1)}Jac_{\zeta(B+2), \cdots, \zeta C}\pi(\psi',\epsilon',(\rho,A,B+2,\zeta,\epsilon_{0}))\biggr).
$$
De plus on a par d\'efinition, pour tout \'el\'ement $X$ du groupe de Grothendieck:
$$
Jac_{\zeta C}Jac_{\zeta(B+1), \cdots, \zeta(C-1)}Jac_{\zeta(B+2), \cdots, \zeta C}X=
Jac_{\zeta(B+1), \cdots, \zeta C}Jac_{\zeta(B+2), \cdots, \zeta C}X.
$$
Par calculer le $Jac$ qui nous int\'eresse, on calcule tous les modules de Jacquet des constituants de $\pi(\psi',\epsilon', \cdots)$ le long du parabolique de Levi $GL((2(C-B)-1)d_{\rho},F) \times G'$, o\`u $G'$ est un groupe de m\^eme type que $G$ mais de rang plus petit. On regarde ces modules de Jacquet dans le groupe de Grothendieck ce sont donc des combinaisons lin\'eaires  de repr\'esentations irr\'eductibles du type $\sigma'\otimes X'$. Ensuite on projette sur le support cuspidal qui nous int\'eresse. Pour faire ce calcul, on reprend les classifications de Zelevinsky des $\sigma'$ soit par segment croissant si $\zeta=+$ soit par segment d\'ecroissant si $\zeta=-$ de telle sorte que la seule repr\'esentation $\sigma'$ qui contribue soit $$
<\rho\vert\,\vert^{\zeta (B+1) },\cdots, \rho\vert\,\vert^{\zeta C}>\times <\rho\vert\,\vert^{\zeta (B+2) },\cdots, \rho\vert\,\vert^{\zeta C}>.
$$
Mais cette induite est irr\'eductible donc isomorphe \`a $$
<\rho\vert\,\vert^{\zeta (B+2) },\cdots, \rho\vert\,\vert^{\zeta C}> \times <\rho\vert\,\vert^{\zeta (B+1) },\cdots, \rho\vert\,\vert^{\zeta C}>.
$$
Mais cela montre que l'on aurait obtenu le m\^eme r\'esultat en calculant $$Jac_{\zeta(B+1), \cdots, \zeta C}Jac_{\zeta(B+2), \cdots, \zeta C}Jac_{\zeta(B+1), \cdots, \zeta C} \pi(\psi',\epsilon',(\rho,A,B+2,\zeta,\epsilon_{0})).$$ Mais ceci fait z\'ero car ce calcul se ''factorise'' par $Jac_{\zeta (B+1)}\pi(\psi',\epsilon',(\rho,A,B+2,\zeta,\epsilon_{0}))$ dont la proposition \ref{proprietedujac} appliqu\'ee par r\'ecurrence donne la nullit\'e.
D'o\`u
$$
(1)=Jac_{\zeta(C+1), \cdots, \zeta A}<\rho\vert\,\vert^{\zeta(B-1)}, \cdots, \rho\vert\,\vert^{-\zeta(C-1)}> \times Jac_{\zeta(B+1), \cdots, \zeta{(C-1)}}Jac_{\zeta(B+2), \cdots \zeta C}\pi(\psi',\epsilon',(\rho,A,B+2,\zeta,\epsilon_{0})).$$Or ni $\zeta(B-1)$ ni $\zeta(C-1)$ ne sont des \'el\'ements du segment $[\zeta(C+1),\zeta A]$ et on obtient donc par les calculs standard de modules de Jacquet (le rappel de \ref{standard} suffit ici):
$$
(1)=<\rho\vert\,\vert^{\zeta(B-1)}, \cdots, \rho\vert\,\vert^{-\zeta(C-1)}> \times
Jac_{\zeta(C+1), \cdots, \zeta A}
Jac_{\zeta(B+1), \cdots, \zeta{(C-1)}}Jac_{\zeta(B+2), \cdots \zeta C}\pi(\psi',\epsilon',(\rho,A,B+2,\zeta,\epsilon_{0})).$$
Or de fa\c{c}on assez formelle:
$$
Jac_{\zeta(C+1), \cdots, \zeta A}
Jac_{\zeta(B+1), \cdots, \zeta{(C-1)}}Jac_{\zeta(B+2), \cdots \zeta C}=Jac_{\zeta(B+1), \cdots, \zeta (C-1)}Jac_{\zeta(B+2), \cdots, \zeta A}$$
Ce qui donne encore:
$$(1)=<\rho\vert\,\vert^{\zeta(B-1)}, \cdots, \rho\vert\,\vert^{-\zeta(C-1)}> \times
Jac_{\zeta(B+1), \cdots, \zeta (C-1)}Jac_{\zeta(B+2), \cdots, \zeta A}
\pi(\psi',\epsilon',(\rho,A,B+2,\zeta,\epsilon_{0})).
$$
On veut encore remplacer $Jac_{\zeta(B+2), \cdots, \zeta A}
\pi(\psi',\epsilon',(\rho,A,B+2,\zeta,\epsilon_{0}))$ par $\pi(\psi',\epsilon',(\rho,A-1,B+1,\zeta,\epsilon_{0})$; si $A=B+1$ de toute fa\c{c}on ce terme n'intervenait pas, si $A\geq B+2$, c'est la proposition \ref{calculdujac1} que l'on applique par r\'ecurrence au bloc $(\rho,A,B+2,\zeta)$ qui en v\'erifie les hypoth\`eses.
D'o\`u finalement
$$
(1)= <\rho\vert\,\vert^{\zeta(B-1)}, \cdots, \rho\vert\,\vert^{-\zeta(C-1)}> \times
Jac_{\zeta(B+1), \cdots, \zeta (C-1)}\pi(\psi',\epsilon',(\rho,A-1,B+1,\zeta,\epsilon_{0})).\eqno(2)
$$
C'est-\`a-dire qu'en calculant $Jac_{\zeta B, \cdots, \zeta A}$ on a remplac\'e $(A,B,C)$ par $(A-1,B-1,C-1)$.

On a aussi par d\'efinition, pour $\eta=\pm$,
$$
Jac_{\zeta B,\cdots,
\zeta A}\pi(\psi',\epsilon',
(\rho,A,B+1,\zeta,\eta
),(\rho,B,B,\zeta, \eta\epsilon_{0}))=
$$
$$
 Jac_{\zeta(B+1),\cdots,
\zeta A}\, Jac_{\zeta B}
$$
$$\pi(\psi',\epsilon',
(\rho,A,B+1,\zeta,\eta
),(\rho,B,B,\zeta, \eta\epsilon_{0})).
$$
Supposons maintenant que les hypoth\`eses de la proposition \ref{calculdujac1} soient satisfaites. En particulier le bloc
$(\rho,B,B,\zeta)$ les v\'erifie; on a alors par r\'ecurrence$$
Jac_{\zeta B}\pi(\psi',\epsilon',
(\rho,A,B+1,\zeta,\eta
),(\rho,B,B,\zeta, \eta\epsilon_{0}))=
\pi(\psi',\epsilon',
(\rho,A-1,B+1,\zeta,\eta
),(\rho,B-1,B-1,\zeta,
\eta\epsilon_{0}).
$$
On en d\'eduit:
$$
Jac_{\zeta B,\cdots, \zeta A)}\pi(\psi',\epsilon', (\rho,A, B+1,\eta),(\rho,B,B,\zeta, \eta\epsilon_{0})=
$$
$$ Jac_{\zeta(B+1),\cdots,
\zeta A}\pi(\psi',\epsilon',
(\rho,A,B+1,\zeta,\eta
),(\rho,B-1,B-1,\zeta,
\eta\epsilon_{0})),
$$
et en appliquant la proposition \ref{calculdujac1} par r\'ecurrence ce que
l'on a le droit de faire
$$
= \pi(\psi',\epsilon',(\rho,A-1,B,\zeta, \eta),(\rho, B-1,B-1,\zeta,\eta
\epsilon_{0})),\eqno(3)
$$c'est-\`a-dire qu'ici aussi, on a obtenu le terme de d\'epart mais en rempla\c{c}ant le couple $(A,B)$ par le couple $(A-1,B-1)$. Ainsi
en mettant ensemble les termes obtenus en (2) ci-dessus sans aucune hypoth\`ese sur $(A,B)$ et les termes (3) qui eux ont \'et\'e obtenus modulo les hypoth\`eses de \ref{calculdujac1} on obtient exactement la d\'efinition de $\pi(\psi',\epsilon',(\rho,A-1,B-1,\zeta,\epsilon_{0}))$ qui elle n\'ecessite aussi l'hypoth\`ese que le morphsime $(\psi',\epsilon',(\rho,A-1,B-1,\zeta,\epsilon_{0}))$ soit de restriction discr\`ete \`a la diagonale c'est-\`a-dire l'hypoth\`ese de \ref{calculdujac1}.   Cela termine la preuve de cette proposition.
\nl
Montrons maintenant la deuxi\`eme proposition; cette proposition a \'et\'e d\'emontr\'ee dans le cas o\`u $\psi$ est \'el\'ementaire (cf. \ref{rappel}) et on suppose donc que $\psi$ n'est pas \'el\'ementaire. On fixe $(\rho,A,B,\zeta)\in Jord(\psi)$ tel que $A>B$ et on pose encore $\epsilon_{0}:=\epsilon(\rho,A,B,\zeta)$.

Soit $x$ tel que $Jac_x \pi(\psi,\epsilon)\neq 0$; n\'ecessairement, en revenant \`a la d\'efinition \ref{definition}

soit il existe $\eta=\pm$ tel que $Jac_x
\pi(\psi',\epsilon',
(\rho,A,B+1,\zeta,\eta),(\rho,B,B,\zeta,\eta\epsilon(\rho,a,b)))\neq 
0$

soit il existe $C\in ]B,A]$ tel que $$Jac_x \biggl(
<\rho\vert\,\vert^{\zeta B},\cdots, \rho\vert\,\vert^{-\zeta C}> \times Jac_{\zeta (B+2),\cdots, \zeta C}
\pi(\psi',\epsilon', (\rho,A,B+2,\epsilon_{0}))\biggr)\neq
0.$$
\nl 
Pour le premier cas, les valeurs de $x$ qui ne satisfont
pas \`a la proposition sont uniquement $\zeta (B+1)$. Dans le
deuxi\`eme cas, celles qui ne
satisfont pas sont $\zeta C$ avec $C\in ]B,A]$ et celles
qui v\'erifient $$Jac_{\zeta (B+2), \cdots, \zeta C, x}
\pi(\psi',\epsilon',(\rho,A,B+2,\zeta ,\epsilon_{0}))\neq
0.$$ 
Dans  cette derni\`ere \'eventualit\'e, on v\'erifie en utilisant la proposition elle-m\^eme par r\'ecurrence que les seules difficult\'es viennent de $x=\zeta(B+1)$ et $x=\zeta (C+1)$.  On va encore
v\'erifier que
$x=\zeta (C+1)$ avec
$C=A$ n'est pas une exception; en effet si $C=A$, on applique \ref{calculdujac1} \`a $(\rho,A,B+2,\zeta)$ pour obtenir
faut $$
 Jac_{\zeta(B+2)\cdots,
\zeta A}\pi(\psi',\epsilon',(\rho,A,B+2,\zeta,\epsilon_{0}))
= \pi(\psi',\epsilon',
(\rho,A-1,B+1,\epsilon(\rho,a,b))) 
$$et en appliquand \ref{proprietedujac} par r\'ecurrence, $Jac_{\zeta (A+1)}\pi(\psi',\epsilon',
(\rho,A-1,B+1,\epsilon(\rho,a,b))) \neq 0$ entra\^{\i}ne qu'il existe $(\rho,A',B',\zeta')\in Jord(\psi)$ avec $\zeta(A+1)=\zeta' B'$ c'est-\`a-dire la propri\'et\'e cherch\'ee.

En r\'esum\'e, les
valeurs qui nous g\^enent sont $x=\zeta (C+1)$ avec $C\in [B,A[$ (c'est-\`a-dire $B$ inclus et $A$ exclu). Supposons d'abord que $C\in ]B,A[$. Alors les 
termes de la somme d\'efinissant $\pi(\psi,\epsilon)$ qui
contribuent au $Jac_{\zeta (C+1)}$ sont exactement au nombre de
2, celui correspondant \`a $C+1$ et celui correspondant \`a
$C$. La contribution du terme correspondant \`a $C$ est:
$$(-1)^{A-C}Jac_{\zeta (C+1)}\biggl(<\rho\vert\,\vert^{\zeta B}, \cdots,\rho\vert\,\vert^{-\zeta C}> \times
Jac_{\zeta(B+2), \cdots, \zeta C}\pi(\psi',\epsilon',(\rho,A,B+2,\zeta,\epsilon_{0}))\biggr)=$$
$$
(-1)^{A-C}<\rho\vert\,\vert^{\zeta B}, \cdots,\rho\vert\,\vert^{-\zeta C}> \times
Jac_{\zeta(B+2), \cdots, \zeta C,\zeta (C+1)}\pi(\psi',\epsilon',(\rho,A,B+2,\zeta,\epsilon_{0}))$$et la contribution du terme correspondant \`a $C+1$ est
$$
(-1)^{A-C-1}
Jac_{\zeta (C+1)}\biggl(<\rho\vert\,\vert^{\zeta B}, \cdots,\rho\vert\,\vert^{-\zeta (C+1)}> \times
Jac_{\zeta(B+2), \cdots, \zeta (C+1)}\pi(\psi',\epsilon',(\rho,A,B+2,\zeta,\epsilon_{0}))\biggr)=$$
$$
(-1)^{A-C-1}
(<\rho\vert\,\vert^{\zeta B}, \cdots,\rho\vert\,\vert^{-\zeta C}> \times
Jac_{\zeta(B+2), \cdots, \zeta (C+1)}\pi(\psi',\epsilon',(\rho,A,B+2,\zeta,\epsilon_{0})).
$$
La somme des 2 contributions fait donc 0 comme annonc\'e. Reste \`a voir le cas o\`u $C=B$ c'est-\`a-dire \`a calculer $Jac_{\zeta(B+1)}\pi(\psi,\epsilon)$.  
Le calcul est
assez difficile car tous les termes contribuent. On fait d'abord le cas o\`u $A=B+1$. Dans ce cas, si $\epsilon_{0}=+$
$$\pi(\psi,\epsilon)= <\rho\vert\,\vert^{\zeta B}, \cdots,   \rho\vert\,\vert ^{-\zeta(B+1)}>\times  \pi(\psi',\epsilon')\ominus_{\eta=\pm}\pi(\psi',\epsilon',(\rho,B+1,B+1,\zeta,\eta),(\rho,B,B,\zeta,\eta))$$
et si $\epsilon_{0}=-$
$$
\pi(\psi,\epsilon)=\oplus_{\eta=\pm }\pi(\psi',\epsilon',(\rho,B+1,B+1,\zeta,\eta),(\rho,B,B,\zeta,-\eta)).
$$
On v\'erifie avec \ref{rappel} que dans le deuxi\`eme cas $Jac_{\zeta(B+1)}=0$. Dans le cas o\`u $\epsilon_{0}=+$, avec \ref{rappel} on v\'erifie que
$$
Jac_{\zeta (B+1)}(\oplus_{\eta=\pm}\pi(\psi',\epsilon',(\rho,A=B+1,A=B+1,\zeta,\eta),(\rho,B,B,\zeta,\eta))=
<\rho\vert\,\vert^{\zeta B}, \cdots,   \rho\vert\,\vert ^{-\zeta B}>\times  \pi(\psi',\epsilon').
$$
Il est alors clair"
 que $Jac_{\zeta (B+1)}\pi(\psi,\epsilon)=0$ dans ce cas aussi.

On suppose maintenant que $A=B+2$. Ici on a:
$$
\pi(\psi,\epsilon)=<\rho\vert\,\vert^{\zeta B}, \cdots, \rho\vert\,\vert^{-\zeta(B+2)}> \times Jac_{\zeta (B+2)}\pi(\psi',\epsilon',(\rho,B+2,B+2,\zeta,\epsilon_{0}))\eqno(4)
$$
$$
- <\rho\vert\,\vert^{\zeta B}, \cdots, \rho\vert\,\vert^{-\zeta(B+1)}> \times 
\pi(\psi',\epsilon',(\rho,B+2,B+2,\zeta,\epsilon_{0}))\eqno(5)
$$
$$
\ominus_{\eta=\pm}\eta \pi(\psi',\epsilon',(\rho,B+2,B+1,\zeta,\eta),(\rho,B,B,\zeta,\eta\epsilon_{0})).\eqno(6)
$$Dans (6) on remplace $\pi(\psi',\epsilon',(\rho,B+2,B+1,\zeta,\eta),(\rho,B,B,\zeta,\eta\epsilon_{0}))$ par sa d\'efinition en utilisant $(\rho,B+2,B+1,\zeta)$ et cette d\'efinition d\'epend de la valeur de $\eta$. Ainsi 
$$
(6)=\ominus <\rho \vert\,\vert^{\zeta (B+1)}, \cdots, \rho\vert\,\vert^{-\zeta (B+2)}\times \pi(\psi',\epsilon',(\rho,B,B,\zeta,\epsilon_{0})) $$
$$\oplus_{\eta=\pm,\eta'=\pm} \pi(\psi',\epsilon',(\rho,B+2,B+2,\zeta, \eta'), (\rho,B+1,B+1,\zeta,\eta\eta'),(\rho,B,B,\zeta,\eta \epsilon_{0})). \eqno(7)
$$
Quand on applique $Jac_{\zeta (B+1)}$ \`a (7), il n'y a pas de contribution des termes tels que $\eta'\neq \epsilon_{0}$ et  on obtient:
$$
\ominus <\rho \vert\,\vert^{\zeta B}, \cdots, \rho\vert\,\vert^{-\zeta (B+2)}\times \pi(\psi',\epsilon',(\rho,B,B,\zeta,\epsilon_{0})) \eqno(8)
$$
$$
\oplus <\rho\vert\,\vert^{\zeta B}, \cdots, \rho\vert\,\vert^{-\zeta B}> \times \pi(\psi',\epsilon', (\rho, B+2,B+2, \zeta, \epsilon_{0})).\eqno(9)
$$
Il est clair que (9)   annule la contribution de $Jac_{\zeta(B+1)}$ appliqu\'e \`a (5). On calcule $Jac_{\zeta(B+1)}$ appliqu\'e \`a (4) en utilisant  2 fois \ref{calculdujac1} d'abord pour $(\rho,B+2,B+2,\zeta)$ puis pour $(\rho,B+1,B+1,\zeta)$
$$
Jac_{\zeta (B+1)}Jac_{\zeta (B+2)}\pi(\psi',\epsilon',(\rho,B+2,B+2,\zeta,\epsilon_{0}))= Jac_{\zeta(B+1)}
\pi(\psi',\epsilon', (\rho,B+1,B+1,\zeta,\epsilon_{0})=$$
$$
\pi(\psi',\epsilon',(\rho,B,B,\zeta,\epsilon_{0}).
$$Et ainsi (8) annule $Jac_{\zeta(B+1)}$ appliqu\'e \`a (4). Cela termine la preuve de ce cas.  

On fait maintenant le cas g\'en\'eral qui est de m\^eme nature. On suppose donc que $A>B+2$ et on \'ecrit:
$$
\pi(\psi,\epsilon)=
\oplus_{C\in ]B,A]} (-1)^{A-C}<\rho\vert\,\vert^{\zeta B}, \cdots, \rho\vert\,\vert^{-\zeta C}>\times Jac_{\zeta(B+2), \cdots, \zeta C}\pi(\psi',\epsilon',(\rho,A,B+2,\zeta,\epsilon_{0}))\eqno(*_{C})
$$
$$
\oplus_{\eta=\pm}(-1)^{[(A-B+1)/2]}\eta^{A-B+1}\epsilon_{0}^{A-B}$$
$$
\pi(\psi',\epsilon',(\rho,A,B+1,\zeta,\eta),(\rho,B,B,\zeta,\eta\epsilon_{0})).\eqno(*_{\eta}).
$$On fixe $C>B+1$ et on r\'ecrit $(*_{C})$ en rempla\c{c}ant $\pi(\psi',\epsilon',(\rho,A,B+2,\zeta,\epsilon_{0}))$ par sa d\'efinition, en  utilisant $(\rho,A,B+2,\zeta)$.
D'o\`u
$$
(*_{C}):=(-1)^{A-C}\oplus_{C'\in ]B+2,A]}(-1)^{A-C'}<\rho\vert\,\vert^{\zeta B}, \cdots, \rho\vert\,\vert^{-\zeta C}>\times $$
$$
Jac_{\zeta(B+2), \cdots, \zeta C}\biggl( <\rho\vert\,\vert^{\zeta (B+2)}, \cdots, \rho\vert\,\vert^{-\zeta C'}>\times Jac_{\zeta(B+4), \cdots, \zeta C'}\pi(\psi',\epsilon',(\rho,A,B+4,\zeta,\epsilon_{0}))\biggr)\eqno(\dag_{C})
$$
$$
(-1)^{A-C}\oplus_{\eta=\pm}(-1)^{(A-B-1)/2]}\eta^{A-B-1}\epsilon_{0}^{A-B}
<\rho\vert\,\vert^{\zeta B}, \cdots, \rho\vert\,\vert^{-\zeta C}>\times $$
$$
Jac_{\zeta(B+2), \cdots, \zeta C}\pi(\psi',\epsilon',(\rho,A,B+3,\zeta,\eta),(\rho,B+2,B+2,\zeta,\eta\epsilon_{0})).\eqno(\dag'_{C})
$$
Dans les formules ci-dessus, $Jac_{\zeta(B+2)}$ se calcule ais\'ement car on a suppos\'e que $C>B+1$
$$
(\dag_{C})= 
(-1)^{A-C}\oplus_{C'\in ]B+2,A]}(-1)^{A-C'}
<\rho\vert\,\vert^{\zeta B}, \cdots, \rho\vert\,\vert^{-\zeta C}>\times $$
$$
Jac_{\zeta(B+3), \cdots, \zeta C}\biggl( <\rho\vert\,\vert^{\zeta (B+1)}, \cdots, \rho\vert\,\vert^{-\zeta C'}>\times Jac_{\zeta(B+4), \cdots, \zeta C'}\pi(\psi',\epsilon',(\rho,A,B+4,\zeta,\epsilon_{0}))\biggr)$$et quand on applique encore $Jac_{\zeta(B+1)}$ si $C\neq B+1$ cela donne encore, pour $C\neq B+1$:
$$
Jac_{\zeta(B+1)}(\dag_{C})= 
(-1)^{A-C}\oplus_{C'\in ]B+2,A]}(-1)^{A-C'}
<\rho\vert\,\vert^{\zeta B}, \cdots, \rho\vert\,\vert^{-\zeta C}>\times $$
$$
Jac_{\zeta(B+3), \cdots, \zeta C}\biggl( <\rho\vert\,\vert^{\zeta B}, \cdots, \rho\vert\,\vert^{-\zeta C'}>\times Jac_{\zeta(B+4), \cdots, \zeta C'}\pi(\psi',\epsilon',(\rho,A,B+4,\zeta,\epsilon_{0}))\biggr)$$
Supposons que $C'>C$, le terme correspondant au couple $(C,C')$ se r\'ecrit:
$$
(-1)^{A-C}(-1)^{A-C'}
<\rho\vert\,\vert^{\zeta B}, \cdots, \rho\vert\,\vert^{-\zeta C}>\times 
<\rho\vert\,\vert^{\zeta B}, \cdots, \rho\vert\,\vert^{-\zeta C'}>\times $$
$$
Jac_{\zeta(B+3), \cdots, \zeta C}Jac_{\zeta(B+4), \cdots, \zeta C'}\pi(\psi',\epsilon',(\rho,A,B+4,\zeta,\epsilon_{0})).
$$
Supposons maintenant que $C'\leq C$, on v\'erifie que cela se r\'ecrit:
$$
(-1)^{A-C}(-1)^{A-C'}
<\rho\vert\,\vert^{\zeta B}, \cdots, \rho\vert\,\vert^{-\zeta C}>\times 
<\rho\vert\,\vert^{\zeta B}, \cdots, \rho\vert\,\vert^{-\zeta (C'-1)}>\times $$
$$
Jac_{\zeta (B+3), \cdots, \zeta (C'-1),\zeta (C'+1), \cdots, \zeta C}
Jac_{\zeta(B+4), \cdots, \zeta C'}\pi(\psi',\epsilon',(\rho,A,B+4,\zeta,\epsilon_{0}))
$$
ce qui est encore:
$$
(-1)^{A-C}(-1)^{A-C'}
<\rho\vert\,\vert^{\zeta B}, \cdots, \rho\vert\,\vert^{-\zeta C}>\times 
<\rho\vert\,\vert^{\zeta B}, \cdots, \rho\vert\,\vert^{-\zeta (C'-1)}>\times $$
$$
Jac_{\zeta (B+3), \cdots, \zeta (C'-1)}
Jac_{\zeta(B+4), \cdots, \zeta C}\pi(\psi',\epsilon',(\rho,A,B+4,\zeta,\epsilon_{0})).
$$Si on remplace ci-dessus $C$ par $C'$ et $C'$ par $C+1$ on obtient l'oppos\'e de la contribution du couple $(C,C')$ trouv\'e pr\'ec\'edemment ou encore
 le terme correspondant au couple $(C,C')$ avec $C'>C>B+1$ a une contribution oppos\'ee \`a celle du couple $C'=C+1 \leq C=C'$. En d'autre termes 
 $$
 \oplus_{C\in ]B+1,A]}Jac_{\zeta(B+1)}(\dag_{C})=0
 $$
 On applique $Jac_{\zeta (B+1)}$ \`a $(\dag'_{C})$ apr\`es avoir utilis\'e \ref{calculdujac1} pour calculer $Jac_{\zeta(B+2)}$. On suppose encore que $C\neq B+1$ et on trouve:
 $$
 Jac_{\zeta(B+1)}(\dag'_{C})=
 (-1)^{A-C}\oplus_{\eta=\pm}(-1)^{[(A-B-1)/2]}\eta^{A-B-1}\epsilon_{0}^{A-B}
<\rho\vert\,\vert^{\zeta B}, \cdots, \rho\vert\,\vert^{-\zeta C}>\times $$
$$
Jac_{\zeta(B+3), \cdots, \zeta C}\pi(\psi',\epsilon',(\rho,A,B+3,\zeta,\eta),(\rho,B,B,\zeta,\eta\epsilon_{0})).
$$
On r\'ecrit $(*)_{\eta}$ en utilisant $(\rho,A,B+1,\zeta)$ pour d\'efinir $\pi(\psi',\epsilon', \cdots)$ qui intervient dans cette expression et on somme tout de suite sur $\eta$. Cela donne:
$$
\oplus_{\eta=\pm}(-1)^{[(A-B+1)/2]}\eta^{A-B+1}\epsilon_{0}^{A-B} \oplus_{C\in ]B+1,A]}(-1)^{A-C}
<\rho\vert\,\vert^{\zeta( B+1)}, \cdots, \rho\vert\,\vert^{-\zeta C}>\times$$
$$ Jac_{\zeta (B+3), \cdots, \zeta C}\pi(\psi',\epsilon',(\rho,A,B+3,\zeta,\eta), (\rho,B,B,\eta \epsilon_{0})\eqno(10)
$$
$$
\oplus_{\eta=\pm, \eta'=\pm}(-1)^{[(A-B+1)/2]}\eta^{A-B+1}\epsilon_{0}^{A-B} (-1)^{[(A-B)/2]}(\eta')^{A-B}(\eta)^{A-B-1}$$
$$\pi(\psi',\epsilon',(\rho,A,B+2,\zeta,\eta'),(\rho,B+1,B+1,\zeta,\eta\eta'),(\rho,B,B,\eta \epsilon_{0})).\eqno(11)
$$Comme $(-1)^{[(A-B+1)/2]}=-(-1)^{[(A-B-1)/2]}$, en appliquant $Jac_{\zeta (B+1)}$ \`a (10), on annule les contributions de $Jac_{\zeta(B+1)}(\dag'_{C}$ pour tout $C\neq B+1$. On a ainsi montr\'e que
$$
Jac_{\zeta(B+1)}\pi(\psi,\epsilon)=Jac_{\zeta(B+1)}\biggl( (-1)^{A-B-1}<\rho\vert\,\vert^{\zeta B}, \cdots, \rho \vert\,\vert^{-\zeta(B+1)}> \times \pi(\psi',\epsilon',(\rho,A,B+2,\zeta,\epsilon_{0}))\eqno(12)$$
$$
\oplus_{\eta=\pm, \eta'=\pm}(-1)^{[(A-B+1)/2]}\eta^{A-B+1}\epsilon_{0}^{A-B} (-1)^{[(A-B+2)/2]}(\eta')^{A-B}(\eta)^{A-B-1}$$
$$\pi(\psi',\epsilon',(\rho,A,B+2,\zeta,\eta'),(\rho,B+1,B+1,\zeta,\eta\eta'),(\rho,B,B,\eta \epsilon_{0}))\biggr).
$$
Comme dans le cas $A=B+2$, les seuls termes associ\'ees \`a $\eta,\eta'$ qui donne une contribution non nulle quand on applique $Jac_{\zeta(B+1)}$ sont ceux pour lesquels $\eta'=\epsilon_{0}$ (cf \ref{rappel} propri\'et\'e 6) le signe est alors $(-1)^{[(A-B+1)/2+[(A-B)/2]}$; en particulier il ne d\'epend pas de $\eta$ et quand on somme sur $\eta$ le $Jac$ est $<\rho\vert\,\vert^{\zeta B}, \cdots, \rho \vert\,\vert^{-\zeta B}> \times \pi(\psi',\epsilon',(\rho,A,B+2,\zeta,\epsilon_{0}))$ au signe pr\`es que l'on vient d'\'ecrire. On remarque que $(-1)^{[(A-B+1)/2+[(A-B)/2]}$ est $(+1)$ si $(A-B)$ est pair et $-1$ sinon, c'est-\`a-dire qu'il vaut $(-1)^{A-B}$. La contribution que l'on vient de trouver est donc l'oppos\'ee de (12), ce qui donne la nullit\'e cherch\'ee.
\nl
Il reste \`a d\'emontrer que
$Jac_{x,x}\pi(\psi,\epsilon)=0$ pour tout $x\in \mathbb{R}$.
Gr\^ace \`a ce que l'on a d\'ej\`a d\'emontr\'e, le seul cas
\`a consid\'erer est $x=\zeta B$. Et clairement, par
l'hypoth\`ese de r\'ecurrence pour tout $\eta=\pm$
$$Jac_{x,x} \pi(\psi',\epsilon', (\rho,A,B+1,\zeta,\eta),
(\rho,B,B,\zeta,\eta\epsilon_{0}))=0$$. 
Soit $C\in ]B,A]$, on  pose $\delta_{C}:=<\rho\vert\,\vert^{\zeta B}, \cdots, \rho\vert\,\vert^{-\zeta C}>$. On a $$
Jac_{\zeta B}Jac_{\zeta(B+2), \cdots, \zeta C}\pi(\psi',\epsilon',(\rho,A,B+2,\zeta,\epsilon_{0})=
Jac_{\zeta(B+2), \cdots, \zeta C}Jac_{\zeta B}\pi(\psi',\epsilon',(\rho,A,B+2,\zeta\epsilon,\epsilon_{0})=0
$$
par la partie de \ref{proprietedujac} d\'ej\`a d\'emontr\'e. De plus $Jac_{\zeta B}\delta_{C}^*=0$ et $Jac_{\zeta B,\zeta B}\delta_{C}=0$. Cela suffit pour obtenir
$$
Jac_{\zeta B,\zeta B}\biggl(\delta_{C}\times Jac_{\zeta (B+2), \cdots, \zeta C}\pi(\psi',\epsilon',(\rho,A,B+2,\zeta,
\epsilon_{0}))\biggr)=0.
$$
Cela termine la preuve.
\subsection{\label{calculdujac2}}
\bf Proposition: \sl soit $(\rho,A,B,\zeta,\epsilon_{0})\in Jord(\psi,\epsilon)$ tel que $A>B$, alors $Jac_{\zeta B,\cdots,
-\zeta A}\pi(\psi,\epsilon)=
\pi(\psi',\epsilon', (\rho, A-1,B+1,\epsilon_{0}).$ En particulier ceci est 0 si $A=B+1$ et $\epsilon_{0}=-$.\rm
\nl
Soit $C\in ]B,A]$; on pose $\delta_{C}:=<\rho\vert\,\vert^{\zeta B}, \cdots, \rho\vert\,\vert^{-\zeta C}>$ et $Y_{C}:=Jac_{\zeta(B+2), \cdots, \zeta C}\pi(\psi',\epsilon',(\rho,A,B+2,\zeta,\epsilon_{0}))$. Et pour $\eta=\pm$, on pose $X_{\eta}:=\pi(\psi',\epsilon',(\rho,A,B+1,\zeta,\eta),(\rho,B,B,\zeta,\eta\epsilon_{0}))$.

On fixe $C$ et on calcule d'abord $Jac_{\zeta B, \cdots, -\zeta A} (\delta_{C}\times Y_{C})$. Pour tout $x\in [\zeta B,-\zeta A]$, $Jac_{x}\delta_{C}^*=0$ et les calculs standard de module de Jacquet (cf. par exemple \ref{standard}) entra\^{\i}nent alors que ce $Jac$ est la somme sur $\ell \in [(B+1),- A]$ de termes de la forme:
$$
Jac_{\zeta B, \cdots, \zeta \ell}\delta_{C}\times Jac_{\zeta (\ell -1), \cdots, -\zeta A}Y_{C}$$o\`u le premier $Jac$ n'intervient pas si $\ell= B+1$ et le deuxi\`eme n'intervient pas si $\ell=-A$. Pour tout $y\in [\zeta(B+2),\zeta A]$ et tout $y'\in [\zeta(\ell-1),-\zeta A]$, on a $\zeta(y-y') \geq (B+2-B)=2$. On a donc:
$$
Jac_{\zeta(\ell-1), \cdots, -\zeta A}Y_{C}=Jac_{\zeta(B+2), \cdots, \zeta A}Jac_{\zeta(\ell-1), \cdots, -\zeta A}\pi(\psi',\epsilon',(\rho,A,B+2,\zeta,\epsilon_{0})).\eqno(1)
$$Supposons que $\ell >-A$ et montrons que (1) est nul.  Gr\^ace \`a \ref{proprietedujac}, la non nullit\'e de (1) n\'ecessite qu'il existe $(\rho,A',B',\zeta')\in Jord(\psi')$ tel que $\zeta' B'=\zeta(\ell-1)$. On rappelle que $\ell-1\in [B,-A]$; le fait que $\psi\circ \Delta$ est discret et que $B'\neq B$ n\'ecessairement, entra\^{\i}ne que $\vert \ell-1 \vert \notin [B,A]$. On pose $x=\zeta(\ell -1)$ et on a s\^urement que (1) se factorise par 
$$Jac_{x, \cdots, -\zeta B}\pi(\psi',\epsilon',(\rho,A,B+2,\zeta,\epsilon_{0})).\eqno(2)$$Le lemme \ref{techniquepoursocle} que l'on d\'emontrera ci-dessous  appliqu\'e \`a $y=-\zeta B$ prouve la nullit\'e de (2) et donc celle de (1). Remarquons d'autre part que si $\ell=-A$, $Jac_{\zeta B, \cdots, -\zeta A=\zeta \ell}\delta_{C}=0$ sauf pour $C=A$ ou cela vaut 1. Ainsi la contribution des termes $\delta_{C}\times Y_{C}$ au calcul du module de Jacquet cherch\'e provient uniquement du terme correspondant \`a $C=A$ et ce terme donne $$Jac_{\zeta (B+2), \cdots, \zeta A}\pi(\psi',\epsilon', (\rho,A,B+2,\zeta,\epsilon_{0})).\eqno(3)
$$
On fixe maintenant $\eta=\pm$ et on d\'emontre plus g\'en\'eralement, 

\sl
soit $(\psi_{1},\epsilon_{1})$ un couple de m\^eme type que $\psi,\epsilon$, en particulier $\psi_{1}\circ \Delta$ est suppos\'e discret; on suppose qu'il existe $B_{1}$ un demi-entier et un signe $\zeta_{1}$ tel que $(\rho,B_{1},B_{1},\zeta_{1}) \in Jord(\psi_{1})$. Alors:
$$
Jac_{\zeta_{1}B_{1}, \cdots, -\zeta_{1}B_{1}}\pi(\psi_{1},\epsilon_{1})=0.
$$\rm
On d\'emontre cette assertion par r\'ecurrence; elle est vraie pour les s\'eries discr\`etes, elle est donc vraie pour les morphismes $\psi_{1}$ qui sont \'el\'ementaires \'etant donn\'e la fa\c{c}on que l'on a de les obtenir \`a partir des s\'eries discr\`etes. On est donc ramen\'e au cas o\`u il existe $(\rho,A',B',\zeta')\in Jord(\psi_{1})$ avec $A'>B'$. On utilise ce quadruplet pour donner la d\'efinition de $\pi(\psi_{1},\epsilon_{1})$; cette d\'efinition fait intervenir des $\pi(\psi'_{1},\epsilon'_{1}, (\rho,A',B'+2,\zeta',\epsilon'_{0}))$ et des $\pi(\psi'_{1},\epsilon'_{1},(\rho,A',B'+1,\zeta',\eta'),(\rho, B',B',\zeta',\eta'\epsilon'_{0}))$ qui v\'erifient les m\^emes hypoth\`eses que $\psi_{1},\epsilon_{1}$. Pour pouvoir conclure par r\'ecurrence, il suffit de v\'erifier que $Jac_{\zeta_{1}B_{1}, \cdots, -\zeta_{1}B_{1}}$ commute aux inductions par des repr\'esentations de la forme $<\rho\vert\,\vert^{\zeta' B'}, \cdots, \rho\vert\,\vert^{-\zeta' C'}>$ o\`u $C'\in ]B',A']$ et des op\'erations $Jac_{\zeta'(B'+2), \cdots, \zeta' A'}$. Cela a \'et\'e v\'erifi\'e dans \ref{independance}. D'o\`u l'assertion.

Ainsi $Jac_{\zeta B, \cdots , -\zeta A}\pi(\psi,\epsilon)=Jac_{\zeta (B+2), \cdots, }\pi(\psi',\epsilon', (\rho,A,B+2,\zeta,\epsilon_{0}))=\pi(\psi',\epsilon',(\rho,A-1,B+1,\zeta,\epsilon_{0}))
$ gr\^ace \`a \ref{calculdujac1} appliqu\'e \`a $(\rho,A,B+2,\zeta,\epsilon_{0})$. Il ne nous reste donc plus qu'\`a d\'emontrer le lemme ci-dessous.
\subsection{\label{techniquepoursocle}}
Le lemme ci-dessous est purement technique et $\rho,\psi,\epsilon$ sont fix\'es comme dans tout ce chapitre. En particulier $\psi\circ \Delta$ est discret.
\nl
\bf Lemme: \sl Soit $y$ un
demi-entier relatif tel que pour tout $(\rho,A',B',\zeta')\in
Jord(\psi)$,
$\vert y\vert \notin [B',A']$.
Alors pour tout
$x$ tel que
$\vert x\vert \leq \vert y\vert$ avec $x-y\in \mathbb{Z}$,
$Jac_{x,\cdots, y}\pi(\psi,\epsilon)=0$.\rm
\nl
On le d\'emontre par r\'ecurrence. On suppose d'abord que $\psi$ est \'el\'ementaire.
L'hypoth\`ese sur $y$ dit simplement que pour tout bloc de Jordan de $\psi$ n\'ecessairement de la forme $(\rho,B',B',\zeta')$ on a $\vert y \vert \neq B'$. De plus on sait que $\vert x\vert \leq
\vert y\vert$. Il faut \'evidemment revenir \`a la
d\'efinition; on peut  supposer, comme nous le
ferons que $\pi(\psi,\epsilon)$ n'est pas cuspidale. On a
alors montr\'e qu'il existe
$\zeta=\pm$ et
$\alpha\in \mathbb{N}$ avec $\alpha>1$ tels que $(\rho,
(\alpha-1)/2,(\alpha-1)/2,\zeta))\in Jord(\psi)$ et
l'une des 2 situations ci-dessous est r\'ealis\'ee (on pose $B:=(\alpha-1)/2$):

 $B>0$ et $(\rho,B-1,B-1,\zeta')\notin Jord(\psi)$ (quelque soit $\zeta'=\pm$) ou $B=1/2$ et $\epsilon(\rho,B,B,\zeta)=+$; on note alors
$\psi_1$ le morphisme qui se d\'eduit de $\psi$ en changeant
le bloc de Jordan 
$(\rho,
B,B,\zeta)$ en
$(\rho,B-1,B-1,\zeta)$ et
$\epsilon_1$ le morphisme qui se d\'eduit naturellement de
$\epsilon$ et alors $\pi(\psi,\epsilon)$ est l'unique
sous-module irr\'eductible de l'induite
$\rho\vert\,\vert^{\zeta B}\times \pi(\psi_1,\epsilon_1)$;

sinon alors $B>1/2$ et $(\rho,B-1,B-1,\zeta) \in Jord(\psi)$ avec  $\epsilon(\rho,B,B,\zeta
)= \epsilon
(\rho,B-1,B-1,\zeta)$ et on note
$\psi_1$ le morphisme qui se d\'eduit de $\psi$ en enlevant 
ces 2 blocs de Jordan et $\epsilon_1$ la restriction
\'evidente de $\epsilon$. Alors, $\pi(\psi,\epsilon)$ est
l'un des 2 sous-modules irr\'eductibles de l'induite:
$$
<\rho\vert\,\vert^{\zeta B},\cdots,
\rho\vert\,\vert^{-\zeta(B-1)}>\times
\pi(\psi_1,\epsilon_1).
$$
De plus, dans les 2 cas, $B$ est tel que
$Jac_z\pi(\psi,\epsilon)=0$ pour tout $z$ tel que $\vert
z\vert <B$. 

\nl
Revenons
\`a la propri\'et\'e que l'on cherche \`a montrer; avec ce
qui pr\'ec\`ede la nullit\'e est claire si $\vert x\vert
<B$. On suppose donc que
$\vert x\vert
\geq B$ et on raisonne par l'absurde en supposant
que $Jac_{x,\cdots, y}\pi(\psi,\epsilon)\neq 0$. Dans le
premier cas, on a soit
$Jac_{x,\cdots, y}\pi(\psi_1,\epsilon_1)\neq 0$ soit il
existe $x_1\in [x,y]$ tel que $\vert x_1\vert=B$
et 
$Jac_{x,\cdots x_1+\gamma, x_1-\gamma, \cdots
y}\pi(\psi_1,\epsilon_1)\neq 0$, o\`u $\gamma=1$ si le segment $[x,y]$ est d\'ecroissant et $-1$ sinon. Cela entra\^{\i}ne encore $Jac_{x_1-\gamma,
\cdots, y}\pi(\psi_1,\epsilon_1)\neq 0$. On a encore $\vert
x_1-\gamma\vert \leq \vert y\vert$ et toutes les
hypoth\`eses de l'\'enonc\'e v\'erifi\'ees. On obtient alors
le r\'esultat par r\'ecurrence. Le deuxi\`eme cas se traite
de la m\^eme fa\c con; on montre encore que si
$Jac_{x,\cdots, y}\pi(\psi,\epsilon)\neq 0$, il existe
$x'$ avec $\vert x'\vert \leq \vert y\vert$ tel que $[x',y]$
soit un sous-intervalle de $[x,y]$ et tel que
n\'ecessairement
$Jac_{x',\cdots y}\pi(\psi_1,\epsilon_1)\neq 0$ (cf.
\ref{notation}). Cela permet de conclure.

On fait maintenant le cas g\'en\'eral; on retraduit l'hypoth\`ese, pour tout $(\rho,\alpha)$ bloc de Jordan de $\psi\circ \Delta$, $\vert y \vert \neq (\alpha-1)/2$. Pour terminer la preuve, il faut consid\'erer le cas o\`u dans $Jord(\psi)$  se trouve $(\rho,A,B,\zeta)$ avec $A>B$.  Et on suppose que $Jac_{x,\cdots,
y}\pi(\psi,\epsilon)\neq 0$. On revient \`a la
d\'efinition de $\pi(\psi,\epsilon)$ avec les notations, pour $C\in ]B,A]$, $\delta_{C},Y_{C}$ de la preuve de \ref{calculdujac2} ainsi que $X_{\eta}$ pour $\eta=\pm$. Comme les morphismes qui interviennent dans la d\'efinition de $X_{\eta}$ on la m\^eme restriction \`a la diagonale de $SL(2,{\mathbb C})$ que $\psi$, on a tout de suite par r\'ecurrence que pour tout $\eta=\pm$, $Jac_{x,\cdots, y}X_{\eta}=0$. Fixons maintenant $C\in ]B,A]$ et on calcule
$
Jac_{x,\cdots, y}(\delta_{C}\times X_{C})$ par les formules standard: cet \'el\'ement est combinaison lin\'eaire d'\'el\'ements index\'es par un d\'ecoupage de l'intervalle $[x,y]$ en trois sous-esembles (dont certains peuvent \^etre vides, not\'es $F_{1},F_{2},F_{3}$ et tel que $Jac_{x\in F_{1}}\delta_{C}\neq 0$, $Jac_{x\in F_{2}}\delta_{C}^*\neq 0$ et $Jac_{x\in F_{3}}Y_{C}\neq 0$. Fixons un tel d\'ecoupage en supposant que le terme correspondant n'est pas 0. Montrons que $y\in F_{3}$; s'il n'en \'etait pas ainsi, on aurait soit $\vert y \vert \in [B,C]$, ce qui est exclu par hypoth\`ese, soit $\vert y\vert <B$. Si cette derni\`ere condition est satisfait tous les points du segment $[x,y]$ sont de valeur absolue strictement inf\'erieur \`a $B$ et il en est donc de m\^eme des \'el\'ements de $F_{1}$ et des \'el\'ements de $F_{2}$; or $F_{1}$ s'il n'est pas vide contient $\zeta B$ qui ne v\'erifie pas cette condition et $F_{2}$ s'il n'est pas vide contient $\zeta A$ qui ne v\'erifie pas non plus cette conditions. Ainsi, on a abouti \`a une contradiction qui prouve que $y\in F_{3}$. On a alors remarqu\'e qu'il existe $x'\in [x,y]$ tel que $Jac_{x',\cdots y}Y_{C}\neq 0$. En revenant \`a la d\'efinition, on a donc:
$$
Jac_{\zeta(B+2), \cdots, \zeta C, x', \cdots , y}\pi(\psi',\epsilon',(\rho,A,B+2,\zeta,\epsilon_{0}))\neq 0.
$$
On a soit $\vert y\vert <B$ soit $\vert y \vert >C$. Dans le premier cas, on a:
$$
Jac_{\zeta(B+2), \cdots, \zeta C, x', \cdots , y}\pi(\psi',\epsilon',(\rho,A,B+2,\zeta,\epsilon_{0}))=
Jac_{\zeta(B+2), \cdots, \zeta C}Jac_{x',\cdots y}\pi(\psi',\epsilon',(\rho,A,B+2,\zeta,\epsilon_{0}))=0
$$
par l'hypoth\`ese de r\'ecurrence. Dans le 2e cas $y$ est soit l'\'el\'ement maximal soit l'\'el\'ement minimal de l'ensemble $[\zeta (B+2), \zeta C]\cup [x',y]$ et il existe donc $x''$ dans cet ensemble tel que $[x'',y]$ soit un segment (croissant ou d\'ecroissant) et 
$$
Jac_{x'',\cdots y}\pi(\psi',\epsilon',(\rho,A,B+2,\zeta,\epsilon_{0}))\neq 0.
$$
Le point est de remarque que $\vert x''\vert <\vert y\vert$ mais cela vient de ce que $\vert y \vert$ est maximal dans l'ensemble $[ (B+2),  C]\cup [\vert x'\vert ,\vert y\vert ]$. Cela termine la preuve. 

\section{Caract\'erisation de $\pi(\psi,\epsilon)$ comme
repr\'esentation \label{caracterisationcommerepresentation}}
Dans tout ce paragraphe on suppose que $\psi\circ \Delta$ est discret. Et
le but de ce paragraphe est de d\'emontrer que les
\'el\'ements $\pi(\psi,\epsilon)$ a priori d\'efinis dans le
groupe de Grothendieck sont en fait une somme de
repr\'esentations irr\'eductibles toutes in\'equivalentes. 
On utilise la notation suivante: soit $Y$ une
repr\'esentation semi-simple et soit $\delta$ une
repr\'esentation irr\'eductible d'un groupe $GL$, on note
$<\delta,Y>$ le socle de l'induite $\delta\times Y$. On \'enonce ici les r\'esultats que l'on d\'emontrera en \ref{demonstration}.

\subsection{Description par r\'ecurrence \label{enonce}}
\bf Th\'eor\`eme: \sl L'\'el\'ement $\pi(\psi,\epsilon)$ du
groupe de Grothendieck est une repr\'esentation dont la
description se fait par r\'ecurrence ainsi; on fixe
$(\rho,a,b)\in Jord (\psi)$ tel que $inf(a,b)\geq 2$.

Dans le cas o\`u $inf(a,b)=2$ et $\epsilon(\rho,a,b)=-$ la
d\'efinition suffit $$\pi(\psi,\epsilon)=\oplus_{\eta=\pm}\pi(\psi',\epsilon',(\rho,a+\zeta_{a,b}1,b-\zeta_{a,b}1,\eta),(\rho,a-1,b-1,-\eta)).$$

 Supposons que $inf(a,b)=2$ et que
$\epsilon(\rho,a,b)=+$, alors $$\pi(\psi,\epsilon)= <<\rho\vert\,\vert^{(a-b)/2},
\cdots,
\rho\vert\,\vert^{-\zeta_{a,b}((a+b)/2-1)}>,
\pi(\psi',\epsilon')>.
$$
Supposons que $inf(a,b)$ est impair
$$
\pi(\psi,\epsilon)=<<\rho\vert\,\vert^{(a-b)/2},\cdots,
\rho \vert\,\vert^{-\zeta_{a,b}((a+b)/2-1)} >,
\pi(\psi',\epsilon',(\rho,inf(a-\zeta_{a,b}2,a),inf(b+
\zeta_{a,b}2,b),\epsilon(\rho,a,b)))>
$$
$$
\oplus
\pi(\psi',\epsilon',\cup_{\alpha\in
[\vert a-b\vert +1,a+b-1]_2}
(\rho,sup(\zeta_{a,b} \alpha,1), sup
(-\zeta_{a,b}\alpha, 1),
\epsilon(\rho,a,b)(-1)^{(\alpha-\vert a-b\vert-1)/2})).
$$
Supposons que $inf(a,b)$ est pair si $\epsilon\neq
(-1)^{inf(a,b)/2}$, alors $$\pi(\psi,\epsilon)=
<<\rho\vert\,\vert^{(a-b)/2}, \cdots,
\rho\vert\,\vert^{-\zeta_{a,b}((a+b)/2-1)}>,
\pi(\psi',\epsilon',
(\rho,inf(a,a-\zeta_{a,b}2),inf(b,b+\zeta_{a,b}2),
\epsilon(\rho,a,b)))$$si $\epsilon(\rho,a,b)=(-1)^{(inf(a,b)/2}$,
alors \`a la repr\'esentation ci-dessus, il faut ajouter
$$\oplus_{\eta=\pm}
\pi(\psi',\epsilon', \cup_{\alpha\in [\vert
a-b\vert+1,a+b-1]_2}
(\rho,sup(\zeta_{a,b}\alpha,1),sup(-\zeta_{a,b}\alpha,1),
(-1)^{(\alpha-\vert a-b\vert-1)/2}\eta))).$$
\rm
De fa\c{c}on assez concise, en utilsant les notations de \ref{autreformulation}, le th\'eor\`eme ci-dessus se r\'ecrit, 

\sl soit $(\rho,A,B,\zeta)\in Jord(\psi)$ avec $A>B$, on note comme ci-dessus $\psi',\epsilon'$ le couple qui se d\'eduit de $\psi,\epsilon$ en enlevant le quadruplet $(\rho,A,B,\zeta)$ et avec la notation $\epsilon_{0}:=\epsilon(\rho,A,B,\zeta)$:
$$
\pi(\psi,\epsilon)= <<\rho\vert\,\vert^{\zeta B}, \cdots, \rho\vert\,\vert^{-\zeta A}> \times \pi(\psi',\epsilon', (\rho,A-1,B+1,\zeta,\epsilon_{0}))>
$$
si $A=B+1$ ce terme n'appara\^{\i}t que si $\epsilon(\rho,A,B,\zeta)=+$ et vaut le socle de l'induite avec $\pi(\psi',\epsilon')$
$$
\oplus_{\eta=\pm; \epsilon_{0}=\eta^{A-B+1}(-1)^{(A-B+1)(A-B)/2}}\pi(\psi',\epsilon',\cup_{C\in [B,A]}(\rho,C,C,\zeta,(-1)^C\eta)).
$$
La somme dans ce dernier terme contient un terme si $A-B$ est pair, il en contient 0 si $A-B$ est impair avec $\epsilon_{0}\neq (-1)^{(A-B+1)/2}$ et 2 dans le cas restant. On appelle cette somme les termes compl\'ementaires.

\subsection{Description en tant que
repr\'esentations\label{explicite}} Dans cette section on
tire quelques cons\'equences du th\'eor\`eme permettant de
mieux cerner la repr\'esentation
$\pi(\psi,\epsilon)$; on pourra appliquer ces propri\'et\'es
par r\'ecurrence. On utilise le langage de \ref{autreformulation} les formules \'etant plus simples bien qu'\`a mon avis moins conceptuelles.
\nl 
{\bf Th\'eor\`eme}: {\sl
Soit $(\rho,A,B,\zeta,\epsilon_{0})\in Jord(\psi,\epsilon)$ tel que $A>B$.
Alors
$
\pi(\psi,\epsilon)=$

$\oplus_{\ell \in [0,[(A-B)/2]],
\eta=\pm;\epsilon_{0}=\eta^{A-B+1}\prod_{C\in [B+\ell,A-\ell] }(-1)^{[C]}
}$
$$
<<\rho\vert\,\vert^{\zeta B},\cdots,\rho\vert\,\vert^{-\zeta
A)}>, \cdots,
<\rho\vert\,\vert^{\zeta(B+\ell)},
\cdots, \rho\vert\,\vert^{-\zeta (A-\ell)})>,
\pi(\psi',\epsilon', \cup_{C\in
[B+\ell, A-\ell]}
(\rho,C,C,\zeta,\eta
(-1)^{[C]})))>;
$$
il n'y a pas d'induction pour $\ell=0$ et pour $\ell=[(A-B+1)/2]$
avec $A-B$ impair il n'y a pas de $\cup$ .}
\nl
\rm
On remarque que si $A-B+1$ est impair la somme sur
$\eta$ contient exactement un terme quelque soit $\ell$, par
contre si $A-B+1$ est pair, la somme est, suivant les
valeurs de $\ell$ pour $\epsilon_{0}$  fix\'e,
soit vide soit contient 2 termes (il y a alternance) .

C'est une application r\'ep\'et\'ee du th\'eor\`eme. 
\subsection{Sans multiplicit\'e\label{sansmultiplicite}}
\bf Proposition: \sl la repr\'esentation $\pi(\psi,\epsilon)$ est sans multiplicit\'e toujours sous l'hypoth\`ese que $\psi\circ \Delta$ est discret.\rm
\nl
La proposition est trivialement vraie si $\psi$ est \'el\'ementaire car alors $\pi(\psi,\epsilon)$ est irr\'eductible. Avec la r\'ecurrence que l'on a mise en place ici, nous allons montrer que cette proposition est cons\'equence du th\'eor\`eme \ref{enonce}. Pour cela, nous allons reformuler ce th\'eor\`eme. On reprend les notations de \ref{autreformulation} et on fixe $(\rho,A,B,\zeta)\in Jord(\psi)$ en supposant que $A>B$. On note $Irr(\pi(\psi,\epsilon)$ l'ensemble des repr\'esentations irr\'eductibles constituant $\pi(\psi,\epsilon)$ compt\'ees avec multiplicit\'e et on g\'en\'eralise cette notation \`a tous les couples $(\psi'',\epsilon'')$ de m\^eme nature que $\psi,\epsilon$. Pour $\eta=\pm$, on note $\psi_{\eta},\epsilon_{\eta}$ le couple qui se d\'eduit de $\psi,\epsilon$ en enlevant $(\rho,A,B,\zeta)$ et en ajoutant les quadruplets $(\rho,C,C,\zeta)$ pour $C\in [B,A]$ en prolongeant $\epsilon$ par $\epsilon_{\eta}(\rho,C,C,\zeta)=(-1)^C \eta$ pour tout $C\in [B,A]$ et ce qui n'est d\'efini, par hypoth\`ese, que si $\epsilon_{0}=(\prod_{C\in [B,A]}(-1)^{[C]})\eta^{(A-B+1}$. Avec cette notation, le th\'eor\`eme de \ref{enonce} dit que $Irr(\pi(\psi,\epsilon))$ est l'union  des $Irr(\pi(\psi_{\eta},\epsilon_{\eta}))$ (pour les $\eta$ convenables) avec l'ensemble des sous-modules irr\'eductibles de l'induite $<\rho\vert\,\vert^{\zeta B}, \cdots, \rho\vert\,\vert^{-\zeta A}>\times \pi(\psi',\epsilon',(\rho,A-1,B+1,\zeta,\epsilon(\rho,A,B,\zeta))$.

Montrons d'abord plus pr\'ecis\'ement que 

\nl\sl
l'application $Jac_{\zeta B, \cdots, -\zeta A}$ induit une bijection de $Irr(\pi(\psi,\epsilon))-\cup_{\eta=\pm}Irr(\pi(\psi_{\eta},\epsilon_{\eta}))$ sur $Irr(
\pi(\psi',\epsilon',(\rho,A-1,B+1,\zeta,\epsilon(\rho,A,B,\zeta)))$. 
Et que de plus $Jac_{\zeta B, \cdots, -\zeta A}$ annule $\cup_{\eta=\pm}Irr(\pi(\psi_{\eta},\epsilon_{\eta}))$.\rm
\nl

On a montr\'e dans la preuve  de \ref{calculdujac2} que $Jac_{\zeta B, \cdots, -\zeta B}$ annule $\cup_{\eta=\pm}Irr(\pi(\psi_{\eta},\epsilon_{\eta}))$; il en est donc de m\^eme \`a fortiori de $Jac_{\zeta B, \cdots, -\zeta A}$. D'o\`u la deuxi\`eme propri\'et\'e de l'assertion.

Montrons la bijection annonc\'ee; par d\'efinition du socle pour tout $X\in Irr(\pi(\psi,\epsilon))-\cup_{\eta=\pm}Irr(\pi(\psi_{\eta},\epsilon_{\eta}))$, 
$
Jac_{\zeta B, \cdots, -\zeta A}X\neq 0.
$
Et m\^eme il existe un \'el\'ement $Y\in Irr(\pi(\psi',\epsilon',(\rho,A-1,B+1,\zeta,\epsilon(\rho,A,B,\zeta)))$ tel que $$X\hookrightarrow <\rho\vert\,\vert^{\zeta B}, \cdots, \rho\vert\,\vert^{-\zeta A}>\times Y.\eqno(1)$$
On note $\sigma_{Y}$ l'induite du membre de droite de (1). On va d\'emontrer que $Jac_{\zeta B, \cdots, -\zeta A}\sigma_{Y}=Y$ ce qui suffira, car alors $\sigma$ n'a qu'un unique sous-module irr\'eductible n\'ecessairement $X$ et l'application qui \`a $Y$ associe l'unique sous-module de $\sigma_{Y}$ est l'inverse de la bijection annonc\'ee.

Les formules standard pour calculer les modules de Jacquet montrent que si l'assertion cherch\'ee n'est pas vraie, il existe $x\in [\zeta B,-\zeta A]$ tel que $Jac_{x,\cdots, -\zeta A}Y\neq 0$. On aurait alors aussi $$Jac_{x, \cdots, -\zeta A}\pi(\psi',\epsilon',(\rho,A-1,B+1,\zeta,\epsilon(\rho,A,B,\zeta))\neq 0,
$$
ce qui est exclu par \ref{techniquepoursocle} puisque $A\notin [B',A']$ pour tout $(\rho,A',B',\zeta')\in Jord(\psi')$ par l'hypoth\`ese $\psi\circ \Delta$ discret.
\nl
L'hypoth\`ese de r\'ecurrence assure qu'un $X$ de $Irr(\pi(\psi,\epsilon))-\cup_{\eta}Irr(\pi(\psi_{\eta},\epsilon_{\eta}))$ intervient sans multiplicit\'e dans $Irr(\pi(\psi,\epsilon))$. On sait aussi que chaque ensemble $Irr(\pi(\psi_{\eta},\epsilon_{\eta})$ est sans multiplicit\'e mais il faut encore d\'emontrer quand 2 valeurs de $\eta$ sont possibles que $Irr(\pi(\psi_{+},\epsilon_{+})) \cap Irr(\pi(\psi_{-},\epsilon_{-}))=\emptyset.$
\nl
Montrons l'assertion suivante par r\'ecurrence; soit $\psi$ comme ci-dessus et $D$ un entier positif ou nul tel qu'il existe un signe $\zeta''$ avec $(\rho,C,C,\zeta'')\in Jord(\psi)$ et soit $\epsilon_{i}$ pour $i=1,2$ des applications de $Jord(\psi)$ dans $\{\pm 1\}$ telles que $\epsilon_{1}(\rho,D,D,\zeta'')\neq \epsilon_{2}(\rho,D,D,\zeta'')$. Alors 
$$Irr(\pi(\psi,\epsilon_{1}))\cap Irr(\pi(\psi,\epsilon_{2}))=\emptyset.$$
On d\'emontre cela par r\'ecurrence, si $\psi$ est \'el\'ementaire cela r\'esulte de ce que l'application qui \`a $\epsilon$ associe $\pi(\psi,\epsilon)$ est injective et \`a valeurs dans l'ensemble des repr\'esentations irr\'eductibles. Si $\psi$ n'est pas \'el\'ementaire, on utilise \ref{enonce} comme on l'a fait ci-dessus en remarquant que l'hypoth\`ese sur $(\rho,D,D,\zeta'')$ perdure. Cela termine la preuve de la proposition.

\subsection{Un cas de pr\'esentation par r\'ecurrence  simple\label{lecastroue}}
On garde les notations, $\psi',\epsilon'$ qui se d\'eduit de
$\psi,\epsilon$ en enlevant le bloc de Jordan fix\'e
$(\rho,A,B,\zeta,\epsilon_{0})$. On suppose que $B\geq 1$ et que $(\rho,A-1,B-1,\zeta) \notin Jord(\psi)$ et on note ${\psi}_{-}$ le morphisme qui se d\'eduit de $\psi,\epsilon$ en rempla\c{c}ant le bloc $(\rho,A,B,\zeta,\epsilon_{0})$ par $(\rho,A-1,B-1,\zeta,\epsilon_{0})$. On suppose que ${\psi}_{-}\circ \Delta$ est encore discret. Cela se traduit aussi par le fait que pour tout quadruplet $(\rho,A',B',\zeta')\in Jord(\psi)$, on a $A'\neq B-1$.
\nl
\bf Proposition: \sl Soit $\psi,\epsilon$ et $(\rho,A,B,\zeta)\in
Jord(\psi)$ v\'erifiant les hypoth\`eses ci-dessus. Alors 
$$
\pi(\psi,\epsilon)=<<\rho\vert\,\vert^{\zeta B}, \cdots,
\rho\vert\,\vert^{\zeta A}>,
\pi(\psi',\epsilon',
(\rho,A-1,B-1,\zeta,\epsilon_{0}))>.$$
\nl
\rm
On d\'emontre cette proposition est cons\'equence de \ref{enonce}. Et on d\'emontre cela par r\'ecurrence pour pouvoir l'utiliser dans les raisonnements par r\'ecurrence. 
On va en fait d\'emontrer quelque chose de plus pr\'ecis que l'\'enonc\'e. On note $Irr(\pi(\psi,\epsilon))$ l'ensemble des sous-modules irr\'eductibles de $\pi(\psi,\epsilon)$; on garde la notation $\psi_{-},\epsilon_{-}$ pour  le couple qui se d\'eduit de $\psi,\epsilon$ en changeant $(\rho,A,B,\zeta,\epsilon_{0})$ en $(\rho,A-1,B-1,\epsilon_{0})$; on d\'efinit donc de fa\c{c}on analogue $Irr(\pi(\psi_{-},\epsilon_{-})$. On va d\'emontrer que

\nl\sl
l'application $Jac_{\zeta B, \cdots, \zeta A}$ induit une bijection entre $Irr(\pi(\psi,\epsilon))$ et $Irr(\pi(\psi_{-},\epsilon_{-}))$, l'application inverse \'etant $Y \mapsto <<\rho\vert\,\vert^{\zeta B}, \cdots, \rho\vert\,\vert^{\zeta A}>\times Y>.$\rm
\nl
Dans ces conditions, ce qu'il faut d\'emontrer et que pour tout $X\in Irr(\pi(\psi,\epsilon))$, $Jac_{\zeta B, \cdots \zeta A}X\neq 0$ et que pour tout $Y\in Irr(\pi(\psi_{-},\epsilon_{-})$, $Jac_{\zeta B, \cdots, \zeta A}
(< \rho\vert\,\vert^{\zeta B}, \cdots, \rho\vert\,\vert^{\zeta A}>\times Y)=Y.$

Montrons la premi\`ere partie, c'est-\`a-dire soit $X\in Irr(\pi(\psi,\epsilon))$ et montrons que $Jac_{\zeta B, \cdots, \zeta A}X\neq 0$. S'il existe $(\rho',A',B',\zeta')\in Jord(\psi)- \{(\rho,A,B,\zeta)\}$ avec $A'>B'$, on d\'emontre le r\'esultat par r\'ecurrence en utilisant \ref{enonce} pour $(\rho',A',B',\zeta')$ au lieu de $(\rho,A,B,\zeta)$. On suppose donc que pour tout $(\rho,A',B',\zeta')$ diff\'erent de $(\rho,A,B,\zeta)$, on a $A'=B'$. Ainsi si $A=B$, on est ramen\'e au cas \'el\'ementaire et le r\'esultat r\'esulte des constructions (cf. \ref{rappel} propri\'et\'es). On suppose donc que $A>B$. Soit $X\in Irr(\pi(\psi,\epsilon))$. On suppose d'abord que $X$ intervient dans le terme compl\'ementaire de \ref{enonce}; on est alors ramen\'e au cas de morphismes \'el\'ementaires et l'on a, pour toute collection de signe $\epsilon_{C}$ index\'ee par $C\in [B,A]$:
$$
\pi(\psi',\epsilon', \cup_{C\in [B,A]} (\rho,C,C,\zeta,\epsilon_{C}) -\hookrightarrow \rho\vert\,\vert^{\zeta B} \times \cdots \times \rho \vert\,\vert^{\zeta A} \times \pi(\psi',\epsilon', \cup_{C\in [B,A]} (\rho,C-1,C-1,\zeta,\epsilon_{C})).$$
Et $X$ est le terme de gauche pour un bon choix de $\epsilon_{?}$. D'o\`u la non nullit\'e annonc\'ee. On consid\`ere finalement le cas o\`u $X$ est un sous-module irr\'eductible de l'induite
$$
<\rho\vert\,\vert^{\zeta B}, \cdots, \rho\vert\,\vert^{-\zeta A}> \times \pi(\psi',\epsilon',(\rho, A-1,B+1,\zeta,\epsilon_{0})).
$$On fixe $X'$ un sous-module irr\'eductible de $\pi(\psi',\epsilon',(\rho, A-1,B+1,\zeta,\epsilon_{0}))$ telle que l'inclusion ci-dessus se factorise par $X'$ et on applique la proposition \`a $X'$ en utilisant $(\rho, A-1,B+1,\zeta)$. On note alors $X''$ tel que $Jac_{\zeta (B+1), \cdots, \zeta (A-1)}X'=X''$ d'o\`u aussi $$X'= <<\rho\vert\,\vert^{\zeta(B-1)}, \cdots, \rho\vert\,\vert^{\zeta(A-1)}>, X''>.
$$
Avec des r\'esultats standard, on sait que
$$
<\rho\vert\,\vert^{\zeta B}, \cdots, \rho\vert\,\vert^{-\zeta A}> \hookrightarrow 
\rho\vert\,\vert^{\zeta B} \times <\rho\vert\,\vert^{\zeta(B-1)}, \cdots, \rho\vert\,\vert^{-\zeta(A-1)}> \times \rho \vert\,\vert^{-\zeta A}.
$$
On en d\'eduit
$$
X\hookrightarrow \rho\vert\,\vert^{\zeta B} \times <\rho\vert\,\vert^{\zeta(B-1)}, \cdots, \rho\vert\,\vert^{-\zeta(A-1)}> \times \rho \vert\,\vert^{-\zeta A }\times <\rho\vert\,\vert^{\zeta(B-1)}, \cdots, \rho\vert\,\vert^{\zeta(A-1)}>\times X'' \simeq
$$
$$
\rho\vert\,\vert^{\zeta B} \times <\rho\vert\,\vert^{\zeta(B-1)}, \cdots, \rho\vert\,\vert^{\zeta(A-1)}>\times 
<\rho\vert\,\vert^{\zeta(B-1)}, \cdots, \rho\vert\,\vert^{-\zeta(A-1)}> \times
 \rho \vert\,\vert^{-\zeta A }\times X''.
 $$
 Le point est de d\'emontrer que l'hypoth\`ese entra\^{\i}ne que l'induite $ \rho \vert\,\vert^{-\zeta A }\times X''$ est irr\'eductible. Cela est \'equivalent \`a d\'emontrer que l'induite $\rho\vert\,\vert^{-\zeta A}\times X''$ est isomorphe \`a l'induite $\rho\vert\,\vert^{\zeta A}\times X''$ (chacune des induites \`a un unique sous-module irr\'eductible car $Jac_{\pm A}X''=0$). On le sait dans le cas o\`u $X''$ est de la forme $\pi(?,?)$ o\`u $(?,?)$ est un morphisme \'el\'ementaire car $(\rho,2A-1)$ n'est pas un bloc de Jordan au sens ordinaire de ce morphisme, par hypoth\`ese. Et on se ram\`ene \`a ce cas, en utilisant \ref{explicite}.  Ayant cela, on revient aux inclusions ci-dessus:
 $$
 X\hookrightarrow \rho\vert\,\vert^{\zeta B}\times <\rho\vert\,\vert^{\zeta(B-1)}, \cdots, \rho\vert\,\vert^{\zeta(A-1)}>\times <\rho\vert\,\vert^{\zeta(B-1)}, \cdots, \rho\vert\,\vert^{-\zeta(A-1)}> \times
 \rho \vert\,\vert^{-\zeta A }\times X''\simeq
 $$
 $$
 \rho\vert\,\vert^{\zeta B}\times <\rho\vert\,\vert^{\zeta(B-1)}, \cdots, \rho\vert\,\vert^{\zeta(A-1)}>\times <\rho\vert\,\vert^{\zeta(B-1)}, \cdots, \rho\vert\,\vert^{-\zeta(A-1)}> \times
 \rho \vert\,\vert^{\zeta A }\times X''\simeq$$
 $$
  \rho\vert\,\vert^{\zeta B}\times <\rho\vert\,\vert^{\zeta(B-1)}, \cdots, \rho\vert\,\vert^{\zeta(A-1)}>\times 
 \rho \vert\,\vert^{\zeta A }\times <\rho\vert\,\vert^{\zeta(B-1)}, \cdots, \rho\vert\,\vert^{-\zeta(A-1)}> \times
 X''.
 $$
 Et cela donne $Jac_{\zeta B, \cdots, \zeta A}X\neq 0$ comme annonc\'e. Comme on sait que $Jac_{x}X=0$ pour tout $x\in ]\zeta B,\zeta A]$, il existe donc $Y$ une repr\'esentation irr\'eductible et une inclusion $X\hookrightarrow <\rho\vert\,\vert^{\zeta B}, \cdots, \rho\vert\,\vert^{\zeta A}> \times Y$; comme en plus $Jac_{\zeta B,\zeta B}X=0$, on a certainement  $Jac_{x}Y=0$ pour tout $x\in [\zeta B,\zeta A]$ d'o\`u ais\'ement $Jac_{\zeta B, \cdots, \zeta A}X=Y$. Ceci termine la preuve.

\nl
Nous utiliserons cette proposition en particulier dans le
cas suivant: $\psi,\epsilon$ est fix\'e avec $(\rho,A,B,\zeta,\epsilon_{0})\in
Jord(\psi,\epsilon)$ tel que $A>B+1$. On remarque que dans le corollaire le groupe qui intervient est de rang strictement plus petit que $G$. On peut donc appliquer le corollaire sans probl\`eme dans les d\'emonstrations par r\'ecurrence comme nous les faisons ici.
\nl
\bf Corollaire: \sl Pour tout $C\in ]B+1,A]$, $Jac_{\zeta (B+2), \cdots,\zeta C }
\pi(\psi',\epsilon', (\rho,A,B+2,\zeta,\epsilon_0))=$
$$
<<\rho\vert\,\vert^{\zeta(C+1)}, \cdots,
\rho\vert\,\vert^{\zeta A}>,\pi(\psi',\epsilon',
(\rho,A-1,B+1,\zeta,\epsilon_0))>.
$$
\rm
La proposition pr\'ec\'edente s'applique pour donner:
$$
\pi(\psi',\epsilon',(\rho,A,B+2,\zeta,\epsilon_0))=
<<\rho\vert\,\vert^{\zeta(B+2)},\cdots,
\rho\vert\,\vert^{\zeta A}>,
\pi(\psi',\epsilon',(\rho,A-1,B+1,\zeta,\epsilon_0))>.
$$
Et le corollaire s'en d\'eduit par un calcul de module de
Jacquet.
\subsection{Propri\'et\'e d'injectivite \label{injectivite}}
Pour $i=1,2$ fixons des couples $(\psi_{i},\epsilon_{i})$ et supposons que $\psi_{1}\circ \Delta=\psi_{2} \circ \Delta=: \psi_{\Delta}$. Dans ces conditions $Cent_{^LG}\psi_{i}$ pour $i=1,2$ sont des sous-groupes de $Cent _{^LG}(\psi_{\Delta})$; cela a donc un sens de consid\'erer leur intersection. On suppose encore ici que $\psi\circ \Delta$ est discret et ici de fa\c{c}on d\'eterminante on suppose que l'application qui \`a un caract\`ere $\epsilon$ du groupe $Cent_{^LG}(\psi \circ \Delta)$ associe la s\'erie discr\`ete $\pi(\psi\circ \Delta,\epsilon)$ est injective. On a aussi besoin de la m\^ eme assertion pour tous les morphismes $\phi'$ de $W_{F}\times SL(2,{\mathbb C})$ dans le $L$ groupe d'un groupe de   m\^eme type que $G$ tel que $Jord(\phi')\subset Jord(\psi\circ \Delta)$.
\nl
\bf Proposition: \sl avec les notations et hypoth\`eses ci-dessus, supposons que la restriction de $\epsilon_{1}$ \`a $Cent_{^LG}(\psi_{1})\cap Cent_{^LG}(\psi_{2})$ ne co\"{\i}ncide pas avec la restriction de $\epsilon_{2}$ \`a ce m\^eme groupe. Alors $Irr(\pi(\psi_{1},\epsilon_{1}))\cap Irr (\pi(\psi_{2},\epsilon_{2}))=\emptyset.$\rm

\nl
Avant de d\'emontrer cette proposition, on va exprimer en termes combinatoires le fait que $\psi_{1}\circ \Delta=\psi_{2}\circ \Delta$; avec les notations de \ref{autreformulation} cela dit exactement que pour tout $\rho$,  on a l'\'egalit\'e ensembliste $\cup_{A_{1},B_{1}}[B_{1},A_{1}]=\cup_{A_{2},B_{2}} [B_{2},A_{2}]$ o\`u $A_{i},B_{i}$ pour $i=1,2$ parcourent l'ensemble des couples tels qu'il existe un signe $\zeta$ v\'erifiant $(\rho,A_{i},B_{i},\zeta)\in Jord(\psi_{i})$.
\nl
On d\'emontre cette proposition en admettant \ref{enonce} par r\'ecurrence d'abord sur le rang du groupe puis  sur $\sum_{(\rho,A,B,\zeta)\in Jord(\psi_{1} \cup Jord(\psi_{2}} A-B$. Si ce nombre est nul, les 2 morphismes $\psi_{1}$ et $\psi_{2}$ sont \'el\'ementaires mais ils ne co\"{\i}ncident pas forc\'ement. Toutefois, les centralisateurs de $\psi_{1}$, de $\psi_{2}$ et de $\psi_{\Delta}$ co\"{\i}ncident. On peut donc dire que $\epsilon_{1}\neq \epsilon_{2}$.   On a rappel\'e en \ref{rappel}, la d\'efinition des nombres $a_{\rho,\psi_{i},\epsilon_{i}}, b_{\rho,\psi_{i},\epsilon_{i}}$ et des signes $\zeta_{\rho,\psi_{i},\epsilon_{i}}$, pour $i=1,2$. Ces d\'efinitions se voient sur les modules de Jacquet, donc si il n'y a pas identit\'e entre les valeurs pour $i=1$ et celles pour $i=2$, il est imm\'ediat que les repr\'esentations correspondantes sont diff\'erentes. Pour simplifier l'\'ecriture, on note ces nombres, $a,b$ et le signe $\zeta$. Si $a>b+2$, $Jac_{\zeta (a-1)/2}\pi_{\psi_{i},\epsilon_{i}}$ est de la forme $\pi(\psi'_{i},\epsilon'_{i})$ o\`u l'on a simplement remplac\'e le bloc $(\rho,(a-1)/2,(a-1)/2,\zeta, \epsilon_{i}(a))$ par $(\rho,(a-3)/2,(a-3)/2,\zeta,\epsilon_{i}(a)$. On obtient alors le r\'esultat par r\'ecurrence. Si $a=b+2$, on d\'efinit $\psi'_{i}, \epsilon'_{i}$ en enlevant les 2 blocs $(\rho,(a-1)/2,(a-1)/2,\zeta)$ et $(\rho,(b-1)/2,(b-1)/2,\zeta)$. 

On remarque que l'hypoth\`ese $\epsilon_{1}\neq \epsilon_{2}$ entra\^{\i}ne que soit $\epsilon'_{1}\neq \epsilon'_{2}$ soit que $a=3,b=1$ et $\epsilon_{1}(\rho,3)\neq \epsilon_{2}(\rho,3)$; en effet dans le cas o\`u $a$ est pair la valeur de $\epsilon_{i}$ sur $(\rho,a)$ est pr\'ecis\'ement $(-1)^{a/2 +1}$ et est donc ind\'ependante de $i$. Si $a$ est impair, on a n\'ecessairement $b\geq 1$ avec nos hypoth\`eses et $\epsilon_{i}(\rho,a)=\epsilon_{i}(\rho,1)(-1)^{(a-1)/2}$. Ainsi $\epsilon_{i}(\rho,a)$ n'est d\'etermin\'e par $\epsilon'_{i}$ que si $(\rho,1)\in Jord(\psi'_{i})$ c'est-\`a-dire $b>1$.

Si $\epsilon'_{1}\neq \epsilon'_{2}$ on obtient le r\'esultat par r\'ecurrence en remarquant que $Jac_{\zeta(a-1)/2, \cdots, -\zeta(b-1)/2}\pi(\psi_{i},\epsilon_{i})=\pi(\psi'_{i},\epsilon'_{i})$. Dans le cas restant, on suppose pour avoir quelque chose \`a d\'emontrer que $\pi(\psi'_{1},\epsilon'_{1})=\pi(\psi'_{2},\epsilon'_{2})$. On commence par v\'erifier que cela entra\^{\i}ne que $\psi'_{1}=\psi'_{2}$; on sait que pour tout $(\rho,\alpha)\in Jord(\psi_{\Delta})$ avec $\alpha>3$, il existe pour $i=1,2$ un signe $\zeta_{\alpha,i}$ tel que $(\rho,\alpha,\zeta_{\alpha,i})\in Jord(\psi'_{i})$ et le point est de d\'emontrer que $\zeta_{\alpha,1}=\zeta_{\alpha,2}$ pour tout $\alpha$. Cela se fait avec les modules de Jacquet, consid\'erons d'abord $\alpha_{0}$ l'\'el\'ement minimal; on a, pour $i=1,2$:
$$
Jac_{\zeta_{\alpha_{0},i}(\alpha_{0}-1)/2}\pi(\psi'_{i},\epsilon'_{i})=\pi(\psi''_{i},\epsilon''_{i}),
$$
o\`u $\psi''_{i}$ s'obtient simplement en rempla\c{c}ant $(\rho,\alpha_{0})$ par $(\rho, \alpha_{0}-2)$ sans rien changer aux signes et 
$$
Jac_{-\zeta_{\alpha_{0},i}(\alpha_{0}-1)/2}\pi(\psi'_{i},\epsilon'_{i})=0.
$$
Cela entra\^{\i}ne d\'ej\`a $\zeta_{\alpha_{0},1}=\zeta_{\alpha_{0},2}$; puis on recommence avec l'\'el\'ement minimal de $Jord_{\rho}(\psi''_{i})-\{\alpha_{0}-2\}$ et avec $\pi(\psi''_{i},\epsilon''_{i})$ qui donne le m\^eme type de relations.

Maintenant que l'on sait que $\psi'_{1}=\psi'_{2}$, on a aussi $\psi_{1}=\psi_{2}$ et on conna\^{\i}t l'injectivit\'e de l'application $\epsilon\mapsto \pi(\psi,\epsilon)$. Cela termine la preuve du d\'ebut de la r\'ecurrence.

On suppose donc que soit $\psi_{1}$ soit $\psi_{2}$ n'est plus \'el\'ementaire et on fixe $B$ un entier minimum avec la propri\'et\'e qu'il existe un entier $A$ et un signe $\zeta$ tel que $A>B$ et  $(\rho,A,B,\zeta) \in Jord(\psi_{i})$ pour une valeur de $i$ au moins. Si $i=1$ et $i=2$ sont possibles alors $A$ est maximal pour le choix de $i$ fait et par sym\'etrie on peut supposer que $i=1$. Supposons qu'il existe $X\in Irr(\psi(\psi_{i},\epsilon_{i}))$ pour $i=1,2$. On applique \ref{enonce} \`a $\psi_{1},\epsilon_{1}$ en utilisant $(\rho,A,B,\zeta)$. Si $Jac_{\zeta B, \cdots, -\zeta A}X=0$, on sait qu'il existe $\eta=\pm$ (avec les notations d\'ej\`a introduites) tel que $X\in Irr(\pi(\psi_{\eta},\epsilon_{\eta})$. On v\'erifie que $Cent_{^LG}(\psi_{\eta})$ contient $Cent_{^LG}(\psi_{1})$ et on v\'erifie que la restriction de $\epsilon_{\eta}$ \`a $Cent_{^LG}(\psi_{1})$ n'est autre que $\epsilon_{1}$ (c'est m\^eme comme cela qu'on a construit $\eta$). On obtient le r\'esultat cherch\'e par r\'ecurrence, en rempla\c{c}ant $\psi_{1},\epsilon_{1}$ par $\psi_{\eta},\epsilon_{\eta}$. On suppose donc que 
$$
Jac_{\zeta B, \cdots, -\zeta A}X:=Y \neq 0.
$$En particulier $Jac_{\zeta B}X \neq 0$.
D'apr\`es \ref{proprietedujac} appliqu\'e \`a $\psi_{2}$, on sait qu'il existe $A_{2}$ un entier tel que $(\rho,A_{2},B,\zeta)\in Jord(\psi_{2})$. Il n'y a aucune raison pour que $A_{2}=A$. On pose maintenant $A=A_{1}$. On v\'erifie comme en \ref{calculdujac2} que la non nullit\'e de $Jac_{\zeta B, \cdots -\zeta B}X$ entra\^{\i}ne que $A_{2}>B$.  D'apr\`es \ref{enonce} pr\'ecis\'e en ****, on sait que pour $i=1,2$:
$$
Jac_{\zeta B, \cdots ,-\zeta A_{i}}X=:X_{1,i} \in Irr(\pi(\psi'_{i},\epsilon'_{i}, (\rho, A_{i}-1,B_{i}+1,\zeta,\epsilon_{i}(\rho,A_{i},B)).
$$Si $A_{1}=A_{2}$, on obtient le r\'esultat par r\'ecurrence puisque les hypoth\`eses sont encore v\'erifi\'ees pour les morphismes intervenant dans le terme de droite ci-dessus.

Supposons donc que $A_{1}\neq A_{2}$. On note $A_{max}$ le plus grand entier tel que pour tout $C \in [B,A_{max}]$ $(\rho,2C+1) \in Jord(\psi_{\Delta})$. On a n\'ecessairement $A_{max}\geq A$; le calcul de $Jord(\psi_{\Delta})$ en fonction de $Jord(\psi_{i})$ et la maximalit\'e de $A_{max}$, assurent que pour $i=1,2$, il existe un d\'ecoupage en sous-intervalles $[B_{i}^k, A_{i}^k]$ o\`u $k\in [1,N_{i}]$ avec $N_{i}$ un entier tel que $[B_{i}^1,A_{i}^1]=[B,A_{i}]$ et pour tout $k\in [1,N_{i}]$, il existe un signe $\zeta_{k,i}$ tel que $(\rho,A_{i}^k,B_{i}^k,\zeta_{k,i})\in Jord(\psi_{i})$. En particulier $\zeta_{1,i}=\zeta$ pour $i=1,2$. On note encore $A_{0}$ le plus petit entier tel qu'il existe pour $i=1,2$ un entier $k_{i}\leq N_{i}$ tel que $A_{0}=A_{k_{1}}=A_{k_{2}}$; comme ces \'egalit\'es sont v\'erifiees pour $A_{max}$, un tel entier existe. Toujours pour $i=1,2$, on noter $\psi^0_{i},\epsilon^0_{i}$ le couple qui se d\'eduit de $\psi_{i},\epsilon_{i}$ en enlevant les blocs $(\rho,A_{i}^j,B_{i}^j,\zeta_{j,i})$ de $Jord(\psi_{i})$. 

On peut appliquer la proposition \ref{lecastroue} \`a $\pi(\psi'_{i},\epsilon'_{i}, (\rho, A_{i}-1,B_{i}+1,\zeta,\epsilon_{i}(\rho,A_{i},B)$ en utilisant $(\rho, A_{i}^2,B_{i}^2,\zeta_{2,i})$ pour calculer $Jac_{\zeta_{2,i}B_{i}^2, \cdots, \zeta_{2,i}A_{i}^2} X_{1,i}$, puis de proche en proche pour calculer:
$$
Jac_{\zeta_{2,i}B_{i}^2}, \cdots, \zeta_{2,i}A_{i}^2, \cdots, \zeta_{k_{i},i}B_{i}^{k_{i}}, \cdots, \zeta_{k_{i},i}A_{i}^{k_{i}}X_{1,i}=Y_{i}
$$
$$\in Irr(\pi(\psi_{i}^{0}, \epsilon_{i}^0, (\rho,A_{i}-1,B_{i}+1,\epsilon_{i}(\rho,A_{i},B_{i}\zeta)), \cup_{j\in [2,k_{i}}](\rho, A_{i}^j-1,B_{i}^j-1,\zeta_{j,i},\epsilon_{i}(\rho,A_{i}^j,B_{i}^j,\zeta_{j,i})). \eqno(1)
$$
Montrons que pour $i=1,2$ et pour tout $j\in [2,k_{i}]$, $\zeta_{j,i}=-\zeta$. En effet, s'il n'en est pas ainsi, on fixe $j$ minimum avec cette propri\'et\'e, d'o\`u $\zeta_{j,i}=\zeta$ et les propri\'et\'es standard de \ref{notation} (1), montrent que la non nullit\'e ci-dessus entra\^{\i}ne a fortiori la non nullit\'e de $Jac_{\zeta_{j,i} B_{i}^j}X_{1,i}$ et en revenant encore \`a la d\'efinition de $X_{1,i}$ la non nullit\'e de $Jac_{\zeta B_{i}^j}X$ (ici on utilise le fait que $B_{i}^j>A_{i}>B_{i}$, d'o\`u $B_{i}^j-B_{i}>1$). Ainsi, en posant $i':= \{1,2\}-\{i\}$, il existe $A'$ tel que $(\rho,B_{i}^j,A',\zeta) \in Jord(\psi_{i'})$ d'apr\`es \ref{proprietedujac}. Ainsi il existe $j'\leq k_{i'}$ tel que $B_{i}^j=B_{i'}^{j'}$. D'o\`u encore:
$$
A_{i}^{j-1}=B_{i}^j- 1=B_{i'}^{j'}-1=A_{i'}^{j'-1}<A^0,$$
ce qui contredit la minimalit\'e de $A^0$. Mais alors ainsi $Y_{i}$ ne d\'epend pas de $i$, on le note $Y$ On remarque que  $\epsilon_{i}$ pour $i=1,2$ sont uniquement d\'etermin\'es par le caract\`ere intervenant dans le membre de droite de (1) donc ces caract\`eres ne co\"{\i}ncident pas sur l'intersection des commutants et on obtient une contradiction par la r\'ecurrence appliqu\'ee \`a $Y$.
Cela termine la preuve.

\section{Preuve de la caract\'erisation comme
repr\'esentation\label{demonstration}}
\subsection{Quelques notations\label{notations}}
On appelle constituant irr\'eductible de
$\pi(\psi,\epsilon)$ une repr\'esentation irr\'eductible qui
intervient avec un coefficient non nul dans l'\'el\'ement
$\pi(\psi,\epsilon)$ vu dans le groupe de Grothendieck.
On fixe $(\rho,A,B,\zeta,\epsilon_{0})\in Jord(\psi,\epsilon)$  avec
$A>B$. N\'ecessairement $\tilde{\pi}$ est un sous-quotient irr\'eductible de l'une
des repr\'esentations suivantes:
$$\exists C\in ]B,A],\qquad
X_C:=\delta_{C}\times Y_C$$
o\`u $\delta_{C}=<\rho\vert\,\vert^{\zeta B}, \cdots,
\rho\vert\,\vert^{-\zeta A}>$ et  $Y_C:=Jac_{\zeta (B+2), \cdots,
\zeta C}\pi(\psi',\epsilon',(\rho,A,B+2,\zeta, \epsilon_{0}))$ o\`u il n'y a pas de $Jac$ pour $C=B+1$
soit
$$\exists \eta\in \{\pm 1\}, \qquad
X_\eta:=\pi(\psi',\epsilon',(\rho,A,B+1,\zeta,\eta),(\rho,B,B,\zeta,\eta\epsilon_{0}
)).$$
On d\'emontre le th\'eor\`eme en
l'admettant 
 par r\'ecurrence pour les repr\'esentations
$\pi(\psi'',\epsilon'')$ telles que soit le rang du groupe correspondant est plus petit que le rang de $G$ soit on a \'egalit\'e et $\sum_{(\rho'',A'',B'',\zeta'')\in Jord(\psi'')} A''-B''<
\sum_{(\rho',A',B',\zeta')\in Jord(\psi)}A'-B'$.

 On remarque que par \ref{explicite} les
repr\'esentations
$X_C$ et $X_\eta$ sont connues par r\'ecurrence.

\subsection{Propri\'et\'es des modules de
Jacquet\label{jacquet}}
\nl
\bf Proposition : \sl Soit $x\in \mathbb{R}$ et soit
$\tilde{\pi}$ un constituant irr\'eductible de
$\pi(\psi,\epsilon)$ alors $Jac_x\tilde{\pi}=0$ sauf s'il
existe $(\rho,A,B,\zeta) \in Jord_\rho(\psi)$ tel que $x=\zeta B$.
De plus $Jac_{x,x}\tilde{\pi}=0$ pour tout $x$\rm
\nl
La proposition est vraie si $\psi$ est \'el\'ementaire au
sens de \ref{rappel} (cf propri\'et\'es de loc.cit.). On la d\'emontre par r\'ecurrence.  On
reprend les notations de
\ref{notations} et comme dans la preuve de
\ref{moduledejacquet}, on v\'erifie que les seules valeurs de
$x$ qui ne satisfont pas directement \`a la conclusion de la
proposition sont, pour le $\rho$ de
\ref{notations}
$$
x=\zeta C, C\in [B,A].
$$Consid\'erons d'abord le cas o\`u $x=\zeta B$; il faut ici
d\'emontrer que $Jac_{\zeta B,\zeta B}\tilde{\pi}=0$. Avec
\ref{notations} il suffit de d\'emontrer cela pour $X_C$ et
pour $X_\eta$ (avec les notations de loc.cit.); pour
$X_\eta$ cela se fait par r\'ecurrence et pour $X_C$ cela a
\'et\'e fait dans la preuve de \ref{moduledejacquet}.

Soit $x=\zeta C$ pour $C\in ]B,A]$, en particulier $x\neq 0$ et on montre
d'abord que:
\nl
\sl pour
tout constituant $\tilde{\pi}$ de $\pi(\psi,\epsilon)$, $Jac_{x,x}\tilde{\pi}=0$ et  soit
$Jac_x\tilde{\pi}=0$ soit $Jac_x\tilde{\pi}$ est
irr\'eductible et dans ce dernier cas, l'induite
$\rho\vert\,\vert^x\times \tilde{\pi}$ a un unique
sous-module irr\'eductible qui est alors $\tilde{\pi}$.\rm
\nl
En effet soit  $\tilde{\pi}$ un constituant irr\'eductible de
$\pi(\psi,\epsilon)$ et montrons que
$Jac_{x,x}\tilde{\pi}=0$; pour cela, on fixe, comme en
\ref{notations}, soit
$C\in ]B,A]$ soit
$\eta$ tel que $\tilde{\pi}$ soit un sous-quotient
irr\'eductible de $X_C$ ou $X_\eta$. Et il suffit donc de
d\'emontrer que $Jac_{x,x}X_C=0$ pour tout $C$ comme ci-dessus et
$Jac_{x,x}X_\eta=0$ pour $\eta=\pm$. On va le faire en
d\'etail pour $x=\zeta(B+1)$ qui est le cas le plus
difficile. Fixons d'abord $\eta=\pm$, on a \'ecrit en
\ref{enonce} une description explicite des constituants de
$X_\eta$: l'un des termes (non irr\'eductible en g\'en\'eral) est $<\sigma,Y>$ o\`u
$\sigma=<\rho\vert\,\vert^{\zeta(B+1)}, \cdots,
\rho\vert\,\vert^{-\zeta A}>$ et $$Y=\pi(\psi',\epsilon',
(\rho, A-1,B+2,\eta), (\rho,B,B,\zeta,
\eta\epsilon_{0})).$$
Les autres termes sont de la forme $\pi(\psi',\epsilon',\cup_{D\in ]B,A]}
(\rho,D,D,\zeta,(-1)^{[D]}\lambda),
(\rho,B,B,\zeta,\eta\epsilon_{0}))$, o\`u $\lambda$ est un signe convenable. On montre que pour
ces derniers termes $Jac_{\zeta (B+1), \zeta(B+1)}=0$ par exemple par r\'ecurrence. Pour le terme $<\sigma,Y>$, on a par la
proposition appliqu\'ee par r\'ecurrence
$Jac_{\zeta (B+1)}Y=0$.
Or on a aussi $Jac_{\zeta (B+1)}\sigma^*=0$ pour un tel
$\sigma$ et il suffit donc de remarquer que
$Jac_{\zeta(B+1),\zeta(B+1)}\sigma=0$ ce qui r\'esulte des
r\'esultats de Zelevinski. Cela prouve notre assertion pour les constitutants de  $X_{\eta}$.

Revenons donc \`a $C$ mais pour le moment, on suppose que 
$C\in ]B+1,A]$; ici on a $Jac_{\zeta(B+1)}
\delta_C=Jac_{\zeta(B+1)}\delta^*_C=0$; c'est l'hypoth\`ese
$C>B+1$. On sait \ref{lecastroue} corollaire que
$$Y_C=<<\rho\vert\,\vert^{\zeta (C+1)},\cdots,
\rho\vert\,\vert^{A}>,
\pi(\psi',\epsilon',(\rho,A-1,B+1,\zeta,\epsilon_{0}))>.$$
Donc l'hypoth\`ese $j>1$ assure que
$Jac_{\zeta(B+1),\zeta(B+1)}Y_C$ se calcule en fonction de
$Jac_{\zeta(B+1),\zeta(B+1)}
\pi(\psi',\epsilon',(\rho,A-1,B+1,\epsilon_{0}))$ qui
est nul par r\'ecurrence.

Reste le cas $C=B+1$, ici on
v\'erifie que
$Jac_{\zeta(B+1)}Y_{C=B+1}=0$: en effet $Y_{B+1}=
\pi(\psi',\epsilon',(\rho,A,B+2,\epsilon(\rho,a,b)))$.
D'o\`u notre assertion en utilisant l'hypoth\`ese de
r\'ecurrence. Ensuite il n'y a plus qu'\`a v\'erifier que
$Jac_{\zeta(B+1)} <\rho\vert\,\vert^{B},\cdots,
\rho\vert\,\vert^{-\zeta (B+1)}>=0$ et
$
Jac_{\zeta(B+1),\zeta(B+1)}<\rho\vert\,\vert^{\zeta B},\cdots,
\rho\vert\,\vert^{-\zeta(B+1)}>^*=0.$
Cela termine la preuve de l'assertion $Jac_{x,x}\tilde(\pi)=0$.

Soit maintenant $\tilde{\pi}$ un constituant irr\'eductible
de $\pi(\psi,\epsilon)$ tel que $Jac_{x}\tilde{\pi}\neq 0$.
Pour un bon choix de repr\'esentation irr\'eductible $X$, on a une inclusion, par r\'eciprocit\'e de
Frobenius:
$
\tilde{\pi}\hookrightarrow \rho\vert\,\vert^x\times X
$et puisque $Jac_{x,x}\tilde{\pi}=0$, n\'ecessairement $Jac_{x}X=0$.
On calcule $Jac_x$ de l'induite de droite et on trouve que
cela vaut $X$; ici on a utilis\'e la non nullit\'e de $x$. On a ainsi montr\'e toutes nos assertions.
\nl
On a d\'emontr\'e en \ref{proprietedujac} pour les valeurs de $x$ \'ecrites
ci-dessus, $Jac_x\pi(\psi,\epsilon)=0$ et on vient de montrer que pour
$\tilde{\pi}$ et
$\tilde{\pi}'$ des constituants irr\'eductibles de
$\pi(\psi,\epsilon)$ s'il y a une intersection (dans
le groupe de Grothendieck) entre
$Jac_x \tilde{\pi}$ et $Jac_x \tilde{\pi}'$ alors ces repr\'esentations qui sont irr\'eductibles co\"{\i}ncident et que $\tilde{\pi}\simeq \tilde{\pi}'$ en tant qu'unique sous-module irr\'eductible de l'induite $\rho\vert\,\vert^x \times Jac_{x}\tilde{\pi}$. Il est donc clair qu'il ne peut pas y avoir de simplifications et  qu'aucun constituant irr\'eductible de $\pi(\psi,\epsilon)$ ne v\'erifie $Jac_x\neq 0$ pour ces valeurs de $x$. Cela termine la preuve.

\subsection{D\'ebut de la preuve}
\subsubsection{Hypoth\`eses pour cette sous-section\label
{hypotheses}}

 La d\'emonstration se fait par r\'ecurrence puisqu'il n'y a rien \`a d\'emontrer dans le cas o\`u pour
tout $(\rho,A,B,\zeta)\in Jord(\psi)$, $A=B$. Elle est
fastidieuse et technique, le lecteur comprendra certainement
que ce que l'on cherche \`a d\'emontrer n'est pas limpide.
Je n'ai pas de r\'ealisation de cette repr\'esentation.
On va commencer par traiter des cas particuliers auxquels la
preuve par r\'ecurrence se ram\`enera. De cette fa\c con
cela sera un peu moins abstrait. On fixe $(\rho,A,B,\zeta,\epsilon_{0})\in Jord(\psi,\epsilon)$ tel
que $A>B$.

Dans toute cette sous-section, on suppose
que $(\rho,A,B,\zeta,\epsilon_{0})$ que nous venons de fixer est le seul triple
de $Jord(\psi)$ v\'erifiant $A>B$. On suppose aussi
que l'ensemble des demi-entiers $D$ pour lesquels il existe un signe $\zeta_{D}$ tels que $D<B$ et $(\rho,D,D,\zeta_{D})$ est soit vide soit est un intervalle commen\c{c}ant \`a $0$ si $B$ est entier et \`a 1/2 sinon; dans ce dernier cas, on suppose que $\epsilon(\rho,1/2,1/2,\zeta_{1/2})=-1$. On suppose aussi que si l'ensemble ci-dessus est un intervalle non nul, $\epsilon$ alterne sur cet intervalle.

 A partir de \ref{1construction} on va supposer que $\zeta=+$ et que $A,B$ sont des entiers; dans \ref{lecasa=b+1}, on ne fait pas ces hypoth\`eses et le lecteur verra ainsi ce que l'on \'evite en les faisant. L'hypoth\`ese sur $\zeta$ permet de fixer la croissance des segments utilises pour la classification des repr\'esentations irr\'eductibles de $GL$. L'hypoth\`ese sur $A,B$ \'evite de devoir pr\'eciser si on commence \`a $0$ ou \`a $1/2$; ici les 2 cas ne sont pas identiques, le cas o\`u $A$ et $B$ sont entiers est plus difficile que l'autre car les repr\'esentations que l'on construits ne sont pas totalement d\'etermin\'ees par leur module de Jacquet uniquement dans ce cas l\`a. En particulier \ref{dernierereduction} n'a pas d'objet si $A,B$ sont des demi-entiers non entiers.

\subsubsection{Le cas le plus simple $A=B+1$ \label{lecasa=b+1}}
 Le cas o\`u $A=B+1$ a un \'enonc\'e
 tr\`es simple rappelons le avant de le d\'emontrer:

\sl
 si
$\epsilon_{0}=+$ $\pi(\psi,\epsilon)$ est  l'unique sous-module
irr\'eductible de $<\rho\vert\,\vert^{(a/2-1}, \cdots,
\rho\vert\,\vert^{-(a/2)}>\times \pi(\psi',\epsilon')$ et en particulier est irr\'eductible. 

si $\epsilon_{0}=-$, $\pi(\psi,\epsilon)=$ 
$\oplus_{\eta=\pm}\pi(\psi',\epsilon',(\rho,B+1,B+1,\zeta,\eta), (\rho,B,B,\zeta,-\eta))$ et en particulier est de longueur 2 exactement.
\rm

Montrons cela; le cas o\`u $\epsilon_{0}=-$ est juste
la d\'efinition. Supposons donc que $\epsilon_{0}=+$.
Par d\'efinition on a:
$$
\pi(\psi,\epsilon)=<\rho\vert\,\vert^{\zeta B}, \cdots,
\rho\vert\,\vert^{-\zeta(B+1))}>\times \pi(\psi',\epsilon')\eqno(1)
$$
$$
\ominus_{\eta=\pm}\pi(\psi',\epsilon',(\rho,B+1,B+1,\zeta,\eta),
(\rho,
B-1,B-1,-\eta)).\eqno(2)
$$
On remarque que notre $B+1$ (resp. $B$) est le $(a_{\rho,\psi'',\epsilon''}-1)/2$ (resp. $(b_{\rho,\psi'',\epsilon''}-1)/2$) de \ref{rappel} pour $\psi'',\epsilon''$ se d\'eduisant de $\psi$ en enlevant le bloc $(\rho,A,B,\zeta,\epsilon_{0})$ et en y ajoutant les 2 blocs $(\rho,B+1,B+1,\zeta,\lambda)$, $(\rho,B,B,\zeta,\lambda)$ pour n'importe quel choix de $\lambda=\pm$. On obtient donc imm\'ediatement les assertions:
$<\rho\vert\,\vert^{\zeta B},\cdots,\vert\,\vert^{-\zeta B}>
\times \pi(\psi',\epsilon')$ est semi-simple de longueur 2
exactement et que (2) est la somme des 2 sous-modules
irr\'eductibles de l'induite $$\rho\vert\,\vert^{\zeta (B+1)}\times
<\rho\vert\,\vert^{\zeta B},\cdots,\vert\,\vert^{-\zeta B}>
\times \pi(\psi',\epsilon').
$$
Il est alors facile de voir que (2) est inclus dans
l'induite de (1) et ici il est d\'ej\`a clair que
$\pi(\psi,\epsilon)$ est une repr\'esentation. On applique
\ref{rappel} (propri\'et\'e 4) \`a $\pi(\psi',\epsilon')$ puisque que le $(a_{\rho,\psi',\epsilon'}-1)/2$ de loc.cit. est ici strictement sup\'erieur \`a $B+1$: soit $\tilde{\pi}$ un sous-quotient de l'induite
de (1), alors il existe un ensemble ${\cal E}$ totalement
ordonn\'e de demi-entiers relatifs tel que $$\{\vert x\vert;
x\in {\cal E}\}= \{\vert x\vert; x\in [B+1,-B]\}$$
et une repr\'esentation irr\'eductible $\sigma$ sous-module
de $\times_{x\in {\cal E}}\rho\vert\,\vert^x$
(repr\'esentation d'un GL convenable)  avec une inclusion
$\tilde{\pi} \hookrightarrow \sigma\times
\pi(\psi',\epsilon')$. On remarque que pour $x$ tel que
$Jac_x \sigma\neq 0$ n\'ecessairement $Jac_x\tilde{\pi}\neq
0$. Mais on a vu en \ref{jacquet} que cela entra\^{\i}ne
que
$x=\zeta B$; on n'a pas exactement vu cela, on pourrait aussi avoir $x=\pm D$ avec $D$ comme ci-dessus mais il faudrait alors, en repartant de la d\'efinition de $\pi(\psi,\epsilon)$,  $Jac_{x}\pi(\psi',\epsilon')\neq 0$ ce qui n'est pas possible \`a cause des hypoth\`eses sur $\epsilon$. On utilise la
classification de Zelevinski \`a l'aide de segments
d\'ecroissants si $\zeta=+$ et croissants si $\zeta=-$ et on \'ecrit  $\sigma$ dans cette
classification; c'est-\`a-dire qu'il existe $v\in
\mathbb{N}$ et pour tout $i\in [1,v]$ des segments
$[x_i,y_i]$ avec la propri\'et\'e de croissance pr\'ecis\'ee, c'est-\`a-dire $\zeta x_{i} \geq \zeta y_{i}$ et pour tout $i\in
[1,v[$  $\zeta x_i\leq \zeta x_{i+1}$. On a certainement $Jac_x\sigma\neq
0$ pour $x=x_1$ et pour $x=x_i$ pour $i$ v\'erifiant pour
tout $j<i$ $x_j\neq x_i-1$. Mais si ${\cal E}$ ne contient
pas $\zeta (B+1)$, $B$ est certainement le plus grand
\'el\'ement de $\zeta {\cal E}$ et alors n\'ecessairement $v=1$
et $[x_1,y_1]=[\zeta B, -\zeta(B+1)]$ ce qui est ce que l'on cherche.
Si ${\cal E}$ contient $\zeta(B+1)$ alors on a   $$x_1=\zeta B, x_{2}=\zeta (B+1), v=2.$$ On a 2 possibilit\'es pour $y_{1},y_{2}$. Soit:
$$y_{1}=\zeta B
y_2=-\zeta (B-1)\eqno(3)
$$ou
$$
y_1=-\zeta B; y_{2}=\zeta(B+1).\eqno(4)
$$
Mais dans le cas (3), $\sigma$ est une induite
irr\'eductible et on a aussi $Jac_{\zeta(B+1)}\sigma\neq 0$ ce
qui est exclu. Reste le cas (4). Ici on utilise le fait que
$\rho\vert\,\vert^{\zeta(B+1)}\times \pi(\psi',\epsilon')$ est
irr\'eductible donc isomorphe \`a
$\rho\vert\,\vert^{-\zeta(B+1)}\times \pi(\psi',\epsilon')$ et on
aurait une inclusion
$$
\tilde{\pi}\hookrightarrow <\rho\vert\,\vert^{\zeta B},\cdots,
\rho\vert\,\vert^{-\zeta B}>\times
\rho\vert\,\vert^{\zeta(B+1)}\times \pi(\psi',\epsilon')$$
$$
\simeq 
<\rho\vert\,\vert^{\zeta B},\cdots,
\rho\vert\,\vert^{-\zeta B}>\times
\rho\vert\,\vert^{-\zeta(B+1)}\times \pi(\psi',\epsilon').
$$
Or l'induite $<\rho\vert\,\vert^{\zeta B},\cdots,
\rho\vert\,\vert^{-\zeta B}>\times
\rho\vert\,\vert^{-\zeta(B+1)} $ est de longueur 2 exactement avec
un quotient qui v\'erifie $Jac_{-\zeta(B+1)}\neq 0$. L'inclusion
ci-dessus se factorise donc par le sous-module, i.e. $<\rho\vert\,\vert^{\zeta B}, \cdots
,\rho\vert\,\vert^{-\zeta(B+1)}>$. Pour terminer la preuve il ne
reste plus qu'\`a remarquer que l'induite:
$$<\rho\vert\,\vert^{\zeta B}, \cdots
,\rho\vert\,\vert^{-\zeta(B+1)}>\times \pi(\psi',\epsilon')$$n'a
qu'un unique sous-module irr\'eductible par r\'eciprocit\'e
de Frobenius.

\subsubsection{1e construction \label{1construction}}

On rappelle que l'on fixe $(\rho,A,B,\zeta) \in Jord(\psi)$ avec $A>B+1$ et on pose $\epsilon_{0}$ la valeur de $\epsilon$ sur cet \'el\'ement de $Jord(\psi)$. On reprend la notation $\psi',\epsilon'$ qui se d\'eduit de $\psi$ en enlevant ce bloc de Jordan.  Montrons
\nl
\bf Lemme: \sl Soit $\tilde{\pi}$ un constituant de $\pi(\psi,\epsilon)$; alors il existe un entier $t_{1}\in [0, [(A-B)/2]]$ et un signe $\lambda$ tel que $\tilde {\pi}$ soit un sous-quotient d'une induite de la forme:
$$
\times_{j\in [0,2t_{1}-1]} \rho\vert\,\vert^{A-2j} \times \cdots \times \rho\vert\,\vert^{-(A-2j-1)}\times 
\pi(\psi',\epsilon', \cup_{C\in [B, A-2t_{1}]}(\rho,C,C,\zeta,(-1)^{[C]}\lambda)),$$
o\`u $[0,2t_{1}-1]$ est un ensemble vide si $t_{1}=0$ et $[B, A-2t_{1}-2]$ est un ensemble vide si $t_{1}=(A-B)/2$. De plus n\'ecessairement $\epsilon_{0}= \prod_{C\in [B,A-2t_{1}]}\biggl((-1)^{[C]}\lambda\biggr)$. \rm
\nl
On prend la d\'efinition de $\pi(\psi,\epsilon)$ soit il existe $\eta=\pm$ tel que $\tilde{\pi}$ soit un sous-quotient de $\pi(\psi',\epsilon',(\rho,A,B+1,\zeta,\eta),(\rho,B,B,\zeta,\eta\epsilon_{0}))$ soit il existe $j\in [1,A-B]$ tel que $\tilde{\pi}$ soit un sous-quotient de l'induite (on utilise tout de suite le corollaire de \ref{lecastroue}):
$$
\rho\vert\,\vert^{\zeta B}\times \cdots \times \rho\vert\,\vert^{-\zeta (B+j)}\times \rho\vert\,\vert^{\zeta(B+j+1)}    \cdots \times \rho\vert\,\vert^{\zeta A}\times \pi(\psi',\epsilon',(\rho,A-1,B+1,\zeta,\epsilon_{0}).
$$
Consid\'erons d'abord ce deuxi\`eme cas; il est encore vrai que $\tilde{\pi}$ est alors sous-quotient de l'induite:
$$
\rho\vert\,\vert^{A}\times \cdots \times \rho\vert\,\vert^{-B}\times \pi(\psi',\epsilon',(\rho,A-1,B+1,\zeta,\epsilon_{0})).
$$
On applique encore \ref{lecastroue} pour remplacer $\pi(\psi',\epsilon',(\rho,A-1,B+1,\zeta,\epsilon_{0})$ par l'induite $\rho\vert\,\vert^{B}\times \cdots \rho\vert\,\vert^{A-1} \times 
\pi(\psi',\epsilon',(\rho,A-2,B,\zeta,\epsilon_{0})$. Et on obtient que 
$\tilde{\pi}$ est alors sous-quotient de l'induite:
$$
\rho\vert\,\vert^{A}\times \cdots \times \rho\vert\,\vert^{-(A-1)}\times \pi(\psi',\epsilon',(\rho,A-2,B,\zeta,\epsilon_{0})).
$$
On applique ensuite le lemme par r\'ecurrence pour obtenir le r\'esultat cherch\'e. Reste le premier cas; on fixe donc $\eta$ tel que $\tilde {\pi}$ est un sous-quotient de $\pi(\psi',\epsilon',(\rho,A,B+1,\zeta,\eta),(\rho,B,B,\zeta,\eta\epsilon_{0}))$.
On applique d'abord le lemme par r\'ecurrence en utilisant $(\rho,A,B+1,\zeta,\eta)$ (cette \'etape est vide si $A=B+1$). Ainsi il existe $\lambda$ et $t'$ tel que 
$\tilde {\pi}$ soit un sous-quotient d'une induite de la forme:
$$
\times_{j\in [0,2t'-1]} \rho\vert\,\vert^{A-2j} \times \cdots \times \rho\vert\,\vert^{-(A-2j-1)}\times 
\pi(\psi',\epsilon', \cup_{C\in ]B, A-2t']}(\rho,C,C,\zeta,(-1)^{[C]}\lambda), (\rho,B,B,\zeta,\eta\epsilon_{0})),$$avec la relation
$$
\prod_{C\in ]B,A-2t]}((-1)^{[C]}\lambda)=\eta.
$$
Si $(-1)^{B+1}\lambda \neq \eta \epsilon_{0},
$ ou si $B+1=A-2t'$,
le lemme s'en d\'eduit automatiquement en posant tout simplement $t_{1}=t'$. Supposons donc que 
$(-1)^{B+1}\lambda = \eta \epsilon_{0}.
$ On d\'emontre ici que $$\pi(\psi',\epsilon', \cup_{C\in ]B, A-2t']}(\rho,C,C,\zeta,(-1)^{[C]}\lambda), (\rho,B,B,\zeta,\eta\epsilon_{0}))\hookrightarrow $$
$$
<\rho\vert\,\vert^{\zeta(B+1)}, \cdots, \rho\vert\,\vert^{-\zeta B}> \times 
\pi(\psi',\epsilon', \cup_{C\in [B+2, A-2t']}(\rho,C,C,\zeta,(-1)^{[C]}\lambda)).$$
Gr\^ace au lemme \ref{lecastroue}, on peut remplacer $\pi(\psi',\epsilon', \cup_{C\in [B+2, A-2t']}(\rho,C,C,\zeta,(-1)^C\lambda))$ par l'induite $$\rho\vert\,\vert^{\zeta(B+2)} \times \cdots \times \rho\vert\,\vert^{\zeta (A-2t')}\times \pi(\psi',\epsilon', \cup_{C\in [B+1,A-1]}(\rho,C,C,\zeta,(-1)^{{[C]}+1}\lambda))
$$ puis encore $\pi(\psi',\epsilon', \cup_{C\in [B+1,A-2t'-1]}(\rho,C,C,\zeta,(-1)^{[C]+1}\lambda))$ par l'induite $$\rho\vert\,\vert^{\zeta(B+1)}\times \cdots  \times \rho\vert\,\vert^{\zeta (A-2t'-1)}\times \pi(\psi',\epsilon', 
 \cup_{C\in [B, A-2t'-2]}(\rho,C,C,\zeta,(-1)^{[C]}\lambda)).
$$On remarque que l'on a $\prod_{C\in [B,A-2t'-2]}((-1)^{[C]}\lambda)=\epsilon_{0}$, donc $
\pi(\psi',\epsilon', 
 \cup_{C\in [B, A-2t'-2]}(\rho,C,C,\zeta,(-1)^{[C]}\lambda))$ a bien la forme qui convient. Il reste \`a modifier les premiers facteurs de l'induite.
Puisque l'on cherche \`a d\'emontrer que $\tilde{\pi}$ est sous-quotient d'une certaine induite, on peut remplacer l'induite $$<\rho\vert\,\vert^{\zeta (B+1)}, \cdots, \rho\vert\,\vert^{-\zeta B}> \times \rho \vert\,\vert^{\zeta(B+2)} \times \rho^{\zeta (A-2t')}\times \rho\vert\,\vert^{\zeta (B+1)}\times \cdots \times \rho
\vert\,\vert^{\zeta (A-2t'-1)}$$
par l'induite $$\rho\vert\,\vert^{A-2t'}\times \cdots \times \rho\vert\,\vert^{-(A-2t'-1)}.$$
On obtient le lemme en posant $t_{1}=t'+1$. Cela termine la preuve.

\subsubsection{2e construction \label{2construction}}
 On fixe encore $\psi,\epsilon$ et $(\rho,a,b)\in
Jord(\psi)$ en supposant que $inf(a,b)\geq 2$ et que
$(\rho,a,b)$ est le seul triple de $Jord(\psi)$ ayant cette
propri\'et\'e; pour simplifier l'\'ecriture on suppose que $a+b$ est pair. Avec les notations de \ref{autreformulation}, cela est \'equivalent \`a fixer un quadruplet $(\rho,A,B,\zeta)\in Jord(\psi)$ avec $A>B$ et \`a supposer que c'est le seul quadruplet avec cette propri\'et\'e; la parit\'e se traduit par le fait que $A$ et $B$ sont des entiers. On note ici $\tilde{\psi}, \tilde{\epsilon}$ le couple qui se d\'eduit de $\psi,\epsilon$ en enlevant $(\rho,A,B,\zeta)$ et tous les quadruplets de la forme $(\rho,C,C,\zeta')$ pour $C< B$. Pour $D$ un demi-entier tel que $A-D\in {\mathbb N}$ et pour $\lambda$ un signe, on note $\tilde{\psi}_{D,\lambda},\tilde{\epsilon}_{D,\lambda}$ le couple qui se d\'eduit de $\tilde{\psi},\tilde{\epsilon}$ en ajoutant cette fois, les quadruplets $(\rho,C,C,\zeta)$ pour $D-C\in {\mathbb N}$  et $\tilde{\epsilon}$ vaut sur un tel quadruplet $(-1)^C\lambda$; on suppose que $\prod_{C\leq D}((-1)^{[C]}\lambda)=\epsilon_{0}$ o\`u $\epsilon_{0}=\epsilon(\rho,A,B,\zeta)$. Pour unifier les \'enonc\'es on pose $\tilde{\psi}_{-1,\lambda},\tilde{\epsilon}_{-1,\lambda}=\tilde{\psi},\tilde{\epsilon}$, o\`u $\lambda$ ne joue aucun r\^ole.  On impose toujours \`a $\lambda$ de v\'erifier:
$$\prod_{C\leq D}((-1)^{[C]}\lambda)=\epsilon_{0}\prod_{(\rho,E,E,\zeta'); E<B}\epsilon(\rho,E,E,\zeta').
$$Ci-dessous, on a des entiers \`a cause de l'hypoth\`ese de parit\'e, en g\'en\'eral on a des demi-entiers.

\bf Lemme: \sl Soit $\tilde{\pi}$ un sous-quotient
irr\'eductible de $\pi(\psi,\epsilon)$. Il existe un ensemble
d'entiers ${\cal E}$ totalement ordonn\'es , un entier
$D$ et un signe $\lambda$ comme ci-dessus v\'erifiant:
$$
\tilde{\pi} \hookrightarrow \times_{x\in {\cal
E}}\rho\vert\,\vert^x\times \pi(\tilde{\psi}_D,\tilde{
\epsilon}_{D,\lambda}),\eqno(1)
$$
$$
{\cal E}\cup \{-{\cal E}\}= \cup_{E\in [B,A]}[-E,E]  \cup _{E<B; (\rho,E,E,\zeta')\in Jord(\psi)} [-E,E]\ominus_{C\leq D} [-C,C]. \eqno(2)
$$
\rm
On montre d'abord qu'il existe un ensemble d'entiers
${\cal E}'$,  v\'erifiant (2) ainsi que $D$ et $\lambda$
comme dans l'\'enonc\'e mais tel que (1) soit affaibli en
$\tilde{\pi}$ est un sous-quotient de l'induite
\'ecrite en (1). L'\'etape faite en
\ref{1construction} ram\`ene le cas de
$(\psi,\epsilon)$ au cas d'un param\`etre
\'el\'ementaire. On peut alors lui appliquer progressivement
les d\'efinitions pour conclure. Il reste \`a appliquer
\ref{rappel} (propri\'et\'e 4) qui oblige \`a changer
\'eventuellement certains signes dans
${\cal E}'$ ainsi que son ordre pour obtenir l'inclusion
annonc\'e en (1). Les op\'erations faites sur ${\cal E}'$ ne
change pas la propri\'et\'e (2). La condition sur le signe est automatique car on ne change pas la valeur de $\prod_{(\rho,A',B',\zeta')}\epsilon(\rho,A',B',\zeta')$ dans toutes ces op\'erations.

\subsection{Un cas particulier important \label{casparticulierprincipal}} 

\subsubsection{1e \'etape \label{1etape}}On fixe $(\rho,A,B,\zeta)$ comme ci-dessus. On fait les hypoth\`eses de \ref{2construction} et on en ajoute m\^eme d'autres. Pour simplifier les notations on suppose que $A,B$ sont des entiers. Si $B=0$, on n'ajoute aucune hypoth\`eses. Par contre si $B>0$, on suppose que $(\rho,C,C,\zeta)\in Jord(\psi)$ pour tout $ 0\leq C<B$ et que $\epsilon$ alterne sur tous ces blocs, c'est-\`a-dire qu'il existe un signe $\lambda'$ tel que $\epsilon(\rho,C,C,\zeta)=(-1)^C\lambda$ pour tout $C$ comme ci-dessus.

On fixe $\tilde{\pi}$ un constituant irr\'eductible de $\pi(\psi,\epsilon)$ et on reprend d'abord la 1e
construction \ref{1construction}, d'o\`u $t_1$ et $\lambda_{1}$.
Ensuite on fait la 2e construction, d'o\`u $D$ et $\lambda$. Mais ici on peut pr\'eciser un peu les rapports entre $t_{1}$ et $D$ (on rappelle que $D$ peut valoir $-1$). On sait que
le signe $\epsilon_{t_{1},\lambda_{1}}$ vaut $\epsilon$ sur les blocs jusqu'\`a $B-1$ et donc alterne par hypoth\`ese. Puis, par construction il alterne  de
nouveau sur les $A-B-2t_1+1$ blocs compris entre
$B$ et  $A-2t_{1 }$. On 
pose $2t:=A-D$ et $2t$ est exactement le
nombre de blocs qui disparaissent de $\psi\circ \Delta$  pour arriver \`a $\tilde{\psi}_{D}$ dans ces
constructions. Donc  soit le signe $\epsilon_{t_{1},\lambda_{1}}$ alterne d\'ej\`a spontan\'ement \`a la fin de la premi\`ere construction et  $t=t_1$, soit il n'y a pas alternance et 
$t=t_1+inf(B, A-B-2t_{1}+1)$. Ci-dessous on passe en revue les cas o\`u l'on peut conclure et ensuite on montrera que ce sont les seuls cas qui peuvent se produire.

\nl
\bf Lemme: \sl avec les notations ci-dessus, supposons que
l'on peut choisir ${\cal E}$ contenant $-\zeta A$  alors
$Jac_{\zeta B, \cdots, -\zeta A}\tilde{\pi}\neq 0$.\rm

\nl
On fixe ${\cal E}$ tel que \ref{2construction} soit satisfait
et  on suppose que 
${\cal E}$ contient
$-\zeta A$. Pour simplifier la d\'emonstration on suppose que $\zeta=+$ pour ne plus devoir l'\'ecrire. Alors il existe $x\in {\cal E}$ tel que $[x,
-A]$ soit un segment et $Jac_{x,\cdots,
-A}\tilde{\pi}\neq 0$ (cf  \ref{notation}). Comme $\vert x\vert \leq
A$, n\'ecessairement $x=B$. D'o\`u le lemme
\nl
Supposons maintenant que $D=A$, c'est-\`a-dire que
$t_1=0$ dans la premi\`ere construction \ref{1construction}
et qu'il n'y a aussi aucune 2e construction, c'est-\`a-dire
que le signe alterne spontan\'ement; il r\'esulte des
d\'efinitions que
$\tilde{\pi}$ est l'une des repr\'esentations
''compl\'ementaires'' de l'\'enonc\'e du th\'eor\`eme celle
qui v\'erifie en plus que le caract\`ere prend des valeurs
diff\'erentes sur le bloc
$(\rho,B,B)$ et $(\rho,B-1,B-1)$, ceci n'est de toute
fa\c con possible que sous les m\^emes conditions que dans
l'\'enonc\'e du th\'eor\`eme \`a cause de la condition de
signe (cf. Remarque de \ref{2construction}).

Supposons
maintenant que $B>0$, en particulier $(\rho,0,0,\zeta) \in Jord(\psi)$. On suppose que
$D=A-2 inf(B, A-B+1)$. On suppose en plus que  si $A-B+1> B$, $\lambda\neq \epsilon(\rho,0,0,\zeta)$.

 On reprend les constructions et en particulier
$t:=(A-D)/2$, notation d\'ej\`a introduite
avant l'\'enonc\'e avec soit
$t=t_1$ dans le cas o\`u $\lambda_{1}\neq
\epsilon(\rho,B-1,B-1)$, soit
$t=t_1+ inf(B,A-B-2t_{1}+1)$ dans le cas oppos\'e puisque
le signe n'alterne plus. Reprenons l'hypoth\`ese $D=A-2inf(B,A-B+1)$ qui dit exactement que $t=inf
(B,A-B+1)$. Ce que l'on veut d\'emontrer ici est que
$t_1=0$ et $\tilde{\pi}$ sera alors aussi une des
repr\'esentations compl\'ementaires, celle que l'on n'a pas trouv\'ee ci-dessus. 
 
   On a s\^rement $t_{1}\leq (A-B+1)/2$ donc si
$A-B+1\leq B$, on a $t=A-B+1 >t_{1}$ d'o\`u $t=t_{1}+ A-B-2t_{1}+1= A-B-t_{1}+1$ d'o\`u encore $t_{1}=0$. Supposons donc que $B< A-B+1$. Ici on a $t=B$ et soit imm\'ediatement $t_{1}=0$ soit $B> A-B-2t_{1}+1$. Dans ce cas, quand on reprend la deuxi\`eme construction (\ref{2construction}) on constate que l'on ne touche pas au bloc $(\rho,0,0,\zeta)$ et donc que $\lambda=\epsilon(\rho,0,0,\zeta)$ ce qui est contraire \`a l'hypoth\`ese et emp\^eche l'in\'egalit\'e $B>A-B-2t_{1}+1$.
    D'o\`u le r\'esultat cherch\'e.
\nl
On reprend les notations qui pr\'ec\`edent, en particulier $(\rho,A,B,\zeta)$ et $\epsilon_{0}$. 
\nl    
\bf Lemme: \sl soit $\lambda$ tel que $\prod_{C\in [B,A]}((-1)^{[C]}\lambda)=\epsilon_{0}$ (on suppose qu'il existe un tel $\lambda$ alors la repr\'esentation $\pi_{\lambda}:=\pi(\psi',\epsilon',\cup_{C\in [B,A]}(\rho,C,C,\zeta,(-1)^{[C]}\lambda))$ intervient avec multiplicit\'e 1 exactement dans $\pi(\psi,\epsilon)$.\rm

On reprend la d\'efinition de $\pi(\psi,\epsilon)$; on consid\`ere d'abord les  sous-quotients intervenant dans $$\sum_{C\in ]\zeta B, \zeta A]}
<\rho\vert\,\vert^{\zeta B}, \cdots, \rho\vert\,\vert^{-\zeta C}> \times Jac_{\zeta(B+2), \cdots, \zeta C} \pi(\psi',\epsilon',(\rho, A,B+2,\zeta,\epsilon_{0}).
$$
 Ces sous-quotients v\'erifient $t_{1}\neq 0$ dans la premi\`ere construction de \ref{1construction}    et ils ne peuvent \^etre isormorphes \`a une repr\'esentation du type $\pi_{\lambda}$ de l'\'enonc\'e. Ce qu'il faut donc d\'emontrer et que $\pi_{\lambda}$ intervient avec multiplicit\'e 1 exactement dans
 $$
 \oplus_{\eta=\pm}(-1)^{[(A-B+1)/2]}\eta^{A-B+1}\epsilon_{0}^{A-B}\pi(\psi',\epsilon',(\rho,A,B+1,\zeta,\eta),(\rho,B,B,\zeta,\eta\epsilon_{0})).\eqno(*)
 $$
 Le cas facile est celui o\`u $A=B+1$; la repr\'esentation $\pi_{\lambda}$ n'existe que si $\epsilon_{0}=-1$ et les 2 valeurs de $\lambda$ sont alors possibles; les repr\'esentations $\pi_{\lambda}$ sont exactement celles qui interviennent dans (*) et le seul point \`a v\'erifier est:
 $$
 (-1)^{[(A-B+1)/2]}\eta^{A-B+1}\epsilon_{0}^{A-B}=(-1)\epsilon_{0}=1.
 $$Dans le cas g\'en\'eral,
 on applique le th\'eor\`eme de \ref{enonce} par r\'ecurrence en utilisant $(\rho,A,B+1,\eta)$; pour $\eta$ fix\'e et pour $\lambda'$ v\'erifiant $\prod_{C\in ]B+1,A]}((-1)^{[C]}\lambda')
=\eta$, la repr\'esentation $\pi_{\eta,\lambda'}:=\pi(\psi',\epsilon', \cup_{C\in ]B+1,A]}(\rho,C,C    ,\zeta,(-1)^{[C]}\lambda'),(\rho,B,B,\zeta,\eta\epsilon_{0}$ intervient avec multiplicit\'e 1 exactement dans $
\pi(\psi',\epsilon',(\rho,A,B+1,\zeta,\eta),(\rho,B,B,\zeta,\eta\epsilon_{0}))$. La repr\'esentation $\pi_{\lambda}$ ne peut \^etre que l'une de ces repr\'esentations $\pi_{\eta,\lambda'}$ et elle est cette repr\'esentation quand on a les \'egalites:
$$
\lambda (-1)^{[B]}=\eta \epsilon_{0}, \qquad (-1)^{[B+1]}\lambda' =- \eta\epsilon_{0}.
$$
Cela entra\^{\i}ne, $\lambda'=\lambda$. On r\'ecrit la propri\'et\'e de $\lambda$:
$$
\prod_{C\in [B,A]}((-1)^{[C]}\lambda)=\epsilon_{0}=\lambda^{A-B+1}(-1)^{(A+B-\delta)(A-B+1)/2},
$$
o\`u $\delta=0$ si $A,B$ sont des entiers et -1 sinon. Consid\'erons d'abord le cas o\`u $A-B+1$ est pair. On a alors, $(-1)^{(A-B+1)/2}=\epsilon_{0}$ et pour chaque valeur de $\lambda$ il existe $\eta$ tel que $\pi_{\lambda}=\pi_{\eta,\lambda}$; le point est de v\'erifier le signe, c'est -\`a-dire:
$$
(-1)^{(A-B+1)/2}\eta^{A-B+1}\epsilon_{0}^{A-B}=(-1)^{(A-B+1)/2}\epsilon_{0}^{A-B}=1,
$$
comme cherch\'e.

Supposons maintenant que $A-B+1$ est impair; on a alors une seule valeur de $\lambda$ possible et elle d\'etermine $\eta$ puisqu'il faut $(-1)^{[B]}\lambda=\eta \epsilon_{0}$. Il faut encore v\'erifier le signe qui vient devant cette repr\'esentation; on remarque que dans $\pi_{\lambda}$ le signe d\'efinissant la repr\'esentation alterne sur tous les blocs de Jord de la forme $(\rho,C,C,\zeta)$ avec $C\in [B,A]$ et en valant $\eta\epsilon_{0}$ sur $(\rho,B,B,\zeta)$ et avec le produit de tous ces signes valant $\epsilon_{0}$. Cela donne la relation:
$$
(\eta\epsilon_{0})^{A-B+1}(-1)^{[(A-B+1)/2]}=\epsilon_{0}.
$$
Ainsi 
$$
 (-1)^{[(A-B+1)/2]}\eta^{A-B+1}\epsilon_{0}^{A-B}=+1,
 $$
 comme cherch\'e. Cela termine la preuve.
 
 \subsubsection{Modules de Jacquet  \label{unicite}} 
On refixe les notations tout en gardant $(\rho,A,B,\zeta) \in Jord(\psi)$ avec $A>B+1$ ainsi que les hypoth\`eses de \ref{hypotheses}. Pour simplifier la d\'emonstration et parceque c'est le cas le plus difficile, on suppose que $A$ et $B$ sont des entiers. Pour $D$ un entier avec $D\in [-1,A]$, on note ${\cal S}_D$ l'ensemble non ordonn\'e $\cup_{C \in ]D,A]}[-C,C]$. Soit $\tilde{\pi}$ un constituant de $\pi(\psi,\epsilon)$. On reprend $t,\lambda$ comme en \ref{2construction} en posant $D=A-2t$. Rappelons que $t$ est un entier.  On a aussi construit l'ensemble ordonn\'e ${\cal E}$ tel que 
 $$
 \tilde{\pi} \hookrightarrow \times {x\in {\cal E}}\rho\vert\,\vert^x \times \pi(\tilde{\psi}_{t},\tilde{\epsilon}_{t,\lambda}).
 $$
 Ici $\tilde{\psi},\tilde{\epsilon}_{t,\lambda}$ se d\'eduit de $\psi,\epsilon$ en enlevant tous les blocs de la forme $(\rho,A',B',\zeta)$ o\`u $B'\leq B$ et en rajoutant les blocs $(\rho,C,C,\zeta)$ o\`u $C\in [0,D]$ (on ne rajoute rien si $D=-1$) et le caract\`ere $\tilde{\epsilon}_{t,\lambda}$ alterne sur tous ces blocs en commen\c{c}ant par $\tilde{\epsilon}_{t,\lambda}(\rho,0,0,\zeta)=\lambda$. On sait aussi que 
 $$
 {\cal E} \cup -{\cal E}={\cal S}_{D}
 $$
 En g\'en\'eral cet ensemble ${\cal E}$ n'est pas unique et on n'a pas encore d\'emontr\'e que $t,\lambda$ le sont. Toutefois, on a:
 \nl
 \bf Lemme: \sl
  Soit $\tilde{\pi}$ un constituant de
$\pi(\psi,\epsilon)$ et soient $t$ et $\lambda$, comme en
\ref{2construction}. Soit aussi un ensemble ${\cal E}$ d'entiers, 
  en valeur absolue inf\'erieurs ou \'egaux \`a
$A$ et tels que
$\tilde{\pi}$ soit un sous-quotient de l'induite
$\times_{x\in {\cal E}}\rho\vert\,\vert^x\times
\pi(\tilde{\psi}_t,\tilde{\epsilon}_{t,\lambda})$.
Alors ${\cal E}\cup -{\cal E}={\cal S}_D$, o\`u $D=A-2t$ \'egalit\'e
d'ensembles non ordonn\'es. De plus $t$ et $\lambda$ sont
uniquement d\'etermin\'es par $\tilde{\pi}$.

\nl
\rm
L'existence a \'et\'e prouv\'ee en \ref{2construction}.
L'unicit\'e r\'esulte d'un calcul de module de Jacquet; en
effet fixons $t$, $\lambda$ et ${\cal E}$ convenant et
un autre triplet $t'$, $\lambda'$, ${\cal E}'$ convenant
aussi. On \'ecrit donc $\tilde{\pi}$ d'une part comme
sous-quotient de $$\times_{x\in {\cal
E}}\rho\vert\,\vert^x\times
\pi(\tilde{\psi}_t,\tilde{\epsilon} _{t,\lambda})
\eqno(1)$$d'autre part comme sous-quotient de 
$$
\times_{x\in {\cal
E}}\rho\vert\,\vert^x\times
\pi(\tilde{\psi}_{t'},\tilde{\epsilon}
_{t',\lambda'}).
\eqno(2)$$
Les formules standards calculant les modules de Jacquet sont
particuli\`erement simples \`a appliquer puisque
$$\forall y;
\vert y\vert
\leq A,
Jac_{y}\pi(\tilde{\psi}_t,\tilde{\epsilon}
_{t,\lambda})=0=
Jac_{y}\pi(\tilde{\psi}_{t'},\tilde{\epsilon}
_{t',\lambda'}).$$ Avant de proc\'eder \`a ces calculs
on r\'eordonne ${\cal E}$ en changeant aussi
\'eventuellement des signes pour les \'el\'ements de cet
ensemble de fa\c con \`a ce que $\tilde{\pi}$ soit un
sous-module de (1); cela est possible
gr\^ace \`a \ref{rappel} (propri\'et\'e 4). Cela nous assure
que
$Jac_{x\in {\cal E}}\tilde{\pi}\neq 0$. On d\'eduit de (1)
que
$Jac_{x\in {\cal
E}}\tilde{\pi}$ est non nul et
isotypique de type $\pi(\tilde{\psi}_t,\tilde{\epsilon}
_{t,\lambda})$ et de (2) que ce module de Jacquet est
certainement nul si ${\cal E}\cup -{\cal E}\neq {\cal E}'\cup
-{\cal E}'$ et s'il n'est pas nul est isotypique de type 
$\pi(\tilde{\psi}_{t'},\tilde{\epsilon}
_{t',\lambda'})$.
Cela donne l'unicit\'e cherch\'ee de
$t$ et $\lambda$. et de 
${\cal E}\cup -{\cal E}$. Puisque pour l'un des choix cette
union est 
$S_D$  cela est vrai pour tous les choix.

 \subsubsection{\label{ordremaximal}}
Ici on fait l'hypoth\`ese de \ref{hypotheses}.  On suppose que $\zeta=+$ sinon il faut le faire intervenir partout et ce qui est croissant devient d\'ecroissant.
 \nl
On fixe  $\tilde{\pi}$ et ${\cal E}, t,\lambda$ tels que
\ref{2construction} (1) et (2) soient v\'erifi\'es. On
renote ici (1):
$$
\tilde{\pi}\hookrightarrow \times_{x\in {\cal
E}}\rho\vert\,\vert^x\times
\pi(\tilde{\psi}_t,\tilde{\epsilon}_{t,\lambda}).
\eqno(1) 
$$
L'ensemble
${\cal E}$ est totalement ordonn\'e, ordre not\'e $>_{\cal
E}$ pour garder
$>$ pour l'ordre sur les nombres. On dit que
${\cal E}$ est muni d'un ordre maximal si ${\cal E}$
s'\'ecrit comme une union de segments croissants, pour $i\in
[1,v]$,
$v$ convenable, de la forme $[x_i,y_i]$ avec $x_1 \geq
\cdots \geq x_v$ et pour $i\leq j$ avec $x_i=x_j$, $y_i\geq
y_j$; l'ordre de ${\cal E}$ \'etant l'ordre $[x_1,y_1],
\cdots, [x_v,y_v]$. A un tel ensemble ordonn\'e de segments,
Zelevinski a associ\'e une repr\'esentation irr\'eductible
$\sigma_{>_{\cal E}}$ du groupe lin\'eaire convenable et de
fa\c con standard, quand on a une inclusion
$$\tilde{\pi}\hookrightarrow \times _{x\in {\cal
E}}\rho\vert\,\vert^x\times \pi_{\delta,\zeta} \eqno(1)$$on
peut r\'eordonner ${\cal E}$ de fa\c con \`a avoir une
inclusion 
$$
\tilde{\pi}\hookrightarrow \sigma_{>_{\cal E}}\times
\pi_{\delta,\zeta}.\eqno(2)
$$
\nl

\bf Lemme: \sl Pour ${\cal E}$ v\'erifiant (1), il existe un
unique ordre maximal sur cet ensemble tel que (2) soit
  satisfait. Pr\'ecis\'ement, $x_1=B$, pour
tout
$i\in [1,v[$, $x_i=x_{i+1}+1$ et $y_1>\cdots >y_v$, ce qui
d\'etermine uniquement l'ordre et induit m\^eme des
restrictions sur la forme de ${\cal E}$.\rm
\nl
On fixe un ordre maximal sur ${\cal E}$, d'o\`u les segments
$[x_i,y_i]$ pour $i\in [1,v]$.

Ayant mis cet ordre, on montre maintenant que pour tout
$i\in [1,v[$, $x{i+1}=x_i-1$. En effet s'il n'en est pas
ainsi, on note ${\cal E}_i$ l'ensemble qui se d\'eduit de
${\cal E}$ en enlevant simplement $x_{i+1}$; c'est encore
naturellement une union de segments ordonn\'es suivant
l'ordre de Zelevinski \`a  ceci pr\`es que le $i+1$-i\`eme
segment peut \^etre plus ''grand'' que certains de  ceux qui
le pr\'ec\`edent; cela se produit si $x_i=x_{i+1}$. Le
$i+1$-i\`eme segment dans ${\cal E}_i$ est alors
$]x_{i+1},y_{i+1}]$ et ce segment n'est pas li\'e au sens
de Zelevinski avec ceux qui le pr\'ec\`ede et ont pour
origine $x$ avec $x=x_i$. On note $\sigma_i$ la
repr\'esentation associ\'ee par Zelevinski \`a ${\cal E}_i$.
On suppose que $i$ est minimum avec la propri\'et\'e
$x_i\neq x_{i+1}+1$ et on montre par les arguments bien
connus dans les groupes lin\'eaires que 
$$  \sigma_{\cal E} \hookrightarrow
\rho\vert\,\vert^{x_{i+1}}\times \sigma_{{\cal E}_i}.$$
Remarquons  que $Jac_{x_{i+1},x_1}\sigma_i\neq 0$ et cela
entra\^{\i}ne donc que
$Jac_{x_{i+1},x_1}\tilde{\pi}\neq 0$. Comme
$Jac_{x_1}\tilde{\pi}\neq 0$, on sait a priori que
$x_1=B$  et la m\^eme propri\'et\'e pour
$x_{i+1}$. On v\'erifie  que
$Jac_{B,B}\tilde{\pi}=0$ gr\^ace \`a
\ref{jacquet}. D'o\`u l'assertion sur les $x_i$.

Montrons
maintenant que pour tout $i\in [1,v[$, $y_i> y_{i+1}$.
Supposons qu'il n'en soit pas ainsi et fixons $i$ minimum
tel que $y_i\leq y_{i+1}$. Alors le i-\`eme segment 
de ${\cal E}$ n'est pas li\'e au i+1-i\`eme. On note ici
${\cal E}_i$ l'ensemble des segments constituant ${\cal E}$
sauf que l'on a enlev\'e $x_{i+1}$ au i+1-i\`eme et que l'on
a permut\'e le i-\`eme et le i+1-\`eme segment. On note
$\sigma_{i}$ la repr\'esentation associ\'ee par Zelevinski
\`a l'ensemble des segments constituant ${\cal E}_i$ et on a
encore une inclusion:
$$
\sigma_{\cal E}\hookrightarrow
\rho\vert\,\vert^{x_{i+1}}\times \sigma_i.
$$
On en d\'eduit que $Jac_{x_{i+1}}\tilde{\pi}\neq 0$ et
$x_{i+1}=B$ ce qui est impossible.
Cela termine la preuve du lemme.
\subsubsection{Non nullit\'e de certains modules de Jacquet\label{descriptiontableau}}
On garde l'hypoth\`ese de \ref{hypotheses}.
On fixe ${\cal E}$ satisfaisant les conditions du paragraphe
pr\'ec\'edent et en particulier ${\cal E}$ peut se voir
comme une union de segments. On
\'ecrit
${\cal E}$ sous forme de tableau o\`u les
\'el\'ements croissent de 1 sur chaque ligne et
d\'ecroissent de 1 sur chaque colonne.
$$
\begin{matrix}
&x_1=B &\cdots &\cdots &\cdots &y_1\\
&\vdots&&&\vdots\\
&x_v&\cdots &y_v
\end{matrix}
$$
On note $s$ le nombre de colonnes et $z_1, \cdots, z_s$ le
dernier \'el\'ement de chaque colonne. Cette pr\'esentation
permet de calculer ais\'ement ${\cal E}\cup -{\cal E}$. On a
$$
{\cal E} \cup -{\cal E}=\cup_{x\in [B,y_1]}[-x,x]
\cup_{i\in [1,s]; z_i\leq 0} [z_i,-z_i] \ominus _{i\in
[1,s];z_i>0}[-z_i+1,z_i-1].
$$
On rappelle que ${\cal E}\cup -{\cal E}=\cup_{z\in
[D+1,A]}[-z,z]$ (o\`u $D=A-2t$). On en d\'eduit donc:
$$
[B,y_1]\cup _{i\in [1,s]; z_i\leq 0}\{-z_i\}=
[D+1,A] \cup_{i\in [1,s];z_i>0}
\{z_i-1\}.\eqno(1)$$

Pour $tilde{\pi}$ un constituant de $\pi(\psi,\epsilon)$
fix\'e, en g\'en\'eral plusieurs choix de ${\cal E}$ sont
possibles. Il n'y a pas d'int\'er\^et \`a les \'ecrire tous
mais 2 vont nous servir plus particuli\`erement. D'abord
remarquons que si ${\cal E}$ n'est pas vide, ce qui est
implicite ici, ${\cal E}$ contient soit $-A$ soit
$A$. On rappelle que par d\'efinition $t$ v\'erifie
$2t=A-D$ et que n\'ecessairement $t\leq A-B+1$.
\nl
Soit $\tilde{\pi}$ fix\'e un composant
irr\'eductible de $\pi(\psi,\epsilon)$; on suppose que le $t$ et $\lambda$ qui lui sont associ\'es, v\'erifient, $t\neq 0$
\nl
\bf Lemme: \sl 
on peut alors choisir  ${\cal E}$ comme ci-dessus convenant pour
$\tilde{\pi}$ avec l'une des propri\'et\'es ci-dessous:

\nl

${\cal E}$ contient
$-A$;

\nl
ou

${\cal E}=\begin{matrix}
&x_1=B &\cdots&\cdots&\cdots&A\\
&\vdots &\vdots &\vdots &\vdots &\vdots\\
&-B+1 &\cdots &\cdots &\cdots &(D+1)
\end{matrix};$

\nl
ou

 ${\cal E}=\begin{matrix}
&x_1=B &\cdots&\cdots&\cdots&A\\
&\vdots &\vdots &\vdots &\vdots &\vdots\\
&-B+1 &\cdots &B-t &\cdots &-(D+1)
\end{matrix};$ 

\nl 
ou

${\cal E}=$
${\cal E}_{D}:=
\begin{matrix}
&x_1=B &\cdots&\cdots&\cdots&A\\
&\vdots &\vdots &\vdots &\vdots &\vdots\\
&B-2t+1 &\cdots &B-t &\cdots &D+1\\
&\vdots &\vdots &\vdots\\
&-B+1 &\cdots &t-B
\end{matrix}
$
\nl
De plus si un ensemble ${\cal E}$ convient qui est tel que
le plus petit \'el\'ement de ${\cal E}$ est inf\'erieur ou
\'egal \`a $-B$, alors $Jac_{B, \cdots,
-A}\tilde{\pi}\neq 0$.\rm 
\nl
On fixe
${\cal E}$ et l'ordre maximal sur
${\cal E}$ tel que
$$\tilde{\pi}\hookrightarrow \sigma_{>_{\cal E}}\times
\pi_{D,\lambda}.$$ 
Et on suppose que ${\cal E}$ est tel que parmi tous les
choix possibles le nombre d'\'el\'ements de ${\cal E}$
strictement positifs est minimal. 

Soit  $i\in [1,v]$ tel que $y_i>0$, pour
tout $j>i$, $y_j\neq y_i-1$ et
$y_i\neq D+1$. Alors on montre que
l'on peut remplacer ${\cal E}$ par un ensemble co\"{\i}ncidant \`a l'ordre pr\`es avec  ${\cal E}-\{y_i\}\cup
\{-y_i\}$ : on proc\`ede comme en \ref{ordremaximal}, en
notant ${\cal E}_i$ l'ensemble des segments qui se d\'eduit
de ${\cal E}$ en rempla\c cant
$y_i$ par $y_i-1$ pour le i-\`eme segment. On note
$\sigma_i$ la repr\'esentation associ\'ee par Zelevinski \`a
cet ensemble ordonn\'e de segments et on a encore:
$\sigma_{\cal E} \hookrightarrow \sigma_i\times
\rho\vert\,\vert^{y_i}$ ici on utilise seulement le fait que
pour tout $j >i$, $y_j<y_i-1$ ou $y_j=y_i$. Comme $y_i\neq
D+1$ la repr\'esentation
$\rho\vert\,\vert^{y_i}\times \pi_{\delta,\zeta}$ est
irr\'eductible (cf. \ref{irreductibilite}) et donc 
isomorphe
\`a
$\rho\vert\,\vert^{-y_i}\times \pi_{\delta,\zeta}$. D'o\`u
l'assertion annonc\'ee. 

On suppose que $\cal E$ ne contient pas $-A$; il contient donc $A$ et n\'ecessairement $y_1=A$. Ainsi le nombre de colonnes du tableau (le $s$ dans ce qui pr\'ec\`ede l'\'enonc\'e du lemme) vaut $A-B+1$, ce que l'on va \'ecrire $b:=A-B+1$. (C'est le $b$ dans l'interpr\'etation en terme de bloc de Jordan venant des repr\'esentations).
On revient aux notations $z_j; j\in
[1,b]$ qui pr\'ec\`edent l'\'enonc\'e et on remarque que
$\{z_j-1; z_j>0\}$ est d'apr\`es ce que l'on vient de voir 
 un
intervalle de la forme $[z,D]$, en acceptant
$z=D+1$ si l'ensemble est vide. Ainsi
\ref{descriptiontableau} (1) que l'on rappelle, en acceptant
$z=D+1$:
$$
[B,A] \cup_{i\in [1,s]; z_i\leq 0}\{-z_i\}=
[D+1,A] \cup [z,D]=
[z,A].
$$
D'o\`u $z\leq B$ et $[z,B[=\{-z_i; z_i\leq 0\}.$On pourrait avoir $z=B$ ce qui veut dire que tous les \'el\'ements de ${\cal E}$ sont strictement positifs mais on obtient tout de suite une impossibilit\'e, il faudrait $B+1 \leq B$!.
Ainsi $z<B$
Cela veut dire que le tableau s'\'ecrit si $z= 0$
$${\cal E}=
\begin{matrix}
&B&\cdots &\cdots &\cdots &A\\
&\vdots &\vdots &\vdots &\vdots &\vdots\\
&B-2t+1 &\cdots &\cdots &\cdots &D+1\\
\end{matrix}
$$Mais on doit avoir $[B,A]\cup [-(B-2t+1),0]=[0,A]$ ce qui force $B-2t+1=-B+1$ et on est dans le 2e cas de l'\'enonc\'e.

Si $z\neq 0$
$${\cal E}=
\begin{matrix}
&B&\cdots &\cdots &\cdots &A\\
&\vdots &\vdots &\vdots &\vdots &\vdots\\
&B-2t+1 &\cdots &z &\cdots &D+1\\
&\vdots &\vdots &\vdots\\
&- B+1 &\cdots &-z
\end{matrix}\eqno(*)
$$
En particulier si
$z=D+1$
 le tableau ${\cal E}$ est rectangulaire si $z=D+1$ et on est dans le 3e cas de l'\'enonc\'e. 

Supposons maintenant que $\{z_j; z_j>0\}$ n'est pas un
ensemble vide; on calcule $z$ en comptant le nombre de
lignes. Ce nombre est $2B$ qui se d\'ecompose en $2t$
''longues'' lignes et $2B-2t$ lignes ''courtes''. Mais
$2z$ est aussi le nombre de lignes ''courtes'' d'o\`u
$z=B-t$
 Remarquons qu'il faut donc
$t\leq B$ et que $t=B$ est \'equivalent \`a ce
que le tableau d\'efinissant ${\cal E}$ soit rectangulaire. 

On d\'emontre maintenant la fin du lemme; en effet
supposons qu'un ensemble ${\cal E}$ qui convient contient un
\'el\'ement inf\'erieur ou \'egal \`a $-B$. On impose
en plus \`a ${\cal E}$ v\'erifiant cette propri\'et\'e
d'avoir un nombre minimal d'\'el\'ements strictement
positifs. On a donc vu que ou bien on a directement
l'assertion, ou bien
${\cal E}$ est celui de l'\'enonc\'e 2e, 3e ou 4e cas. Dans ces  cas, le plus petit
\'el\'ement de ${\cal E}$ est $-B+1$; ce qui est donc
exclu.

\subsubsection{ \label{nullitepourcertaintableau}}
On prend les hypoth\`eses et les notations de la sous-section pr\'ec\'edente; en particulier $\tilde{\pi}$ est fix\'e, d'o\`u $t,\lambda$ et $D=A-2t$. On a aussi d\'efini ${\cal E}_{D}$.

\nl
\bf Lemme: \sl On suppose que
$D< B-1$ et qu'il existe $C\in ]B,A]$ tel que $\tilde{\pi}$ est un sous-quotient de l'induite:
$$X_{C}:=
<\rho\vert\,\vert^{\zeta B}, \cdots, \rho\vert\,\vert^{-\zeta C}>\times Jac_{\zeta(B+2),\cdots, \zeta C}\pi(\psi',\epsilon',(\rho, A-1,B+1,\zeta,\epsilon_{0}))>.
$$
 Alors
$Jac_{x\in {\cal E}_{D}}\tilde{\pi}=0$.
\rm
\nl
 Dans la preuve, on suppose que $\zeta=+$.
Avant de d\'emontrer ce lemme on va faire de fa\c con un peu
g\'en\'erale le calcul suivant. Soit $\ell$ un entier non
nul et on pose ${\cal A}_\ell$:
$$
\begin{matrix}
&B &B+1 &\cdots &A\\
&\vdots&\vdots&\vdots&\vdots\\
&B-\ell+1 &B-\ell+2 &\cdots &A-\ell+1\\
&B-\ell\\
&\vdots\\
&-B+1
\end{matrix}
$$
Alors, si $\ell\leq B-C+1$, on pose $$
{\cal A}_{C,\ell}:=
\begin{matrix}
&B+2 &\cdots &\cdots &\cdots &A\\
&\vdots &\vdots &\vdots &\vdots &\vdots\\
&B+3-\ell &\cdots &\cdots &\cdots &A-\ell+1\\
&B+2-\ell &\cdots &C-\ell
\end{matrix}
$$et
$$Jac_{x\in {\cal A}_{\ell}}X_C=
<\rho\vert\,\vert^{-B},
\cdots,
\rho\vert\,\vert^{-C+\ell}> \times Jac_{x\in {\cal
A}_{C,\ell}}
\pi(\psi',\epsilon',(\rho,A,B-2,\zeta=+,\epsilon_{0}))
$$o\`u le premier facteur de l'induite n'appara\^{\i}t pas
si $\ell=C-B+1$.

Si $\ell > C-B+1$, 
$$Jac_{x\in {\cal A}_{\ell}}X_C= 0
$$
Pour faire ces calculs, on commence, comme on en a le droit,
par calculer $Jac_{B, \cdots, -B+1}X_C$ et cela
vaut $$X'_C:=<\rho\vert\,\vert^{-B}, \cdots,
\rho\vert\,\vert^{-C}>\times Y_C,$$ o\`u 
 $$
 Y_{C}=
  Jac_{(B+2),\cdots,  C}\pi(\psi',\epsilon',(\rho, A-1,B+1,\zeta=+,\epsilon_{0}))>.
  $$
Il faut encore
calculer $Jac_{x\in {\cal B}}$ du r\'esultat, $X'_C$, o\`u
$${\cal B}:= \begin{matrix}
&B+1 &\cdots &A\\
&\vdots &\vdots &\vdots \\
&B+2-\ell &\cdots &A-\ell+1.
\end{matrix}
$$
Dans le tableau ci-dessous le $X$ dit simplement que l'on
omet l'\'el\'ement situ\'e en principe \`a cette place et on
pose, pour $\ell'\leq C-B+1$
$$
{\cal C}_{C,\ell'}:=
\begin{matrix}
&B+2 &\cdots &C\\
&B+1 &\cdots &C-1 &X &C+1 &\cdots
&A-1\\
&\vdots &\vdots &\vdots &\vdots &\vdots &\vdots &\vdots \\
&B+2-\ell' &\cdots &C-\ell' &X &C
-\ell' +2 &\cdots &A-\ell'+1
\end{matrix}
$$On montre que pour $\ell' \leq C-B+1$, $$Jac_{x\in {\cal
B}_{\ell'}}X'_C =<\rho\vert\,\vert^{-B}, \cdots,
\rho\vert\,\vert^{-C+\ell'}>\times Jac_{x
\in {\cal C}_{C,\ell'}}
\pi(\psi',\epsilon',(\rho,A,B+2,\epsilon_{0})).
$$Le terme de droite intervient clairement dans le terme de
gauche mais en principe interviennent aussi les termes
index\'es par $\ell''$ suivants o\`u $\ell''\leq \ell'$:
$$
<\rho\vert\,\vert^{-B},\cdots,
\rho\vert\,\vert^{-C+\ell''}> \times Jac_{x\in {\cal
C}_{C,\ell',\ell''}}
\pi(\psi',\epsilon',(\rho,A,B+2,\zeta=+,\epsilon_{0})),
$$
o\`u 
$
{\cal C}_{C,\ell',\ell''}:=$:
$$
\begin{matrix}
&B+2 &\cdots &C\\
&B+1 &\cdots &C-1 &X &C+1 &\cdots
&A\\
&\vdots &\vdots &\vdots &\vdots &\vdots &\vdots &\vdots \\
&B+2-\ell'' &\cdots &C-\ell''+1
&X&C-\ell'' +3&\cdots &A-\ell''+2\\
&B+1-\ell'' &\cdots &C-\ell''
&C-\ell''+1&C-\ell'' +2&\cdots
&A-\ell''+1\\
&\vdots &\vdots &\vdots &\vdots &\vdots &\vdots &\vdots \\
 &B-\ell '+2
&\cdots &C-\ell' &C-\ell'+1 &C -\ell'
+2 &\cdots &A-\ell'+1
\end{matrix}
$$
Mais pour un tel choix de $\ell''$, on va montrer que 
$Jac_{x\in {\cal C}_{C,\ell',\ell''}}
\pi(\psi',\epsilon',(\rho,a+2,b-2,\epsilon(\rho,a,b)))=0$.
En effet, on consid\`ere uniquement la partie du tableau:
$$
{\cal C}'':= 
\begin{matrix}
&B+2 &\cdots &C\\
&B+1 &\cdots &C-1\\
&\vdots &\vdots &\vdots\\
&B+2-\ell''&\cdots &C-\ell''+1\\
&B+1-\ell''&\cdots&C-\ell''&C-\ell''+1.
\end{matrix}
$$
Et le $Jac$ dont on cherche \`a montrer la nullit\'e se
factorise par $Jac_{x\in {\cal C}''}$. Or une
repr\'esentation irr\'eductible $\sigma$ du $GL$
correspondant dont le support cuspidal est l'ensemble des
$\rho\vert\,\vert^x$ pour $x$ parcourant ${\cal C}''$
d'apr\`es la classification de Zelevinski ne peut pas avoir
$Jac_{x}\sigma\neq 0$ pour uniquement $x=B+2$. D'o\`u la
contradiction, il y aurait une  valeur de $x\neq B+2$,
$x\in {\cal C}''$ telle que $Jac_x
\pi(\psi',\epsilon',(\rho,A,B+2,\zeta=+,\epsilon_{0}))\neq
0$. Ceci est exclu par \ref{moduledejacquet} d'o\`u la
nullit\'e cherch\'ee.

On obtient, comme annonc\'e,  le
cas $\ell \leq C-B+1$ en faisant
$\ell'=\ell$ dans ce qui pr\'ec\`ede et en remarquant de
fa\c con tout
\`a fait formelle que  l'on peut remplacer
${\cal C}_{C,\ell}$ par ${\cal A}_{C,\ell}$.

Supposons maintenant $\ell>C-B+1$ et il suffit m\^eme de
montrer la nullit\'e pour $\ell=C-B+2$. On a d'apr\`es ce qui
pr\'ec\`ede:
$$
Jac_{x\in  {\cal A}_{\ell}}X_j= Jac_{x\in {\cal A}'}\pi
(\psi',\epsilon',(\rho,A,B+2,\zeta=+,\epsilon_{0})),
$$
o\`u $$
{\cal A}':=
\begin{matrix}
&B+2 &\cdots &C\\
&B+1 &\cdots &C-1 &X &C+1 &\cdots
&A\\
&\vdots &\vdots &\vdots &\vdots &\vdots &\vdots &\vdots \\
&B+1-(C-B)&\cdots &B-1 &X &B+1 &\cdots
&A-(C-B)\\
&B-(C-B)&\cdots &B-2&B-1 &B &\cdots
&A-(C-B)-1
\end{matrix}
$$
et l'argument d\'ej\`a donn\'e ci-dessus (pour $\ell''$)
donne la nullit\'e cherch\'ee. Revenons \`a l'\'enonc\'e du
lemme; quand $D+1\leq B-1$ (c'est-\`a-dire $D<B-1$), le calcul de $Jac_{x\in
{\cal T}_D}\tilde{\pi}$ se factorise par $Jac_{x\in
{\cal C}_{C,\ell'}}\tilde{\pi}$ pour tout $\ell'$ tel que
$A-\ell'+1 \geq D+1$. On peut donc prendre
$\ell'=A-D > A-B+1\geq C-B+1$. D'o\`u la nullit\'e
annonc\'ee.

\bf Remarque: \sl le m\^eme argument s'applique pour ${\cal
T}_{\pm}$ o\`u $${\cal T}_{\pm}=
\begin{matrix}
B &\cdots & A \\
\vdots &\vdots &\vdots\\
-(B-1) &\cdots &\pm (D+1)
\end{matrix}
$$
\rm
\subsubsection{}
\bf Corollaire: \sl Soit $\tilde{\pi}$ un composant
irr\'eductible de $\pi(\psi,\epsilon)$. Si $D$
associ\'e \`a $\tilde{\pi}$ en \ref{unicite} v\'erifie
$D < B-1$ alors $\tilde{\pi}$ v\'erifie soit
$Jac_{\zeta B, \cdots, -\zeta A}\tilde{\pi}\neq 0$ soit
est l'une des repr\'esentations compl\'ementaires de
\ref{enonce}.
\rm
\nl On suppose comme pr\'ec\'edemment que $\zeta=+$.
On fixe $\tilde{\pi}$ tel que le $D$ qui lui soit associ\'e v\'erifie $D<B-1$ et on suppose aussi que 
$Jac_{ B,\cdots, - A}\tilde{\pi}= 0$. Il faut d\'emontrer qu'un tel $\tilde{\pi}$ est une des repr\'esentations comp\'ementaires de \ref{enonce}. Avec le lemme \ref{nullitepourcertaintableau} et la remarque qui le suite, on sait dej\`a gr\^ace \`a \ref{descriptiontableau} que $\tilde{\pi}$ est un sous-quotient de l'une des repr\'esentations $X_{\eta}:=pi(\psi',\epsilon',(\rho,A,B+1,\zeta,\eta,(\rho,B,B,\zeta,\eta \epsilon_{0}))),$
o\`u $\eta$ est un signe convenable. On fixe $\eta$ convenant.

 On note
$\psi'',\epsilon''$ le couple qui se d\'eduit de
$\psi',\epsilon'$ en enlevant
$(\rho,B-1,B-1,\zeta,\epsilon(\rho,B-1,B-1,\zeta))$. Et on montre que
$
Jac_{B,\cdots, -B+1}X_{\eta}=0$ si
$\eta\epsilon_{0}\neq \epsilon(\rho,B-1,B-1,\zeta)$ et vaut
$\pi(\psi'',\epsilon'', (\rho,A,B+1,\zeta,\eta)))$ sinon (c'est essentiellement ce qui se passe pour les morphismes \'el\'ementaires rappel\'e en \ref{rappel}).

 Si $\eta\epsilon_{0}\neq \epsilon(\rho,B-1,B-1,\zeta)$, il r\'esulte de \ref{descriptiontableau} que $D=A$ et \ref{1construction} ainsi que \ref{2construction} donnent le r\'esultat cherch\'e. Supposons donc que l'on a \'egalit\'e, ce qui d\'efinit $\eta$.
 
Ainsi il existe $\tilde{\pi}'$ un constituant irr\'eductible 
de $\pi(\psi'',\epsilon'',(\rho,a+1,b-1,\eta)))$ et une
inclusion:
$$
\tilde{\pi}\hookrightarrow \times_{x\in
[B,-B+1]}\rho\vert\,\vert^x \times \tilde{\pi}'.
$$
On associe \`a $\tilde{\pi}'$ des donn\'ees $D',\lambda'$ et un ensemble ${\cal E}'$ comme dans \ref{unicite}; n\'ecessairement
$D',\lambda$ sont aussi $D,\lambda$  les donn\'ees
associ\'ees \`a $\tilde{\pi}$ (cf \ref{unicite}). On peut appliquer
\ref{enonce}
\`a $(\psi'',\epsilon'',(\rho,A,B+1,\zeta,\eta))$. Ainsi soit
${\cal E}'$ contient $-\zeta A$ soit $\tilde{\pi}'$ est
l'une des repr\'esentations compl\'ementaires. Si ${\cal E}'$ contient $-  A$ il existe alors ${\cal E}$ convenant pour $\tilde{\pi}$ et contenant aussi $- A$ et on a vu que cela entra\^{\i}ne que $Jac_{ B, \cdots, - A}\tilde{\pi}\neq 0$, ce qui est contraire \`a l'hypoth\`ese. Ainsi $\tilde{\pi}'$ est l'une des repr\'esentations compl\'ementaire pour $(\psi'',\epsilon'', (\rho,A,B+1,\zeta,\eta))$.
 
Concr\`etement, cela veut dire qu'il existe un signe, $\kappa$ tel que $\tilde{\pi}'=\pi(\psi'',\epsilon'', \cup_{C\in [B+1,A]}(\rho,C,C,\zeta,(-1)^C \kappa))$ et $\tilde{\pi}$ d'apr\`es \ref{rappel} est l'une des 2 repr\'esentations, index\'ees par $\kappa'=\pm$ $$\pi(\psi'',\epsilon'',\cup_{C\in [B+1,A]}(\rho,C,C,\zeta,(-1)^C \kappa), (\rho,B,B,\zeta,\kappa'),(\rho,B-1,B-1,\zeta,\kappa').$$
On aura montr\'e que $\tilde{\pi}$ est l'une des repr\'esentations compl\'ementaires si l'on v\'erifie que $(-1)^{B+1}\kappa\neq \kappa'$. Si l'on a \'egali\'e, on a aussi $Jac_{(B+1), \cdots, -\ B}\tilde{\pi}\neq 0$ d'o\`u a fortiori $Jac_{(B+1)}\tilde{\pi}\neq 0$ ce qui est exclu par \ref{proprietedujac}. D'o\`u le corollaire.

\subsubsection{}

\bf Lemme: \sl Soient $D$ et $\lambda$ comme en
\ref{unicite} et on suppose que $D\geq B-1$ ou encore $D+1 >B-1$.
Soit ${\cal T}={\cal E}_{D}$ d'o\`u aussi
 la repr\'esentation $\sigma_{{\cal T}}$.
Alors l'induite
$\sigma_{{\cal T}}\times
\pi(\tilde{\psi}_D,\tilde{\epsilon}_{D,\lambda})$ a
un unique sous-module irr\'eductible.
\rm
\nl
Le tableau ${\cal T}$ a $2B$ lignes chacune form\'ee d'un segment $[x_{i},y_{i}] $ pour $i\in [1,2B]$.
Pour tout $i\in [1,2B]$, on note $\sigma_{\geq i}$ la
repr\'esentation du groupe lin\'eaire convenable qui dans la
classification de Zelevinski correspond aux segments
$[x_j,y_j]$ pour $j\geq i$. On montre progressivement que
pour tout $i\in [1,2B]$, 
$$
Jac_{x\in \cup_{j<i}[x_j,y_j]}\sigma_{\cal T}\times
\pi_{\delta,\zeta}= \sigma_{\geq i}\times \pi_{\psi_{D},\epsilon_{D,\lambda}}.
$$
Pour d\'emontrer cela, on suppose d'abord que $i\leq 2t$; on
a alors $y_i=D+1 +2t-i \geq B$ d'apr\`es
l'hypoth\`ese sur $D$. Si l'assertion n'est pas vraie
pour un tel
$i$, il faut n\'ecessairement qu'il existe $z$ tel que $[z,
y_i]$ soit un segment et $Jac_{z,\cdots, y_i}
\sigma_{\geq i}^*\neq 0$. Mais ceci n'est  possible que si
$y_i \leq -x_v$ qui est le plus grand \'el\'ement de $-{\cal
T}$. Sous l'hypoth\`ese $D+1 \geq B$ on est
assur\'e que $y_{2t}\leq y_i$ mais 
$$y_{2t} =D+1   \geq B > -x_v=B-1.$$

Pour $i>2t$ l'argument est diff\'erent; ici on utilise le
fait que tout \'el\'ement de $[x_i,y_i]$ est strictement
plus petit que $-y_v$; cela se voit sur le tableau. 
\nl
 
Dans tous les cas on note $\pi_{{\cal T},D,\lambda}$ la
repr\'esentation ainsi d\'efinie. On a \'evidemment accept\'e $A-D=2B$ cas o\`u ${\cal E}_{D}$ est rectangulaire.

Ce qui nous reste \`a d\'emontrer maintenant est que la repr\'esentation $\pi_{{\cal T},D,\lambda}$ n'est pas un constituant de $\pi(\psi,\epsilon)$ sauf si $A-D=2B$ et $\lambda=-\epsilon(\rho,0,0,\zeta)$. \subsubsection{\label{lemmeconclusif}}
On garde les hypoth\`eses et notations pr\'ec\'edentes, on rappelle que
$2t=A-D$; on garde l'hypoth\`ese $\zeta=+$ pour clarifier la situation. Pour simplifier les notations, on pose $\pi_{D,\lambda}:=\pi(\psi_{D},\epsilon_{D,\lambda})$. Pour $C\in ]B,A]$, on a d\'efini $X_{C}=<\rho\vert\,\vert^{B}, \cdots, \rho\vert\,\vert^{-C}\times Y_{C}$ o\`u $Y_{C}=Jac_{B+2, \cdots, C}\pi(\psi',\epsilon', (\rho,A,B+2,\zeta=+,\epsilon_{0})$. Pour $\eta=\pm$, on pose $X_{\eta}:=\pi(\psi',\epsilon', (\rho,A,B+1,\zeta=+,\eta),(\rho,B,B,\zeta=+,\eta\epsilon_{0}))$.
\nl
\bf Lemme: \sl On suppose que $D+1>B-1$ et on fixe $C\in ]B,A]$ et $\eta=\pm$.

(i)Alors $\pi_{D,\lambda}$ n'est pas
sous-quotient de $Jac_{x\in {\cal T}_D} X_C$
sauf si $C-B+1=2t$. 

(ii)Supposons que $C-B+1=2t$. On note ${\cal T}'_{D}$ le tableau
$$\begin{matrix}
&B+1 &\cdots &\cdots &\cdots &\cdots &A\\
&\vdots &\vdots &\vdots &\vdots &\vdots &\vdots\\
&B-2t+3 &\cdots &B-t+1&\cdots &\cdots
&D+1\\ 
&B-2t+2 &\cdots &B -t &\cdots &B-1\\
&\vdots &\vdots &\vdots\\
&-B+2 &\cdots &-B+t
\end{matrix}
$$
Alors,
$$Jac_{x\in {\cal T}_D}X_C= Jac_{x\in {\cal T}'_{D}}
\pi(\psi',\epsilon',(\rho,A-1,B+1,\epsilon_{0})).\eqno(1)
$$

(iii) Supposons encore $C-B+1=2t$. On note ${\cal
T}_{A-2,B,D}$ l'analogue de ${\cal T}_D$ mais en
rempla\c cant le couple $(A,B)$ par $(A-2,B)$; $D$
est sans changement d'o\`u $t$ devient $t-1$. La
multiplicit\'e de
$\pi_{D,\lambda}$ comme sous-quotient de $Jac_{x\in {\cal
T}_D}X_C$ est la m\^eme que sa multiplicit\'e comme
sous-quotient de
$Jac_{x\in {\cal T}_{A-2,B,D}}\pi(\psi',\epsilon',
(\rho,A-2,B,\zeta=+,\epsilon_{0}))$.

(iv)Soit $\eta=\pm$ et ${\cal T}_{A-2,B-1,D}$
l'analogue de ${\cal T}_D$ mais en y rempla\c cant
$(A,B)$ par $(A-2,B-1)$. La multiplicit\'e de
$\pi_{D,\lambda}$ dans
$Jac_{x\in {\cal T}_D}X_{\eta}$
est $0$ si $\eta\epsilon_{0}\neq
\epsilon(\rho,B-1,B-1,\zeta=+)$ et \'egal \`a la multiplicit\'e de 
$\pi_{D,\lambda}$ comme sous-quotient de $Jac_{x\in
{\cal T}_{A-2,B-1,D}}\pi(\psi'',\epsilon'',
(\rho,A-2,B-1,\zeta=+,
\eta))$, o\`u ici
$\psi'',\epsilon''$ se d\'eduit de $\psi',\epsilon'$ en
enlevant le bloc $(\rho,B-1,B-1,\zeta=+)$.

(v) La multiplicit\'e de $\pi_{D,\lambda}$ comme
sous-quotient de 
$Jac_{x\in {\cal T}_D}\pi(\psi,\epsilon)$ est 0 sauf
si $t=B$et
$\lambda=-\epsilon(\rho,1,1)$ o\`u elle vaut 1.
\rm
\nl  
(i) 
 est un calcul fait en \ref{nullitepourcertaintableau} en
posant avec les notations de loc. cit $\ell=2t$: on y a
montr\'e que si
$\ell >C-B+1$, le module de Jacquet  est nul. Si
$\ell< C-B+1$ dans le calcul de ce module de Jacquet il y a
une induction avec le facteur
$<\rho\vert\,\vert^{-B}, \cdots, \rho\vert\,
\vert^{-C+\ell-1}$ qui fait que l'ensemble ${\cal
E}$ associ\'e \`a $\tilde{\pi}$ ne peut \^etre r\'eduit
aux \'el\'ements de ${\cal T}_{D}$.  D'o\`u (i).
Supposons donc que $2t=\ell=C-B+1$, le (ii) est le calcul
d\'ej\`a fait  comme expliqu\'e ci-dessus.

Pour montrer (iii), il faut transformer ${\cal T}'_{D}$ qui est
dans le r\'esultat de (ii) en
${\cal T}_{A-2,B,D}$. Pour cela, on note
momentan\'ement
${\cal T}''$ le tableau qui se d\'eduit de ${\cal T}'$ en
supprimant $B-t+1, \cdots, B-1$ de la ligne qui
commence par $B-2t+3$ et on a ais\'ement, pour toute
repr\'esentation $X$
$$
Jac_{x\in {\cal T}'_{D}}X=Jac_{x\in [B-t+1,B-1]}
Jac_{x\in {\cal T}''}X.
$$
Mais $\times_{x\in
[B-t+1,B-1]}\rho\vert\,\vert^x \times
\pi_{D,\lambda}$ est isomorphe \`a $\times_{x\in
[-B+1,-B+t-1]}\rho\vert\,\vert^x\times
\pi_{D,\lambda}$. Ainsi, on peut remplacer ${\cal T}'$
par ${\cal T}''$ auquel on rajoute comme derni\`ere ligne le
segment $[-B+1, -B+t-1]$. Cela termine la preuve
de (iii).

Pour d\'emontrer (iv), on calcule d'abord 
$$
Jac_{x\in
[B,-B+1]}\pi(\psi',\epsilon',(\rho,A,B+1,\zeta=+,\eta), 
(\rho,B,B,\zeta=+,\eta\epsilon_{0})).\eqno(2).
$$
Cela vaut 0 si $\eta\epsilon_{0}\neq
\epsilon(\rho,B-1,B-1,\zeta=+)$ et
$\pi(\psi'',\epsilon'',(\rho,A,B+1,\eta))$ sinon. Ensuite,
on calcule le Jac suivant les 2 premi\`eres lignes en
utilisant \ref{lecastroue} et on trouve directement (iv).

Montrons maintenant (v). Supposons que $\tilde{\pi}$ soit un
constituant irr\'eductible de $\pi(\psi,\epsilon)$ avec
$D$ v\'erifiant $D+1\geq B$; on suppose
aussi que $Jac_{x\in {\cal E}_D}\tilde{\pi}\neq 0$. On
montre que $Jac_{B,\cdots, -A}\tilde{\pi}=0$.
En effet, on note $\sigma'$ la repr\'esentation
irr\'eductible qui correspond au tableau ${\cal E}_D$
dont on a enlev\'e le dernier \'el\'ement des $2t$
premi\`eres lignes. On a ais\'ement
$$
\tilde{\pi}\hookrightarrow \sigma'\times
<\rho\vert\,\vert^{A}, \cdots,
\rho\vert\,\vert^{D+1}>\times \pi_{D,\lambda},
$$
pour $\zeta$ convenable. Or la repr\'esentation $
<\rho\vert\,\vert^{A}, \cdots,
\rho\vert\,\vert^{D+1}>\times \pi_{D,\lambda}$ a
un unique sous-quotient irr\'eductible, $X$,  qui v\'erifie
$Jac_{A, \cdots, D+1}X\neq 0$; c'est  la
repr\'esentation qui correspond au paquet \'el\'ementaire se
d\'eduisant de $\psi_D, \epsilon_{D,\lambda}$, en
rempla\c cant le bloc $(\rho,D,D,\zeta=+)$ par $(\rho,A,A,\zeta=+)$
sans changer le caract\`ere sur ce bloc. Mais cette repr\'esentation
v\'erifie $Jac_{-A}X=0$. Il est n\'ecessaire que
$\tilde{\pi}\hookrightarrow \sigma'\times X$ et on en
d\'eduit l'assertion $Jac_{B, \cdots,
-A}\tilde{\pi}=0$.

Revenons au calcul de la multiplicit\'e de $\pi_{D,\lambda}$ en tant que sous-quotient de $Jac_{x\in {\cal E}_{D}}\pi(\psi,\epsilon)$. D'apr\`es ce qui pr\'ec\`ede c'est exactement la multiplicit\'e avec laquelle il intervient dans $$(-1)^{D-B+1}Jac_{x\in {\cal T}_{A-2,B,D}}\pi(\psi',\epsilon',(\rho,A-2,B,\zeta=+,\epsilon_{0}))$$ additionn\'ee de celle avec laquelle il intervient dans $$
(-1)^{[A-B+1)/2]}\eta^{(A-B+1)}\epsilon_{0}^{A-B} Jac_{x\in {\cal T}_{A-2,B-1}}\pi(\psi'',\epsilon'',(\rho,A-2,B-1,\zeta=+,\eta))$$ o\`u $\eta=\epsilon_{0}\epsilon(\rho,B-1,B-1,\zeta=+)$.

Comme on va le voir, on connait ces multiplicit\'es en appliquant par r\'ecurrence \ref{enonce} et on va pouvoir montrer les propri\'et\'es cherch\'ees.

\nl
Par r\'ecurrence, on
connait la structure de
$\pi(\psi',\epsilon',(\rho,A-2,B,\zeta=+,\epsilon_{0}))$; on va la d\'ecrire et appliquer $Jac_{x\in {\cal T}_{A-2,B,D}}$ et compter la multiplicit\'e de $\pi_{D,\lambda}$ dans le r\'esultat.

Il
y a les sous-modules qui v\'erifient $Jac_{B, \cdots,
-(A-2)}\neq 0$ mais ceux-l\`a donnent 0 quand on leur applique
$Jac_{x\in {\cal T}_{A-2,B,D}}$ d'apr\`es ce que l'on
vient de voir (en rempla\c cant encore $(A,B)$ par
$(A-2,B)$;

et il y a les termes dits compl\'ementaires de \ref{enonce},
c'est-\`a-dire les repr\'esentations $$\pi_{\eta'}:=\pi(\psi',\epsilon',\cup_{C\in [B, A-2]}(\rho,C,C,\zeta=+,(-1)^C\eta'))$$ o\`u $\lambda'=\pm$ v\'erifie $\times_{C\in [B,A-2]}((-1)^C\eta')=\epsilon_{0}$. Il faut distinguer suivant les valeurs de $\eta'$ possibles:

pour $\eta'$ tel que $(-1)^B\eta' \neq \epsilon(\rho,B-1,B-1,\zeta=+)$; dans ce cas $\pi_{\eta'}$ est spontan\'ement de la forme $\pi_{D,\lambda}$ mais pour la valeur $D=A-2$ et $\lambda=\epsilon(\rho,0,0,\zeta=+)$. Cette repr\'esentation intervient effectivement avec multiplicit\'e 1 mais aucune autre valeur de $D$ n'appara\^{\i}t. La condition pour que cela puisse se produire est donc:
$$
\times_{C\in [B,A-2]}((-1)^C(-1)^{B-1}\epsilon(\rho,B-1,B-1,\epsilon=+))=\epsilon_{0}.\eqno(1)
$$

pour $\eta'$ tel que $(-1)^B\eta'=\epsilon(\rho,B-1,B-1,\zeta=+)$; dans ce cas le $D$ et le $\lambda$ associ\'e a cette repr\'esentation v\'erifie $D=A-2-2inf(B,A-1-B)$. En revenant \`a la notation $t$ telle que $A-D=2t$, on va avoir $t=1+inf(B,A-B-1)$. Mais on sait a priori que $t\leq B$ (dans le tableau ${\cal T}_{D}$ il y a au plus $2B$ lignes et au moins $2t$ lignes). Ainsi n\'ecessairement $A-B+1 <B$. D'o\`u encore $D=A-2-2(A-B-1)=B- (A-B)$ et $D+1=B-(A-B-1)\leq B$ contrairement \`a l'hypoth\`ese de l'\'enonc\'e. Ces repr\'esentations ne nous int\'eressent donc pas.

\nl
Par r\'ecurrence, on conna\^{\i}t aussi la structure de $\pi(\psi'',\epsilon'',(\rho,A-2,B-1,\eta))$ pour $\eta$ v\'efiant
$\eta\epsilon_{0}=\epsilon(\rho,B-1,B-1,\zeta=+)$; comme ci-dessus, nous n'avons \`a nous pr\'eoccuper que des repr\'esentations dites compl\'ementaires dans \ref{enonce}. Il s'agit des repr\'esentations $\pi_{\eta''}:=\pi(\psi'',\epsilon'', \cup_{C\in [B-1,A-2]}(\rho,C,C,\zeta=+,(-1)^C\eta''))$, o\`u $\eta''$ v\'erifient
$\times_{C\in [B-1,A-2]}((-1)^C\eta'')=\eta$. Il faut distinguer suivant les valeurs de $\eta''$ possibles:

pour $\eta''$ tel que $(-1)^{B-1}\eta''\neq \epsilon(\rho,B-2,B-2,\zeta=+)$ (on suppose ici que $B\geq 2$ et on verra les ''petits cas'' plus loin); dans ce cas $\pi_{\eta''}$ est spontan\'ement de la forme $\pi(D,\lambda)$ mais pour la valeur $D=A-2$ et $\lambda=\epsilon(\rho,0,0,\zeta=+)$. Le probl\`eme ici est donc d'\'eliminer cette repr\'esentation avec son analogue trouv\'ee ci-dessus. D\'ej\`a elle n'intervient ici que si $$\prod_{C\in [B-1,A-2]}((-1)^C(-1)^{B}\epsilon(\rho,B-2,B-2,\zeta=+))=\epsilon_{0}\epsilon(\rho,B-1,B-1,\zeta=+)$$
On rappelle que $\epsilon(\rho,B-2,B-2,\zeta=+)=-\epsilon(\rho,B-1,B-1,\zeta=+)$ par hypoth\`ese. Cette condition se r\'ecrit donc
$$\prod_{C\in [B-1,A-2]}((-1)^C(-1)^{B-1}\epsilon(\rho,B-1,B-1,\zeta=+))=\epsilon_{0}\epsilon(\rho,B-1,B-1,\zeta=+)$$ ou encore, en changeant l'intervalle du produit:
$$
\prod_{C\in [B,A-2]}((-1)^C(-1)^{B-1}\epsilon(\rho,B-1,B-1,\zeta=+))=$$
$$
(-1)^{B-1}(-1)^{B-1}\epsilon(\rho,B-1,B-1,\zeta=+)\epsilon_{0}\epsilon(\rho,B-1,B-1,\zeta=+)=\epsilon_{0}$$
Ceci est donc exactement la  condition (1).
La multiplicit\'e avec laquelle elle intervient dans $Jac_{x\in {\cal T}_{D}}\pi(\psi,\epsilon)$ est donc directement 0 si (1) n'est pas satisfait et sinon,  il faut se rappeler les signes, cf. ci-dessus et que $D=A-2$,
$$
(-1)^{A-B-1} + (-1)^{[(A-B+1)/2]}\eta^{(A-B+1)}\epsilon_{0}^{A-B}.\eqno(2)$$
On remplace dans (1) $\epsilon(\rho,B-1,B-1,\zeta=+)$ par $\eta\epsilon_{0}$. La condition (1) se r\'ecrit
$$(\eta \epsilon_{0})^{A-2-B+1}\prod_{C'\in [1,(A-B-1)]}(-1)^{C'}=\epsilon_{0}$$
ou encore $\eta^{A-B+1}\epsilon_{0}^{A-B}(-1)^{[A-B]/2}=1$. Or $(1)^{[(A-B+1)/2]+[(A-B)/2}$ vaut $1$ si $A-B$ est pair et $-1$ si $A-B$ est impair. En d'autres termes cela vaut $(-1)^{A-B}$. La multiplicit\'e cherch\'ee en (2) est donc 0.

Il reste \`a regarder le cas o\`u $\eta''$ est tel que $(-1)^{B-1}\eta''=\epsilon(\rho,B-2,B-2,\zeta=+)$. Dans ce cas, le $D$ et le $\lambda$ associ\'e a cette repr\'esentation v\'erifie $D=A-2-2inf(B-1,A-B)$ et $\lambda$ est d\'etermin\'e par la valeur de $D$. C'est-\`a-dire exactement $\lambda=\epsilon(\rho,0,0,\zeta=+)$ si $A-B< B-1$ et l'oppos\'e si $A-B \geq B-1$; il suffit maintenant de remarquer qu'ici $D=A-2inf(B,A-B+1$ et on trouve la contribution d'une des repr\'esentations  compl\'ementaires de $\pi(\psi,\epsilon)$ (celle qui v\'erifie $\eta(-1)^B=\epsilon(\rho,B-1,B-1,\zeta=+)$) dont on a d\'ej\`a prouv\'e qu'elle intervient avec multiplicit\'e 1 exactement (\ref{1etape}).

Il nous reste \`a regarder le cas $B=1$; cela se fait comme ci-dessus mais en plus simple car il n'y a qu'un type de repr\'esentations compl\'ementaire pour une repr\'esentations $\pi(\psi',\epsilon', (\rho, A',0,\zeta=+,\epsilon'_{0}))$ qui correspond \`a $D=A'$. Cela termine la preuve.

\subsubsection{Fin de la preuve de ce cas particulier}
\bf Proposition: \sl sous les hypoth\`eses de ce paragraphe, 
le th\'eor\`eme \ref{enonce} est vrai.
\rm
\nl
Fixons $\tilde{\pi}$ un constituant de $\pi(\psi,\epsilon)$;
on note $D,\lambda$ les nombres qu'on lui a associ\'es en
\ref{2construction} et \ref{unicite}.

On a d\'ej\`a montr\'e que $\tilde{\pi}$ satisfaisait aux
propri\'et\'es de \ref{enonce} sauf \'eventuellement si
$D+1\geq B$ et $\tilde{\pi}\simeq \pi_{{\cal
T}_D,D,\lambda}$. On veut \'eliminer ce cas sauf si
$t=B$ et $\lambda=-\epsilon(\rho,0,0,\zeta=+)$ qui est l'une des
repr\'esentations compl\'ementaires de \ref{enonce}. On
suppose donc encore que le couple $(D,\lambda)$ ne
v\'erifie pas cette derni\`ere condition et on a alors d\'emontr\'e que
la multiplicit\'e de
$\pi_{D,\lambda}$ dans
$Jac_{x\in {\cal T}_D}\pi(\psi,\epsilon)$ est 0. Il
faut en d\'eduire le r\'esultat cherch\'e; pour cela il
suffit de montrer que si $\tilde{\pi}'$ est un consituant
irr\'eductible de $\pi(\psi,\epsilon)$ tel que $Jac_{x\in
{\cal T}_D}\tilde{\pi}'$ contient $\pi_{D,\lambda}$
comme sous-quotient alors $\tilde{\pi}'\simeq \pi_{{\cal
T},D,\lambda}$. Fixons donc $\tilde{\pi}'$ comme
ci-dessus dont le module de Jacquet contient comme
sous-quotient
$\otimes_{x\in {\cal T}_D}\rho\vert\,\vert^x\otimes
\pi_{D,\lambda}$. On sait d'apr\`es \ref{descriptiontableau} et \ref{unicite} que soit $\tilde{\pi}'$ est un
sous-module de $$<\rho\vert\,\vert^{B}, \cdots,
\rho\vert\,\vert^{-A }>\times \pi(\psi',\epsilon',
(\rho,A-1,B+1,\epsilon_{0})))\eqno(1)$$
soit $\tilde{\pi}'=\tilde{\pi}_{{\cal T},\delta,\zeta}$, ce
que  nous voulons. Mais le d\'ebut de la preuve de (v) du lemme pr\'ec\'edent \'elimine la premi\`ere possibilit\'e, (1). D'o\`u le corollaire.

\subsection{Extension}
Ici on \'etend la preuve ci-dessus  au cas o\`u cas tr\`es voisin o\`u on suppose:
\nl
$(\rho,A,B,\zeta=+,\epsilon_{0})$ est l'unique bloc de Jordan  v\'erifiant $A-B>0$

\nl
il existe $C_{0} \leq B$ et un signe $\epsilon$ tel que $(\rho,C,C,\zeta_{C},(-1)^C\epsilon) \in Jord(\psi,\epsilon)$ pour tout $C<C_{0}$;

on peut toujours supposer $\zeta_{C}=\zeta$ (cf \ref{propriete}) et ce sont donc les hypoth\`eses de \ref{hypotheses}.

On fixe encore $\tilde{\pi}$ une composante irr\'eductible
de $\pi(\psi,\epsilon)$, on suppose encore que $\zeta=+$ et on montre comme en
\ref{ordremaximal}, que soit
$$Jac_{B, \cdots, -A}\tilde{\pi}\neq 0$$ soit 
$Jac_{x\in {\cal T}}\tilde{\pi} \neq 0$, o\`u
$$
{\cal T}:=
\begin{matrix}
&B &\cdots &A\\
&\vdots &\vdots &\vdots\\
&C_{0}+1 &\cdots &C_{0}+A-B+1
\end{matrix}
$$
Or, posons $s:=B-C_{0}$; en particulier le
tableau ci-dessus a $s$ lignes. On note $\sigma_{T}$ la repr\'esentation associ\'ee par Zelevinsky \`a l'ensemble des segments qui constiuent les lignes de ${\cal T}$.
On montre d'abord l'assertion suivante: soit $\tilde{\pi}$ un constituant irr\'eductible de $\pi(\psi,\epsilon)$ tel que $Jac_{x\in {\cal T}}\tilde{\pi}\neq 0$; alors $Jac_{x\in {\cal T}}\tilde{\pi}$ est irr\'eductible, notons le $X$. Et $\tilde{\pi}$ est l'unique sous-repr\'esentation irr\'eductible de l'induite $\sigma_{T}\times X$. 

On le montre en fixant un sous-quotient irr\'eductible $Y$ de $Jac_{x\in {\cal T}}$ tel que $\tilde{\pi}$ soit un sous-module de l'induite $\sigma_{T}\times Y$; l'existence d'un tel $Y$ est facile par r\'eciprocit\'e de Frobenius et utilisant bien s\^ur le fait que $Jac_{x}\tilde{\pi}\neq 0$ avec $\vert x\vert \leq A$ n\'ecessite $x=B$. Le point est de montrer que $Jac_{y}Y=0$ pour tout $y\in [C_{0}+1,A]$. Or s'il n'en est pas ainsi, il existerait une repr\'esentation $Y'$ et une inclusion 
$$
\tilde{\pi}\hookrightarrow \sigma_{T}\times \rho\vert\,\vert^y \times Y'.
$$
Cette inclusion se factorise par un sous-quotient irr\'eductible de $\sigma_{T}\times \rho\vert\,\vert^y$, not\'e $\sigma'$. Mais il r\'esulte des classifications de Zelevinski qu'un tel $\sigma'$ va avoir la propri\'et\'e qu'il existe $y'< B$, $y'>0$ tel que $Jac_{y'}\sigma'\neq 0$. Ceci donne une contradiction $Jac_{y'}\tilde{\pi} \neq 0$. Cela prouve que $Jac_{x\in {\cal T}}\sigma_{T}\times Y=Y$ d'o\`u l'irr\'eductibilit\'e de $Jac_{x\in {\cal T}}$ et aussi l'unicit\'e du sous-module irr\'eductible de cette induite.

Cette assertion montre que l'application $Jac_{x\in {\cal T}}$ induite une bijection entre l'ensemble des constituants irr\'eductibles de $\pi(\psi,\epsilon)$ tel que ce $Jac$ soit non nul et l'ensemble des sous-quotients irr\'eductibles de $Jac_{x\in {\cal T}}\pi(\psi,\epsilon)$.

On calcule facilement $Jac_{x\in {\cal T}}\pi(\psi,\epsilon)$ c'est le calcul de \ref{lecastroue} et on trouve:
$$
Jac_{x\in {\cal T}} \pi(\psi,\epsilon)=\pi(\psi',\epsilon',(\rho,A-s,B-s,\zeta=+,\epsilon_{0})).\eqno(1)
$$
Le membre de droite de (1) est pr\'ecis\'ement le cas que l'on vient de traiter.
Ses constituants irr\'eductibles sont donc de 2 types, les $X_{\eta}$ appel\'es termes compl\'ementaires et ceux qui v\'erifient $Jac_{B-s, \cdots, -A+s}\neq 0$. Soit d'abord $\tilde{\pi}'$ un tel constituant v\'erifiant $Jac_{B-s, \cdots, -A+s}\tilde{\pi'}\neq 0$; on a d\'ej\`a vu qu'alors $Jac_{B-s,\cdots, -A+s}\tilde{\pi}$ est un constituant irr\'eductible de $\pi(\psi',\epsilon', A-s-1, B-s+1,\zeta=+,\epsilon_{0}))$ que l'on note $\tilde{\pi}''$. Soit $\tilde{\pi}$ le constituant de $\pi(\psi,\epsilon)$ qui  correspond, c'est -\`a-dire:
$$
\tilde{\pi}\hookrightarrow \sigma_{T}\times <\rho\vert\,\vert^{B-s}, \cdots, \rho\vert\,\vert^{-A+s}> \times \tilde{\pi}''.
$$
On note ${\cal T}'$ le tableau ${\cal T}$ auquel on a enlev\'e la premi\`ere et la derni\`ere colonne et on note $\sigma_{{\cal T}'}$ la repr\'esentation associ\'ee. Alors on a:
$$
\sigma_{T}\hookrightarrow <\rho\vert\,\vert^{B},\cdots, \rho\vert\,\vert^{B-s+1}>\times \sigma_{T'} \times <\rho\vert\,\vert^{A}, \cdots, \rho\vert\,\vert^{A-s+1}>.
$$
On utilise l'isomorphisme:
$$\sigma_{T'}\times <\rho\vert\,\vert^{A}, \cdots, \rho\vert\,\vert^{A-s+1}> \times 
<\rho\vert\,\vert^{B-s}, \cdots, \rho\vert\,\vert^{-A+s}>$$
$$
\simeq
<\rho\vert\,\vert^{B-s}, \cdots, \rho\vert\,\vert^{-A+s}> \times \sigma_{T'} \times 
<\rho\vert\,\vert^{A}, \cdots, \rho\vert\,\vert^{A-s+1}>
$$
qui vient des r\'esultats standard de Zelevinsky. On a besoin de l'isomorphisme:
$$
<\rho\vert\,\vert^{A}, \cdots, \rho\vert\,\vert^{A-s+1}>\times \tilde{\pi}'' \simeq 
<\rho\vert\,\vert^{-A+s-1}, \cdots, \rho\vert\,\vert^{-A}>\times \tilde{\pi}'',
$$
qui vient du fait que $\tilde{\pi}''$ est un constituant de $\pi(\psi',\epsilon', (\rho, A-s-1,B-s+1,\zeta=+,\epsilon_{0})$:  $\rho\vert\,\vert^x \times \tilde{\pi}''$ est donc irr\'eductible pour tout $x\in [A-s+1,A]$ puisqu'un tel $x$ n'est pas de la forme $C+1$ avec $C\in [A',B']$ o\`u $(\rho,A',B',\zeta') $ un bloc de Jordan pour $\psi' \cup (\rho, A-s-A,B-s+1,\zeta=+)$ (cf. \ref{rappel} que l'on \'etend facilement \`a notre cas, en utilisant \ref{explicite} par exemple).

On a encore $\sigma_{{\cal T}'}\times <\rho\vert\,\vert^{-A+s-1}, \cdots, \rho\vert\,\vert^{-A}> \simeq
<\rho\vert\,\vert^{-A+s-1}, \cdots, \rho\vert\,\vert^{-A}> \times \sigma_{{\cal T}'}$.
En remettant tous cela ensemble on obtient $Jac_{B, \cdots -A}\tilde{\pi}\neq 0$.

\nl

Il ne nous reste donc plus qu'\`a consid\'erer les termes dits compl\'ementaires qui interviennent dans le terme de droite de (1). Fixons donc $\eta'$ tel que $\times_{C\in [B-s,A-s]}((-1)^C\eta')=\epsilon_{0}$ et posons $$X'_{\eta'}=\pi(\psi',\epsilon',\cup_{C\in [B-s,A-s]}(\rho,C,C,\zeta=+,(-1)^C\eta')).$$ Et on note
$X_{\eta'}$ l'unique repr\'esentation irr\'eductible telle que:
$$
X_{\eta'} \hookrightarrow \sigma_{{\cal T}}\times X'_{\eta'}.
$$
Une application directe de \ref{rappel} montre que $X_{\eta'}=\pi(\psi',\epsilon',\cup_{C\in [B,A]}(\rho,C,C,\zeta=+,(-1)^{C-s}\eta'))$. En posant $\eta=(-1)^s\eta'$ on obtient exactement la repr\'esentation compl\'ementaire pour $\pi(\psi,\epsilon)$ correspondant \`a cet $\eta$; la relation que doit satisfaire $\eta$ est \'evidemment satisfaite gr\^ace \`a celle satisfaite par $\eta'$. Cela termine la preuve.
\subsection{R\'eduction}
\subsubsection{Inversion des
socles\label{echangedessocles}}
Avant de pouvoir faire les r\'eductions, on a besoin du
lemme technique ci-dessous. Il emploie les notations
suivantes: $[x,y]$ est un segment croissant ou d\'ecroissant
et 
$\delta:=<\rho\vert\,\vert^x, \cdots, \rho\vert\,\vert^y>$
est soit un module de Speh soit une s\'erie discr\`ete. De
m\^eme $[x',y']$ est un segment croissant ou d\'ecroissant
et $\delta':=<\rho'\vert\,\vert^{x'}, \cdots,
\rho'\vert\,\vert^{y'}>$. On suppose que $-y\notin [x,y]$ et
que $-y'\notin [x',y']$. On suppose aussi $[x,y] \subset
[x',y']$ et $x',y'\notin [x,y]$ ce qui ne suppose pas que
$[x,y]$ et
$[x',y']$ ait la m\^eme propri\'et\'e de croissance. 

Soit $X$ une repr\'esentation semi-simple et on suppose que
pour tout $z\in [x,y]$ et  pour tout $z'\in [x',y']$ on ait
$Jac_{z,\cdots, y}X=0$ et $Jac_{z',\cdots,y'}X=0$. Alors
\nl
\bf Lemme: \sl On a l'\'egalit\'e des socles:
$
<\delta, <\delta',X>>=<\delta',<\delta,X>>.
$
\rm
\nl
Pour d\'emontrer cela, on peut supposer comme nous le ferons
que $X$ est irr\'eductible. On remarque que l'induite
$\delta\times \delta'$ du GL convenable est irr\'eductible
gr\^ace \`a l'hypoth\`ese $[x,y]\subset [x',y']$. On va
d\'emontrer que l'induite
$\delta'\times \delta\times X$ a un unique sous-module
irr\'eductible; cela suffira car  tout sous-module
irr\'eductible de
$\delta'\times <\delta,X>$ co\"{\i}ncidera avec ce
sous-module qui vaudra donc $<\delta',<\delta,X>>$. Mais
$\delta\times
\delta'\times X
\simeq
\delta'\times \delta\times X$. La premi\`ere induite a donc
elle aussi un unique sous-module irr\'eductible, ce sous-module
co\"{\i}ncide n\'ecessairement  avec
$<\delta,<\delta',X>$. D'o\`u l'\'egalit\'e du lemme.
Montrons donc l'unicit\'e. Par r\'eciprocit\'e de Frobenius,
il suffit de montrer que
$$Jac_{x',\cdots, y',x, \cdots, y}\delta'\times \delta\times
X=X.$$

On v\'erifie d'abord que
$Jac_{x',\cdots,y'}\delta'\times \delta\times
X=\delta\times X$. Il est clair que $\delta\times X$ est un
constituant du module de Jacquet cherch\'e et si ce n'est
pas le seul, il y a un d\'ecoupage de $[x',y']$ en 3
ensembles $E_1\cup E_2\cup E_3$ ordonn\'es par l'ordre induit
du segment tel que
$Jac_{z\in E_1}
\delta'\neq 0$, $Jac_{z\in E_2}(\delta')^*\neq 0$ et
$Jac_{z\in E_3} \delta\times X\neq 0$. Or $y'$ ne peut
\^etre ni dans $E_2$ ni dans $E_3$ d'apr\`es les
hypoth\`eses faites. Donc $y'\in E_1$. Mais la non nullit\'e
de $Jac_{z\in E_1}\delta'$ avec $y'\in E_1$ force
$E_1=[x',y']$. Il ne reste plus qu'\`a calculer $Jac_{x,
\cdots, y}\delta\times X$, ce qui se fait de fa\c con
totalement analogue pour trouver $X$.
Cela termine la preuve.
\subsubsection{Un r\'esultat technique}
Soit $Y'$ une repr\'esentation irr\'eductible de $G(n')$,
c'est-\`a-dire d'un groupe de m\^eme type que $G$ et de rang
$n'$. On fixe
$\rho$ et des demi-entiers
$a,b,a',b'$. On suppose que
$a-b$ et $a'-b'$ sont des entiers relatifs et on pose
$\delta:=<\rho\vert\,\vert^a,\cdots, \rho\vert\,\vert^b>$ et
$\delta':=<\rho\vert\,\vert^{a'},\cdots,\rho\vert\,
\vert^{b'}>$. On a en vue un r\'esultat du genre: soit $Y$
un sous-module irr\'eductible de $\delta'\times Y'$, alors
il existe une bijection (naturelle) entre les sous-quotients
irr\'eductibles de $\delta\times Y$ et les sous-quotients
irr\'eductibles de $\delta\times Y'$; la bijection \'etant
donn\'ee par $Jac_{a',\cdots, b'}$. Un tel
r\'esultat est tout \`a fait faux en g\'en\'eral, la
condition minimum pour qu'il puisse \^etre vrai  est que
les induites dans $GL(d_\rho(a+a'-b-b'+2))$, $\delta\times
\delta'$ et $\delta\times (\delta')^*)$ sont
irr\'eductibles. Nous ferons
donc cette
hypoth\`ese  dans tout ce paragraphe
et nous ferons 
 aussi comme hypoth\`ese dans
tout ce paragraphe que  pour tout $x\in [a',b']$,
$Jac_{x,\cdots, b'}Y'=0$ et que soit $\{a,-b\}\cap
[a',b']=\emptyset$ soit $\{b',-b'\}\cap [a,b]=\emptyset$.

\subsubsection{Calcul \'el\'ementaire \label{standard}}
On va constamment utiliser de proche en proche le calcul
\'el\'ementaire suivant: soit $X$ une repr\'esentation (non
n\'ecessairement irr\'eductible) de $G$ et soit
$\alpha,\beta$ des demi-entiers tels que $\alpha-\beta\in
\mathbb{Z}$. Soit encore $y$ un demi-entier. Alors
$Jac_y(<\rho\vert\,\vert^\alpha,\cdots,\rho\vert\,
\vert^\beta>\times X)$ est, dans le groupe de Grothendieck,
la somme d'au plus  3 repr\'esentations:

$
<\rho\vert\,\vert^{\alpha+\zeta 1},\cdots
,\rho\vert\,\vert^\beta>\times X$ si $\alpha=y$ et alors
$\zeta$ est le signe de $\alpha-\beta$

$
<\rho\vert\,\vert^\alpha,\cdots,\rho\vert\,\vert^{\beta
-\zeta 1}\times X$ si $-\beta=y$ et $\zeta$ est alors comme
ci-dessus

$
<\rho\vert\,\vert^\alpha,\cdots,
\rho\vert\,\vert^\beta>\times Jac_y X$.

\subsubsection{Cons\'equences des hypoth\`eses}

Montrons tout de suite la cons\'equence des hypoth\`eses
faites ci-dessus:

\sl pour tout $x\in [a',b']$, $Jac_{x,\cdots,
b'}(\delta\times Y')=0$. 

\rm Supposons d'abord que $
\{a,-b\}\cap
[a',b']=\emptyset$  alors $Jac_{x,\cdots, b'}(\delta\times
Y')\simeq \delta\times Jac_{x,\cdots,b'} Y'$ d'apr\`es
\ref{standard}, de proche en proche seul le 3e cas est
possible.

Sous l'autre hypoth\`ese,
$\{b',-b'\}\cap [a,b]=\emptyset$, on n'a pas l'isomorphisme
pr\'ec\'edent mais
$Jac_{x,\cdots,b'}(\delta\times Y')\neq 0$ n\'ecessite qu'il
existe un sous-ensemble $E \subset [x,b']$ (qui peut a priori \^etre vide) tel que
$Jac_{x'\in E}Y'\neq 0$, et une d\'ecomposition de
$[x,b']-E$ en deux sous-ensembles $E_+$ et $E_-$ tels que
$Jac_{x_+\in E_+}\delta\neq 0$ et $Jac_{x_-\in
E_-}\delta^*\neq 0$. Cela entra\^{\i}ne d\'ej\`a que $E_+$
est de la forme $[a, a_1]$ et que $E_-$ est de la forme
$[-b,b_1]$. Comme par hypoth\`ese ni $b'$ ni $-b'$ ne sont
dans $[a,b]$, $E$ contient $b'$ ce qui est impossible par hypoth\`ese et
prouve l'assertion
\nl
Montrons encore:

 si
$-b'\notin [a',b']$, $Jac_{a',\cdots,b'}(\delta\times
\delta'\times Y')  =\delta\times Y'$

si $-b'\in [a',b']$ mais $[a',-b'[\cap \{a,-b\}=\emptyset$
et $\forall y \in [a',-b'[$, $Jac_yY'=0$; alors
$Jac_{a',\cdots, b'}(\delta\times \delta' \times Y'$ est
exactement la somme de 2 copies de $\delta\times Y'$.

On raisonne comme dans le 2e cas
ci-dessus en remarquant que $\delta\times \delta'\times
Y'\simeq \delta'\times \delta\times Y'$: le calcul du
module de Jacquet se fait en consid\'erant les d\'ecoupages
de
$[a',b']$ en 3 sous-ensembles, $F_+,F_-,F$ tels que
$Jac_{y_+\in F_+}\delta' \neq 0$, $Jac_{y_-\in
F_-}(\delta')^*\neq 0$ et
$Jac_{y\in F}\delta\times Y'\neq 0$. Ainsi $F_+$ est de la
forme $[a',a_1]$, $F_-$ est de la forme $[-b',b_1]$ et $F$
est le compl\'ementaire. Supposons d'abord que
$-b' \notin [a',b']$. Dans ce cas $F_-$ est vide et $F$ est
un segment d'extr\^emit\'e $b'$ ou est vide. Or on sait que
si $F$ est non vide
$Jac_{x\in F}(\delta\times Y')=0$. Donc $F$ est aussi vide
et il n'y a qu'un choix de d\'ecoupage. D'o\`u le r\'esultat.

On suppose maintenant que $-b'\in [a',b']$; alors
donc
$F$ vaut $(]a_1,-b'[\cup ]b_1,b'[)$. D'abord on
montre que dans le d\'ecoupage ci-dessus, si $F_-$ est non
vide, n\'ecessairement
$b_1=b'$. En effet, 
$Jac_{x\in ]a_1,-b'[\cup ]b_1,b']}X \simeq Jac_{x\in
]b_1,b']\cup ]a_1,-b'[}X $, pour toute repr\'esentation $X$, 
l'isomorphisme venant de ce que pour tout
\'el\'ement
$x_1$ de $]b_1,b']$ et tout \'el\'ement $x_2$ de $]a_1,
-b'[$, $x_1 \neq x_2\pm 1$. Ainsi, il faut en particulier
$Jac_{x\in ]b_1,b']} (\delta\times Y')\neq 0$ ce qui
entra\^{\i}ne que $b_1=b'$ si $F_-$ est non vide. Continuons
de supposer que $F_-$ est non vide et montrons que
$]a_1,-b'[$ est vide. Notons $y$ le premier \'el\'ement de
cet intervalle s'il est non vide; on doit avoir
$Jac_y(\delta\times Y')\neq 0$. On a \'enum\'erer les
situation o\`u cela pouvait se produire et les hypoth\`eses
ont \'et\'e mises pour pr\'ecis\'ement \'eliminer ces cas.
Ainsi, il y a au plus un d\'ecoupage possible avec $F_-\neq
\emptyset$. Ce d\'ecoupage fonctionne d'ailleurs tr\`es
bien. Il n'y a aussi qu'un d\'ecoupage avec $F_-=\emptyset$,
c'est celui o\`u $F_+=[a',b']$ (sinon il faudrait $Jac_{y\in
]a_1,b']}(\delta\times Y')\neq 0$). Ensuite le r\'esultat
cherch\'e s'en d\'eduit.

\nl
On suppose ici que $-b'\notin [a',b']$. Comme cas
particulier (o\`u $[a,b]$ est vide) de ce que l'on a vu
ci-dessus, on sait que $Jac_{x\in [a',b']}\delta'\times
Y'\simeq Y'$ est irr\'eductible. En particulier l'induite
$\delta'\times Y'$ a un unique sous-module irr\'eductible,
que l'on note
$Y$ et $Y$ intervient avec multiplicit\'e 1 comme
sous-quotient de cette induite.

\subsubsection{Enonc\'e du r\'esultat technique et d\'emonstration\label{technique}}
{\bf Lemme}: \sl On fixe $a,b,a',b'$ comme ci-dessus en
supposant que $-b'\notin [a',b']$ et soient
$Y',Y$ comme ci-dessus  L'application qui \`a un
sous-quotient irr\'eductible $\theta$ de $\delta\times Y$
associe
$Jac_{a',\cdots, b'}\theta$ d\'efinit une bijection de
l'ensemble des sous-quotients irr\'eductibles de
$\delta\times Y$ dans l'ensemble des sous-quotients
irr\'eductibles de $\delta\times Y'$.\rm
\nl
D'abord on consid\`ere les morphismes d'entrelacement
(standard) entre induites d\'ependant d'un param\`etre $s\in
{\mathbb C}$,
$$M(s):
(\delta')^*\vert\,\vert^s\times \delta\times Y' \rightarrow
\delta'\vert\,\vert^{-s}\times \delta\times Y';
$$
$$M'(s): (\delta')^*\vert\,\vert^s\times Y' \rightarrow
\delta'\vert\,\vert^{-s}\times Y'.$$
On normalise $M'(s)$ de sorte que l'op\'erateur normalis\'e,
$N'(s)$, soit holomorphe non nul en $s=0$. A ce moment l\`a,
en utilisant les normalisations de Langlands pour les groupes
lin\'eaires, on normalise $M(s)$ de fa\c con \`a ce que
$M(s)$ s'\'ecrive comme compos\'e de l'op\'erateur,
$N_{(\delta')^*\times \delta}(s)$ qui
\'echange $(\delta')^*\vert\,\vert^s \times \delta$ en
$\delta\times (\delta')^*\vert\,\vert^s$ puis $N'(s)$ puis
l'op\'erateur $N_{\delta\times \delta'}(-s)$ qui \'echange
$\delta\times
\delta'\vert\,\vert^{-s}$ en $\delta'\vert\,\vert^{-s}\times
\delta$. Les op\'erateurs dans les groupes lin\'eaires sont
holomorphes en $s=0$ et bijectifs gr\^ace \`a l'hypoth\`ese
d'irr\'eductibilit\'e faite. 

On v\'erifie que $(\delta')^*\times Y'$ a $Y$ pour unique
quotient irr\'eductible car ces quotients irr\'eductibles
sont les sous-modules irr\'eductibles de l'induite
$\delta'\times Y'$ (on utilise le fait que la repr\'esentation duale d'une repr\'esentation irr\'eductible d'un groupe classique tel que consid\'er\'e ici est la repr\'esentation de d\'epart tordue \'eventuellement par un automorphisme ext\'erieur venant du groupe des similitudes). De plus comme
$Jac_{a',\cdots,b'}\delta'\times Y'$ est irr\'eductible, $Y$
intervient avec multiplicit\'e 1 comme sous-quotient dans
cette induite. La non nullit\'e de $M'(0)$ entra\^{\i}ne
alors que son image est exactement $Y$. Ainsi $N(s)$ est
d\'efini en $s=0$ et a une image qui est isomorphe \`a
$\delta\times Y$. On note $N(0)$ l'op\'erateur ainsi
d\'efini; on n'a pas le choix $N(s)$ est le produit de
l'op\'erateur d'entrelacement standard par une fonction
m\'eromorphe de $s$. On filtre $\delta\times Y'$ par des
sous-repr\'esentations ${ V}_i$ pour $i$ parcourant un
intervalle de $\mathbb{N}$, filtration croissante de tel
sorte que les sous-quotients soient irr\'eductibles. Ainsi,
pour tout $i$, $N(0)$ induit une application de
$(\delta')^*\times V_i$ dans $\delta'\times V_i$ et par
passage au quotient, si l'on note $\theta_i$ le
sous-quotient irr\'eductible de cette filtration au cran
$i$, une application de $(\delta')^*\times \theta_i$ dans
$\delta'\times \theta_i$. On v\'erifie que $\delta'\times
\theta_i$ a un unique sous-module irr\'eductible et qu'il
intervient avec multiplicit\'e 1 comme sous-quotient de
l'induite: en effet, puisque $\theta_i$ est un sous-quotient
de $\delta\times Y'$, pour tout $x\in [a',b']$,
$Jac_{x,\cdots, b'} \theta_i$ est un sous-quotient de
$Jac_{x,\cdots, b'} (\delta\times Y)$ et vaut donc 0. Comme
ici, $-b'\notin [a',b']$, cela entra\^{\i}ne que
$Jac_{a',\cdots, b'} (\delta'\times \theta_i)= \theta_i$.
D'o\`u les assertions et cela entra\^{\i}ne aussi que
$\tilde{\theta}_i$ est l'unique quotient irr\'eductible de
$(\delta')^*\times
\theta_i$. Ainsi l'image de
$N(0) (\delta')^*
\times
\theta_i$ est soit 0 soit l'unique sous-module
irr\'eductible de $\delta'\times \theta_i$ que nous noterons
$\tilde{\theta}_i$. Ainsi l'image de l'application de
d\'epart
$N(0)$, admet une filtration dont les sous-quotients sont
certains des $\tilde{\theta}_i$ pr\'ec\'edemment d\'efinis.
Comme on sait, a priori, que l'image est isomorphe \`a
$\delta\times Y$, tous les $\tilde{\theta}_i$ doivent
intervenir (il suffit par exemple de comparer les modules de
Jacquet
$Jac_{a',\cdots, b'}$). La bijection du lemme est l'inverse
de l'application $\theta_i \mapsto \tilde{\theta}_i$. Cela
termine la preuve du lemme.

\subsubsection{Premi\`ere r\'eduction \label{1reduction}}
On suppose ici que $Jord(\psi)$ contient 2 quadruplets
$(\rho,A,B,\zeta)$ et $(\rho',A',B',\zeta')$. Dans cette partie on d\'emontre \ref{enonce}
par r\'ecurrence sous cette hypoth\`ese. C'est le cas o\`u
$\rho'=\rho$ qui est le plus difficile et pour \'eviter les
fautes de frappe, on suppose donc que $\rho'=\rho$. Il n'y a
aucune difficult\'e dans la d\'emonstration qui va suivre,
uniquement un probl\`eme de notations. On reprend d'abord
les notations, pour $C'\in ]B',A']$, o\`u $\epsilon'_{0}=\epsilon(\rho,A',B',\zeta)$
$$
\delta_{C'}:=<\rho\vert\,\vert^{\zeta' B'}, \cdots, \rho\vert\,\vert^{-\zeta' C'}>.
$$
On note $\psi'',\epsilon''$ le couple qui se d\'eduit de $\psi,\epsilon$ en enlevant \`a la fois $(\rho,A,B,\zeta,\epsilon_{0}$ et $(\rho,A',B',\zeta',\epsilon'_{0}$.  On remarque que l'hypoth\`ese que $\psi\circ \Delta$ est discret assure que si $A>A'$ alors $[B',-A']\subset [-B,B]$ alors que si $A'>A$ l'inverse se produit, $[B,-A]\subset [-B',B']$. On pose $\delta=<\rho\vert\,\vert^{\zeta B}, \cdots, \rho\vert\,\vert^{-\zeta A}>$ et pour $C'$ comme ci-dessus:
$$
\delta'_{C'}:=<\rho\vert\,\vert^{\zeta' (C'+1)}, \cdots, \rho\vert\,\vert^{\zeta' A'}.
$$
On sait gr\^ace \`a \ref{lecastroue} que l'on peut appliquer par r\'ecurrence que l'on a:
$$
Jac_{\zeta'(B'+2), \cdots \zeta' C'}\pi(\psi'',\epsilon'', (\rho,A,B,\zeta,\zeta_{0}), (\rho,A',B'+2,\zeta',\epsilon'_{0}))=$$
$$ <\delta'_{C'}, \pi(\psi'',\epsilon'',(\rho,A,B,\zeta,\epsilon_{0}),(\rho, A'-1,B'+1,\zeta',\epsilon'_{0})).\eqno (1)
$$
La preuve, dans son esprit, est simple: on utilise $(\rho,A',B',\zeta')$ pour donner la d\'efinition de $\pi(\psi,\epsilon)$; cette d\'efinition fait intervenir des $\pi(\tilde{\psi},\tilde{\epsilon})$ pour les quels $(\rho,A,B,\zeta)$ est un bloc de Jordan; on applique alors \ref{enonce} par r\'ecurrence \`a ces repr\'esentations. Puis ensuite on ''fait commuter'' pour  revenir en arri\`ere. Pr\'ecis\'ement:
$$
\pi(\psi,\epsilon)= \oplus_{C'\in ]B',A']}(-1)^{A'-C'}\delta_{C'}\times <\delta'_{C'}, \pi(\psi'',\epsilon'',(\rho,A,B,\zeta,\epsilon_{0}),(\rho, A'-1,B'+1,\zeta',\epsilon'_{0}))$$
$$
\oplus_{\eta'=\pm}(-1)^{[(A'-B'+1)/2]}(\eta')^{A'-B'+1}(\epsilon'_{0})^{A'-B'}
\pi(\psi'',\epsilon'',(\rho,A,B,\zeta,\epsilon_{0}),(\rho,A',B'+1,\zeta',\eta'),(\rho,B',B',\zeta',\eta'\epsilon'_{0})).
$$
On applique donc \ref{enonce} par r\'ecurrence en utilisant $(\rho,A,B,\zeta)$ et $\pi(\psi,\epsilon)$ est donc la somme des termes suivants:
$$
\oplus_{C'\in ]B',A']}(-1)^{A'-C'}\delta_{C'}\times <\delta'_{C'},<\delta,\pi(\psi'',\epsilon'',(\rho,A-1,B+1,\zeta,\epsilon_{0})(\rho,A'-1,B'+1,\zeta',\epsilon'_{0}))>> \eqno(2)
$$
$$
\oplus_{C'\in ]B',A']; \eta=\pm}(-1)^{A'-C'}\delta_{C'}\times <\delta'_{C'}, \pi(\psi'',\epsilon'', \cup_{C\in [B,A]}(\rho,C,C,\zeta,(-1)^{[C]}\eta),(\rho,A'-1,B'+1,\zeta',\epsilon'_{0}))> \eqno(3)
$$
$$
\oplus_{\eta'=\pm}(-1)^{[(A'-B'+1)/2]}(\eta')^{A'-B'+1}(\epsilon'_{0})^{A'-B'} $$
$$<\delta,\pi(\psi'',\epsilon'', (\rho,A-1,B+1,\zeta,\epsilon_{0}), (\rho,A',B'+1,\zeta',\eta'),(\rho,B',B',\zeta',\eta'\epsilon'_{0}))>\eqno(4)
$$
$$
\oplus_{\eta'=\pm,\eta=\pm}(-1)^{[(A'-B'+1)/2]}(\eta')^{A'-B'+1}(\epsilon'_{0})^{A'-B'} $$
$$
\pi(\psi'',\epsilon'', \cup_{C\in [B,A]}(\rho,C,C,\zeta,(-1)^{[C]}\eta),(\rho,A',B'+1,\zeta',\eta'),(\rho,B',B',\zeta',\eta'\epsilon'_{0})). \eqno(5)
$$
Les termes (3) et (5) donnent imm\'ediatement:
$$\oplus_{\eta=\pm}
\pi(\psi'',\epsilon'',\cup_{C\in [B,A]}(\rho,C,C,\zeta,(-1)^{[C]}\eta),(\rho,A',B',\zeta',\epsilon'_{0})).
$$
Le point est donc de montrer que dans (2) et (4) on peut faire sortir $<\delta,?$. On montre que pour $C'\in ]B',A']$:
$$
<\delta'_{C'},<\delta,\pi(\psi'',\epsilon'',(\rho,A-1,B+1,\zeta,\epsilon_{0}),(\rho,A'-1,B'+1,\zeta',\epsilon'_{0}))>>=$$
$$<\delta,<\delta'_{C'},\pi(\psi'',\epsilon'',(\rho,A-1,B+1,\zeta,\epsilon_{0}),(\rho,A'-1,B'+1,\zeta',\epsilon'_{0}))>>.
$$
Cela r\'esulte de \ref{echangedessocles} en tenant compte du fait que les hypoth\`eses de ce lemme sont satisfaites, pour l'inclusion entre les segments, on l'a dit ci-dessus et pour la nullit\'e des modules de Jacquet, on l'a prouv\'e en \ref{techniquepoursocle}. On pose:
$$
Y_{C'}:=<\delta'_{C'},\pi(\psi'',\epsilon'',(\rho,A-1,B+1,\zeta,\epsilon_{0}),(\rho,A'-1,B'+1,\zeta',\epsilon'_{0}))>.$$
On veut encore que l'ensemble des sous-quotients irr\'eductibles  de $\delta_{C'}\times <\delta,Y_{C'}>$ co\"{\i}ncide avec l'ensemble des sous-quotient de la forme  $<\delta,Z>$, o\`u $Z$ est un sous-quotient irr\'eductible de $\delta_{C'}\times Y_{C'}$. C'est l'objet du lemme technique \ref{technique} dont la d\'emonstration est report\'ee \`a la fin du papier;  le couple $\delta,\delta'$ de loc.cite est ici le couple $\delta_{C'},\delta$, donc le segment $[a,b]$ de loc.cite est ici $[\zeta' B',-\zeta'C']$, le segment $[a',b']$ de loc.cite est ici $[\zeta B,-\zeta A]$ et $Y'$ de loc.cite est $Y_{C'}$. La nullit\'e de $Jac_{x, \cdots, -\zeta A}Y_{C'}$ pour tout $x\in [\zeta B,-\zeta A]$ r\'esulte imm\'ediatement de \ref{techniquepoursocle}; les autres hypoth\`eses sur les segments et leurs extr\'emit\'es sont imm\'ediates \`a v\'erifier.

En regroupant maintenant (2) et (4), on obtient directement 
$$
<\delta,\pi(\psi',\epsilon',(\rho,A-1,B+1,\zeta,\epsilon_{0}),(\rho,A',B',\zeta',\epsilon'_{0})>.
$$
Ceci termine la preuve.
\subsection{deuxi\`eme r\'eduction \label{2reduction}}
On suppose donc maintenant que pour tout $(\rho',A',B',\zeta')\in
Jord(\psi)$, $A'=B4$ sauf pour un quadruplet, encore
not\'e $(\rho,A,B,\zeta)$.

Ici on suppose qu'il existe un tel
$C'$ tel que $(\rho,C',C',\zeta',\epsilon'_{0})\in Jord(\psi,\epsilon)$ pour des signes convenables mais que $C'\geq 1$ et $(\rho,C'-1,C'-1) \notin Jord(\psi\circ \Delta)$.

 On note
$\tilde{\psi},\tilde{\epsilon}$ l'analogue de
$\psi,\epsilon$ quand on remplace $C'$ en $C'-1$
sans changer $\zeta',\epsilon'$.  Pour $C\in ]B,A]$, on reprend la notation $X_{C}$ et pour $\eta=\pm$ la notation $X_{\eta}$ de \ref{notations}.
 
 On v\'erifie alors ais\'ement en utilisant
\ref{technique} que
$Jac_{\zeta' C'}$ induit une bijection entre
l'ensemble des sous-quotients irr\'eductibles de $X_C$
pour tout $C\in ]B,A]$ et l'ensemble analogue quand on remplace
$\psi,\epsilon$ par $\tilde{\psi},\tilde{\epsilon}$. Ceci
est aussi vrai pour les repr\'esentations $X_\eta$ avec
$\eta=\pm$ intervenant dans la d\'efinition de
$\pi(\psi,\epsilon)$. Ainsi on v\'erifie que 
$$
\pi(\psi,\epsilon)=<\rho\vert\,
\vert^{\zeta' C'},\pi(\tilde{\psi},
\tilde{\epsilon})>,
$$
donc en particulier $\pi(\psi,\epsilon)$ est une somme de
repr\'esentations irr\'eductibles. On applique \ref{enonce}
\`a $\pi(\tilde{\psi},\tilde{\epsilon})$ (ici encore
$\epsilon_0=\epsilon(\rho,A,B,\zeta)$ et  $\delta=<\rho\vert\,\vert^{\zeta B}, \cdots, \rho\vert\,\vert^{-\zeta A}>$
 d'o\`u $
\pi(\tilde{\psi},\tilde{\epsilon})=$
$$
<\delta,
\pi(\tilde{\psi}',\tilde{\epsilon}',(\rho,A-1,B+1,\epsilon_0))
\eqno(1)
$$
$$\oplus_{\eta=\pm} \pi(\tilde{\psi}',\tilde{\epsilon}',
\cup_{C\in [A,B]}(\rho,C,C,\zeta,(-1)^{[C]}\eta)).\eqno(2)
$$
Pour conclure on utilise l'\'echange des socles
(\ref{echangedessocles}) appliqu\'e \`a $\delta$ comme ci-dessus et $\delta'= \rho\vert\,\vert^{\zeta' C'}$. On obtient donc 
 $$
 (1)= <\delta,<\rho\vert\,\vert^{\zeta' C'},\pi(\tilde{\psi}',\tilde{\epsilon}', (\rho,A-1,B+1,\zeta,\epsilon_{0})>>$$
 Mais $<\rho\vert\,\vert^{\zeta' C'},\pi(\tilde{\psi}',\tilde{\epsilon}', (\rho,A-1,B+1,\zeta,\epsilon_{0})>=\pi(\psi',\epsilon',(\rho, A-1,B+1,\zeta,\epsilon_{0}))$. De m\^eme:
 $$
 <\rho\vert\,\vert^{\zeta' C'}, \pi(\tilde{\psi}',\tilde{\epsilon}',
\cup_{C\in [A,B]}(\rho,C,C,\zeta,(-1)^{[C]}\eta))>=
  \pi( {\psi}',{\epsilon}',
\cup_{C\in [A,B]}(\rho,C,C,\zeta,(-1)^{[C]}\eta))$$
  pour tout $\eta=\pm$ possible. D'o\`u \ref{enonce} sous ces hypoth\`eses.
\subsection{Troisi\`eme  r\'eduction \label{3reduction}}
Il nous reste donc \`a voir le cas o\`u $(\rho,A,B,\zeta,\epsilon_{0})$ est l'unique \'el\'ement de $Jord(\psi,\epsilon)$ tel que $A>B$ et o\`u pour tout demi-entier $C'\in [1, B-1]$, si $(\rho,C',C',\zeta')\in Jord(\psi)$, alors il existe un signe $\zeta''$ tel que $(\rho,C'-1,C'-1,\zeta'')\in Jord(\psi)$. On a r\'egl\'e en \ref{casparticulierprincipal} o\`u $\epsilon$ alterne sur tous ces blocs de Jordan. On a donc \`a r\'egler le cas o\`u il existe $C'$ tel que $(\rho,C',C',\zeta')\in Jord(\psi)$ et $(\rho,C'-1,C'-1,\zeta'')\in Jord(\psi)$, le caract\`ere $\epsilon$ prenant la m\^eme valeur sur ces 2 blocs. En fixant $C'$ minimal avec cette propri\'et\'e; on peut imposer que $\zeta''=\zeta'$  car la d\'efinition de $\pi(\psi,\epsilon)$ ne d\'epend pas de $\zeta''$ car le signe est altern\'e sur les blocs plus petit ou egaux \`a $C'-1$ en commen\c{c}ant par -1 si $C'$ est un demi-entier non entier. On fixe donc $C',\zeta',\epsilon'_{0}$ la valeur de $\epsilon$ sur $(\rho,C',C',\zeta')$. On va ramener ce cas \`a un dernier cas qui sera trait\'e dans le paragraphe suivant.

\nl
On pose $\pi_\pm(\psi,\epsilon):=\pi(\psi,\epsilon)\oplus
\pi(\psi,\epsilon_-)$ o\`u $(\psi,\epsilon_-)$ se d\'eduit
de $\psi,\epsilon$ en changeant simplement $\epsilon$ sur les
2 blocs $(\rho,C',C',\zeta')$ et
$(\rho, C'-1,C'-1,\zeta')$ en
son oppos\'e. On note $\tilde{\psi},\tilde{\epsilon}$ le couple qui se
d\'eduit de $\psi,\epsilon$ en enlevant les 2 blocs $(\rho,C',C',\zeta')$ et
$(\rho, C'-1,C'-1,\zeta')$. On
revient
\`a la d\'efinition de
$\pi(\psi,\epsilon)$ et $\pi(\psi,\epsilon_-)$. On pose ici $\tilde{\delta}:= <\rho\vert\,\vert^{\zeta' C'}, \cdots, \rho\vert\,\vert^{-\zeta'(C'-1)}>$.

Et on montre
en utilisant \ref{rappel} et \ref{technique} que
$\pi_\pm(\psi,\epsilon)$ est l'ensemble des constitutants de
$\tilde{\delta}times
\pi(\tilde{\psi},\tilde{\epsilon})$
v\'erifiant $Jac_{\rho\vert\,\vert^{\zeta' C'}}\neq
0$. 

On applique \ref{enonce} par r\'ecurrence \`a $\pi(\tilde{\psi},\tilde{\epsilon})$ d'o\`u, en posant $\delta:=<\rho\vert\,\vert^{\zeta B}, \cdots \rho\vert\,\vert^{-\zeta A}>$:
$$
\pi(\tilde{\psi},\tilde{\epsilon})=<\delta,\pi(\tilde{\psi}',\tilde{\epsilon}', (\rho, A-1,B+1,\zeta,\epsilon_{0}))> \oplus_{\eta=\pm} \pi(\tilde{\psi}',\tilde{\epsilon}',\cup_{C\in [A,B]}(\rho, C,C,\zeta,(-1)^{[C]}\eta)).
$$
On applique encore \ref{technique} en prenant pour le couple $(\delta,\delta')$ de loc. cite $(\tilde{\delta},\delta)$ et cela nous donne une bijection entre les sous-quotients irr\'eductibles de $\tilde{\delta}\times \pi(\tilde{\psi},\tilde{\epsilon})$ et ceux de 
$\tilde{\delta}\times \pi(\tilde{\psi}',\tilde{\epsilon}', (\rho, A-1,B+1,\zeta,\epsilon_{0}))$ cette bijection \'etant donn\'ee par $<\delta, ?>$. Cette bijection est compatible \`a l'op\'eration $Jac_{\zeta' C'}$. On en d\'eduit alors avec des notations que l'on esp\`ere \'evidente que :
$$
\pi_{\pm}(\psi,\epsilon)=<\delta, \pi_{\pm}(\psi',\epsilon',(\rho, A-1,B+1,\zeta,\epsilon_{0})> 
\eqno(1)
$$
$$\oplus_{\eta}
\pi_{\pm}(\psi',\epsilon', \cup_{C\in [B,A]}(\rho,C,C,\zeta,(-1)^{[C]}\eta)).\eqno(2)
$$
Il faut revenir de $\pi_\pm(\psi,\epsilon)$ \`a
$\pi(\psi,\epsilon)$, c'est-\`a-dire qu'il faut d'abord
d\'emontrer qu'il n'y a pas de simplification entre les
composants irr\'eductibles de $\pi(\psi,\epsilon)$ et ceux
de $\pi(\psi,\epsilon_-)$ et ensuite il faut s\'eparer les
composants respectifs de fa\c con compatible \`a (1) et (2)
ci-dessus. On peut le faire en utilisant les modules de Jacquet dans les 2 cas que l'on va d\'etailler ci-dessous et il restera un cas \`a traiter par d'autres m\'ethodes.

On suppose  ici que $C'\geq 3/2$. Alors le caract\`ere $\epsilon_{-}$ v\'erifie $$\epsilon_{-}(\rho,C'-2,C'-2,\zeta')=\epsilon(\rho, C'-2,C'-2,\zeta')=-\epsilon(\rho,C'-1,C'-1,\zeta')=\epsilon_{-}(\rho,C'-1,C'-1,\zeta').$$
On applique \ref{technique} aux couples $(\delta=<\rho\vert\,\vert^{\zeta B}, \cdots, \rho\vert\,\vert^{-\zeta   C}>, \delta'=\rho\vert\,\vert^{\zeta'(C'-1)})$, pour tout  $C\in ]B,C]$  fix\'e et  \`a la repr\'esentation 
$$
Y=Jac_{\zeta' (B+1), \cdots, \zeta C} \pi(\psi',\epsilon'_{-}, (\rho, A-1,B+1,\zeta,\epsilon_{0})), Y'=Jac_{\zeta'(C'-1)}Y.
$$
Et on montre ainsi que tous les constituants irr\'eductibles de $\pi(\psi,\epsilon_{-})$ v\'erifient $Jac_{\zeta' (C'-1)}\neq 0$. Par contre on a directement que tous les constituants irr\'eductibles de $\pi(\psi,\epsilon)$ v\'erifient $Jac_{\zeta' (C'-1)}=0$ (c'est vrai pour tous les termes de la d\'efinitions de $\pi(\psi,\epsilon)$. Ainsi $Jac_{\zeta' (C'-1)}$ permet de faire la s\'eparation annonc\'ee.

Supposons maintenant que $C'\leq 1$, c'est-\`a-dire $C'=1$.
On peut encore utiliser les modules de Jacquet si
$(\rho,2,2,\zeta'' )\in Jord(\psi)$ avec $\zeta''$ un signe convenable; en effet si $\zeta''=\zeta'$, on utilise
$Jac_{\zeta' 2}$ tandis que si $\zeta''=-\zeta'$, on utilise $Jac_{-\zeta' 2, -\zeta' 1}$, c'est-\`a-dire que dans la d\'emonstration ci-dessus, on remplace $\zeta'(C'-1)$ soit par $\zeta'(C'+1)$ soit par $-\zeta'(C'+1), -\zeta' C'$.

\subsubsection{Dernier cas\label{dernierereduction}}
Il reste exactement un cas qui n'est pas r\'egl\'e par les r\'eductions pr\'ec\'edentes: $(\rho,A,B,\zeta,\epsilon_{0})$ est l'unique \'el\'ement de $Jord(\psi)$ v\'erifiant $A>B$ de plus $B>1$ et pour tout $C<B$ et tout signe $\zeta_{C}$, $(\rho,C,C,\zeta_{C})\notin Jord(\psi)$ sauf exactement si $C=0$ et $C=1$ (pour un bon choix de $\zeta_{C}$. On peut alors prendre $\zeta_{0}=\zeta_{1}=:\zeta'$. Et  $\epsilon$ prend la m\^eme valeur sur $(\rho,C,C,\zeta')$ pour $C=0,1$ et on note $\epsilon'_{0}$ cette valeur.

La m\'ethode est ici du m\^eme ordre que dans
\ref{casparticulierprincipal} mais il y a quelques petits
changements! On fixe $\tilde{\pi}$ un constituant
irr\'eductible de
$\pi(\psi,\epsilon)$.

Pour $D$ un entier  avec $D\leq A$ et pour  $\lambda$ un signe, on reprend la notation
$\tilde{\psi}_D,\tilde{\epsilon}_{D,\lambda}$ de
loc. cite. 
\nl
\bf Lemme: \sl (i)
Il existe $\delta$ et $\zeta$ ainsi qu'un ensemble
totalement ordonn\'e, 
${\cal E}$ d'entiers tous de valeur absolue
inf\'erieure ou
\'egale \`a $A$ tels que:
$$
\tilde{\pi}\hookrightarrow \times_{x\in {\cal E}}
\rho\vert\,\vert^x\times
\pi(\tilde{\psi}_D,\tilde{\epsilon}_{D,\lambda}).
\eqno(1)$$
(ii)Soit ${\cal E}$ comme en (i) et supposons que cet
ensemble contient $-\zeta A$. Alors $Jac_{\zeta B, \cdots,
-\zeta A}\tilde{\pi}\neq 0$.
\nl
\rm
Pour (i) on applique d'abord \ref{1construction} d'o\`u $t_{1}$ et un signe $\lambda_{1}$; ensuite comme $\epsilon$ prend la m\^eme valeur sur $(\rho,0,0,\zeta')$ et $(\rho,1,1,\zeta')$, on peut enlever ces 2 blocs; on se retrouve avec une repr\'esentation \'el\'ementaire mais dont les premiers blocs sont $\cup_{C\in [B, A-2t_{1}]}(\rho, C,C,\zeta,(-1)^{C}\lambda_{1})$. On peut encore appliquer \ref{rappel} pour passer de $\cup_{C\in [B, A-2t_{1}]}$ \`a $\cup_{C}\in [0,A-2t_{1}-B]$. C'est-\`a-dire que $D=A-2t_{1}-B$ convient avec $\lambda=(-1)^B\lambda_{1}$.

 On remplace ensuite sous-quotient en sous-module (quitte \`a changer les signes et l'ordre des \'el\'ements de $\cal E$ en utilisant \ref{rappel} (propri\'et\'e 4).

Pour (ii), l'hypoth\`ese assure qu'il existe $x_0\in {\cal
E}$ tel que $[x_0,-\zeta A]$ soit un segment et $Jac_{x_0,
\cdots, -\zeta A}\tilde{\pi}\neq 0$ (cf \ref{notation}). Ici $x_0$ vaut soit
$\zeta B$ et (ii) est d\'emontr\'e soit $x_0=\zeta' 1$.
Supposons donc que $x_0=\zeta' 1$; pour $C\in ]B,A]$, on
calcul $Jac_{\zeta' 1, \cdots, -\zeta A}X_C$ o\`u:
$$
X_{C}= <\rho\vert\,\vert^{\zeta B}, \cdots, \rho\vert\,\vert^{-\zeta C}>\times Jac_{\zeta (B+1), \cdots, \zeta C}\pi(\psi',\epsilon',(\rho,A,B+2,\zeta,\epsilon_{0})).
$$
On note $\psi'',\epsilon''$ le couple qui se d\'eduit de $\psi',\epsilon'$ en enlevant les 2 blocs $(\rho, C,C,\zeta',\epsilon'_{0})$ pour $C=0,1$ et on a une inclusion (cf. \ref{rappel} et par exemple \ref{explicite} que l'on peut utiliser par r\'ecurrence pour l'\'etendre au cas non \'el\'ementaire)
$$
\pi(\psi',\epsilon',(\rho,A,B+2,\zeta,\epsilon_{0})) \hookrightarrow <\rho \vert\,\vert^{\zeta' 1}, \rho> \times \pi(\psi'',\epsilon'', (\rho, A,B+2,\zeta,\epsilon_{0})).
$$
D'o\`u $X_{C} \hookrightarrow 
<\rho\vert\,\vert^{\zeta B}, \cdots, \rho\vert\,\vert^{-\zeta C}>\times <\rho \vert\,\vert^{\zeta' 1}, \rho> \times 
Jac_{\zeta (B+1), \cdots, \zeta C}\pi(\psi'',\epsilon'',(\rho,A,B+2,\zeta,\epsilon_{0})).$
On calcule $Jac_{\zeta' 1, \cdots, -\zeta A} $ de l'induite ci-dessus. On obtient:

si $\zeta'=-\zeta$, $<\rho\vert\,\vert^{\zeta B}, \cdots, \rho\vert\,\vert^{-\zeta C}>\times \rho \times 
Jac_{-\zeta 2, \cdots , -\zeta A}
Jac_{\zeta (B+1), \cdots, \zeta C}\pi(\psi'',\epsilon'',(\rho,A,B+2,\zeta,\epsilon_{0})).$

si $\zeta'=\zeta$, $<\rho\vert\,\vert^{\zeta B}, \cdots, \rho\vert\,\vert^{-\zeta C}>\times
Jac_{-\zeta 1, \cdots , -\zeta A}
Jac_{\zeta (B+1), \cdots, \zeta C}\pi(\psi'',\epsilon'',(\rho,A,B+2,\zeta,\epsilon_{0})).$

Mais une non nullit\'e de ces termes force, $Jac_{-\zeta y}Jac_{\zeta (B+1), \cdots, \zeta C}\pi(\psi'',\epsilon'',(\rho,A,B+2,\zeta,\epsilon_{0})) \neq 0$ pour $y=2$ ou $1$ d'o\`u encore
$Jac_{-\zeta y} \pi(\psi'',\epsilon'',(\rho,A,B+2,\zeta,\epsilon_{0})) \neq 0$ ce qui est impossible car $-\zeta y \neq \zeta (B+2)$. Cela termine la preuve du lemme.
\nl
 Soit ${\cal E}$, $D,\lambda$ satisfaisant
au (1) du lemme pr\'ec\'edent. Quitte \`a r\'eordonner
${\cal E}$, on \'ecrit ${\cal E}$ comme union de segments
croissants avec les propri\'et\'es de \ref{ordremaximal} et
en notant
$\sigma_{\cal E}$ la repr\'esentation associ\'ee par
Zelevinski \`a cette union, (1) devient une inclusion:
$$
\tilde{\pi}\hookrightarrow \sigma_{\cal E}\times
\pi(\psi_D, \epsilon_{D,\lambda}).\eqno(2)
$$
Ecrivons donc ${\cal E}$ sous la forme $\cup_{i\in [1,v]}
[x_i,y_i]$ o\`u $v$ est un entier convenable; les segment sont  par hypoth\`ese croissants si $\zeta=+$ et d\'ecroissant si $\zeta=-$. Pour fixer les id\'ees, on va donc supposer que $\zeta=+$. On impose alors, comme on en a le droit que $x_1\geq
\cdots \geq x_v$. On peut reprendre les arguments de
\ref{ordremaximal}; s'il existe $i\in [1,v]$ tel
que pour tout
$j<i$,
$x_j\neq x_i-1$, alors $Jac_{x_i}\tilde{\pi}\neq 0$. Mais
ici cela prouve uniquement que pour un tel $i$, $x_i= B$ 
ou
$x_{i}=\zeta' 1$. Ceci s'applique \'evidemment pour $x_1$.

Supposons que ${\cal E}$ comme ci-dessus a aussi la
propri\'et\'e d'avoir un nombre d'\'el\'ements positifs
minimal et qu'il ne contient pas $-A$. On va trouver les formes particuli\`eres que peuvent avoir ${\cal E}$. 
On doit encore avoir ${\cal E}\cup -{\cal E}\cup_{x\in [0,D]}[-x,x]=\cup_{
x\in [B,A]} [-x,x]  \cup [-1,1] \cup \{0\}$. Ceci se r\'ecrit $$
{\cal E}\cup -{\cal E}=\cup_{x\in [D+1,A] }[-x,x] - \cup_{x\in ]1,B[}[-x,x].$$

 On montre comme en **** que sous ces hypoth\`eses,
pour tout $i\in [1,v]$ tel que $y_i>0$ alors il existe $j>i$
avec $y_j=y_i-1$ sauf \'eventuellement pour la plus petite
valeur de $y_i>0$ qui est n\'ecessairement alors
$D+1$.

De plus, on a avec la notation $t_1$ de \ref{1construction},
$$ A-D= 2t_1+B.\eqno(**)$$

Premier cas: supposons que $x_1=\zeta' 1$; alors pour
tout $i\geq 1$, $x_i\neq  B$. Comme on a $Jac_{x_1,x_1}
\tilde{\pi}=0$ gr\^ace
\`a \ref{jacquet}, on a d'apr\`es ce qui pr\'ec\`ede pour
tout
$i\in [1,v[$,
$x_i=x_{i+1}+1$. Comme dans \ref{ordremaximal}, on v\'erifie
alors que n\'ecessairement $y_1 \geq \cdots \geq y_v$. On
r\'ecrit
${\cal E}$ sous forme de tableau comme en loc.cite; toutefois
ici le nombre de colonnes est ici $A- \zeta' 1+1$. On note $w$
le nombre de colonnes et on r\'eutilise la notation
$z_1, \cdots, z_w$ pour les \'el\'ements en bout de chaque
colonne.

Si $\zeta'=-$,  la premi\`ere colonne est $[-1,z_{1}]$ la deuxi\`eme $[0, z_{2}]$ et les autres de la forme $[a,z_{a+2}]$ avec $a\in [1,A]$ et l'\'egalit\'e sur ${\cal E}\cup{\cal -E}$ devient:
$$[1,A] \cup_{i\in [1,w]; z_i\leq 0} \{-z_i\} \cup ]1,B[=
[D+1,A] \cup_{i\in [1,w];z_i>0} \{z_i-1\}.
$$
D'apr\`es ce que l'on a vu $\cup_{i\in [1,w];z_i>0}
\{z_i-1\}$ est un intervalle \'eventuellement vide (c'est
clair si on
\'ecrit
${\cal E}$ sous forme de tableau) dont la plus
grande extr\'emit\'e ne peut \^etre que $D$. Cela
prouve que l'ensemble de droite n'a pas de multiplicit\'e.`Comme $z_{1} \in [-1,-A[$ l'ensemble de gauche a de la multiplicit\'e, ce qui donne une contradiction.

Si $\zeta'=+$, cela se traduit par:
$$
[1, A]\cup_{i\in [1,w]; z_i\leq 0} \{-z_i\}\cup ]1,B[=
[D+1,A] \cup_{i\in [1,w];z_i>0} \{z_i-1\}.
$$
L'argument donn\'e ci-dessus est toujours valable pour
savoir qu'il n'y a pas de multiplicit\'e dans l'ensemble de
droite. Ainsi $B=2$ et $\cup_{i\in [1,w]; z_i\geq 0} \{-z_i\}$ a au
plus 1 terme $\{0\}$. Donc le nombre de lignes du tableau
est au plus 2; or d'apr\`es (**) il est d'au moins 2 lignes,
c'est \`a dire que le tableau a exactement 2 lignes et que
$t_1=0$ puisque $B= 2$.  On conclut alors que $\tilde{\pi}$ est
n\'ecessairement l'une des repr\'esentations
compl\'ementaires de \ref{enonce}.

Deuxi\`eme cas, $x_1= B$: supposons d'abord que $B>2$. La premi\`ere ligne de ${\cal E}$ est $[B,A]$ et on peut calculer le $Jac$ suivant cette ligne; c'est un calcul qui a \'et\'e fait en \ref{lecastroue}, cela revient \`a remplacer dans chaque terme d\'efinissant $\pi(\psi,\epsilon)$, le couple $(A,B)$ par le couple $(A-1,B-1)$. Ainsi on se ram\`ene facilement au cas o\`u $B=2$. L'int\'er\^et est que, sous cette hypoth\`ese,
 $$
 {\cal E}\cup {\cal -E}= [D+1,A].
 $$
Supposons d'abord que ${\cal E}$ peut encore s'\'ecrire sous forme de tableau dont les lignes sont des segments croissants et les colonnes des segments d\'ecroissants; c'est ce dernier point qui n'est pas automatique:
$$
\begin{matrix}
& B &\cdots &\cdots &\cdots &A\\
&\vdots &\vdots &\vdots &\vdots &\vdots\\
&x_{2t_1+B} &\cdots &\cdots &\cdots &D+1\\
&\vdots &\vdots &\vdots&\vdots\\
&x_v &\cdots &y_v
\end{matrix}
$$Ici  $x_{2t_1+B}=B-2t_{1}-B+1=-2t_{1}+1$. On cherche \`a
d\'emontrer que $t_1=0$.  Et on a encore
$[B=2,A]
\cup_{i; z_i\leq 0} \{-z_i\}=
[D+1,A]\cup_{i;z_i>0}\{z_i-1\}$ et l'ensemble de droite n'a pas de multiplicit\'e. Ainsi $x_v\geq -1$ et  $1-2t_1 \geq -1$ c'est-\`a-dire $t_1=0$ ou $1$. Et on veut d\'emontrer que $t_1=0$. Il reste donc \`a \'eliminer le cas o\`u $t_1=1$; supposons donc que $t_{1}=1$  alors $x_v=x_{2t_1+2}=-1$, le tableau repr\'esentant ${\cal E}$ est rectangulaire avec 4 lignes.  Comme le bout de chaque ligne est sup\'erieur ou \'egal \`a $D+1$, $Jac_x \sigma_{\cal E}^*=0$ pour tout $x>-(D+1)$. Supposons que $\zeta' 1 >-(D+1)$, on a d'une part avec \ref{technique} appliqu\'e comme dans les r\'eductions ci-dessus que $Jac_{\zeta' 1}\tilde{\pi}\neq 0$ et d'autre part  avec l'inclusion
$\tilde{\pi}\hookrightarrow \sigma_{\cal E}\times \pi(\tilde{\psi}_D,\tilde{\epsilon}_{D,\lambda})$et ce qui pr\'ec\`ede $Jac_{\zeta' 1}\pi(\tilde{\psi}_D,\tilde{\epsilon}_{D,\lambda}\neq 0$ ce qui est exclu. On peut quand m\^eme avoir $\zeta' 1 \leq -(D+1)$ ce qui est en fait $\zeta'=-$ et $D=0$ mais puisque $Jac_{-1}\tilde{\pi}\neq 0$, il faudrait que l'inclusion ci-dessus se factorise en
$\tilde{\pi}\hookrightarrow \rho\vert\,\vert^{-1}\times Y$,
o\`u $Y$ est un sous-quotient irr\'eductible de $\times_{x\in {\cal E}-\{-1\}}\rho\vert\,\vert^x \times \pi(\tilde{\psi}_{D},\tilde{\epsilon}_{D,\lambda})$. On remarque que ${\cal E}$ est en fait tr\`es particulier sous nos hypoth\`eses car n\'ec\'essairement $A=D+2t_{1}+2=0+2+2=4$ et $B=2$. D'o\`u
$$
{\cal E}':= {\cal E}-\{-1\}:=
\begin{matrix}
&2 &3 &4\\
&1 &2 &3\\
&0 &1 &2\\
&-1 &0&\\
\end{matrix}
$$
Quitte \`a changer l'ordre dans ${\cal E}'$ et des signes, on \'ecrit $Y$ comme sous-module de $\times_{x\in {\cal E}'}\rho\vert\,\vert^x \times \pi(\tilde{\psi}_{D},\tilde{\epsilon}_{D,\lambda})$ et un ${\cal E}$ qui convient est donc $\{-1\}\cup {\cal E}'$. Par minimalit\'e du nombre d'\'el\'ements positif dans notre ${\cal E}$ de d\'epart, ${\cal E}'$ ne peut \^etre constitu\'e que d'\'el\'ements positifs ou nul. Donc $-1$ a \'et\'e chang\'e en 1 et $\sigma_{\cal E}' \hookrightarrow \rho\vert\,\vert \times \sigma''$ pour $\sigma''$ une repr\'esentation bien choisi et comme l'induite $\rho\vert\,\vert^{-1}\times \rho\vert\,\vert$ est irr\'eductible, on a une inclusion:
$$
\tilde{\pi}\hookrightarrow \rho\vert\,\vert \times \pi'',
$$
pour une repr\'esentation $\pi''$ convenable, ce qui entra\^{\i}ne $Jac_{1}\tilde{\pi}\neq 0$; ceci est contradictoire avec \ref{proprietedujac} puisque $\zeta'=-1$ ici.

On en est donc revenu au cas o\`u ${\cal E}$  ne s'\'ecrit pas comme tableau comme pr\'ec\'edemment. On v\'erifie que ${\cal E}$ s'\'ecrit comme union  d'un tableau:
$${\cal T}=
\begin{matrix}
&2 &\cdots &\cdots &\cdots &A\\
&\vdots &\vdots &\vdots &\vdots &\vdots\\
&x_{2t_1+2} &\cdots &\cdots &\cdots &D+1\\
&\vdots &\vdots &\vdots&\vdots\\
&x_v &\cdots &y_v
\end{matrix}
$$
avec l'ensemble $\{\zeta' 1, 0\}$. On a encore en
consid\'erant ${\cal E}\cup -{\cal E}$ mais en appelant ici
$z_i$ les derniers \'el\'ements de chaque colonne de ${\cal
T}$:
$$
\{0,1\} \cup [2,A] \cup_{i;z_i\leq 0}
\{-z_i\}=[D+1, A] \cup_{i;z_i>0} \{z_i-1\}.
$$
Et toujours le fait que l'ensemble de droite n'ayant pas de
multiplicit\'e, celui de gauche ne peut en avoir. Donc ici
l'ensemble $\{i;z_i\leq 0\}$ est n\'ecessairement vide. Cela
force $x_{2t_1+2} \geq 1$ d'o\`u, par un calcul d\'ej\`a
fait $t_1=0$. Ce qui est ce que l'on cherche.

\section{Rappels\label{rappel}}

\subsection{Les s\'eries discr\`etes et hypoth\`eses}

Dans le cas des s\'eries discr\`etes on a d\'emontr\'e en
\cite{europe} et \cite{ams} les r\'esultats suivants que je vais rappel\'e ci-dessous. Soit
$\psi,\epsilon$ un couple form\'e d'un morphisme de
$W_F\times SL(2,{\mathbb C})$ avec les propri\'et\'es
d'alg\'ebricit\'e usuelles dans
$^LG$ tel que
$Cent_{^LG}(\psi)$ soit fini et $\epsilon$ un caract\`ere de
$Cent_{^LG}\psi$ dans $\pm 1$ tel que la restriction de
$\epsilon$ au centre de $^LG$ soit $\sharp$ (cf.
\ref{groupe}). En plongeant $^LG$ dans un GL par son
application naturelle, on d\'ecompose $\psi$ en somme de
repr\'esentation irr\'eductible de la forme $\rho\otimes
\sigma_a$ o\`u
$\rho$ est une repr\'esentation irr\'eductible de $W_F$,
n\'ecessairement autoduale et
$\sigma_a$ est une repr\'esentation irr\'eductible de
$SL(2,{\mathbb C})$ n\'ecessairement uniquement
d\'etermin\'e par sa dimension, not\'ee $a$ en indice; la
parit\'e de $a$ est uniquement d\'etermin\'ee par $\rho$.
C'est de l'alg\`ebre lin\'eaire. On note
$Jord(\psi):=\{(\rho,a)\}$ pour $(\rho,a)$ apparaissant dans
la d\'ecomposition ci-dessus; le caract\`ere $\epsilon$
s'identifie alors naturellement \`a une application de
$Jord(\psi)$ dans $\{\pm 1\}$. A un tel couple $\psi,\epsilon$ on associe un couple de m\^eme type mais cuspidale, c'est ce qui correspond \`a ce que l'on a appel\'e en loc.cit. le support cuspdial partiel. Une fa\c{c}on de le construire, sans utiliser de r\'ecurrence, peut se faire ainsi.

On construit d'abord un couple
$\psi_{0},\epsilon_{0}$ relatif \`a un groupe de
m\^eme type que $G$ mais de rang \'eventuellement plus petit
tel que $Jord(\psi_{0})$ se d\'eduise de $Jord(\psi)$ en
enlevant des blocs $(\rho_i,a_i)$ pour $i\in [1,t]$ tous
distincts o\`u
$t$ est un entier pair \'eventuellement $0$ soumis aux
conditions ci-dessous, $\epsilon_0$ est la restriction de
$\epsilon$ aux blocs restant; les conditions sont que pour
tout
$i\in [1,t/2]$,

$\rho_{2i-1}=\rho_{2i}$,
$\epsilon(\rho,a_{2i-1})=\epsilon(\rho,a_{2i})$, 
$a_{2i}<a_{2i-1}$ (on accepte
$a_{2i}=0$ au plus une fois en posant
$\epsilon(\rho,a_{2i})=+$ mais il faut alors que $a_{2i-1}$
soit pair) et pour tout
$b\in ]a_{2i},a_{2i-1}[$ tel que
$(\rho_{2i-1},b)\in Jord(\psi)$, $\exists j<i$ tel que
$\rho_{2j-1}=\rho_{2i-1}$ et $b=a_{2j}$ ou $a_{2j-1}$.
\nl
On remarque que $\psi_0,\epsilon_0$ n'est pas uniquement
d\'etermin\'e par $\psi,\epsilon$. Mais on fait un tel choix
en supposant que $\vert Jord(\psi_0)\vert$ est minimal. Pour
$\rho$ une repr\'esentation autoduale de
$W_F$, on note
$Jord_\rho(\psi_0):=\{a\in \mathbb{N}; (\rho,a)\in
Jord(\psi_0)$ et on pose $\delta_\rho:=\vert
Jord_\rho(\psi_0)\vert$ et $\eta_\rho$ la parit\'e de tout
\'el\'ement de $Jord_\rho(\psi_0)$, c'est-\`a-dire
$\eta_\rho=0$ si ces nombres sont pairs et $1$ sinon. On
note encore $\epsilon_\rho$ la valeur de $\epsilon_0$ sur
l'\'el\'ement minimal de $Jord_\rho(\psi_0)$; par
construction ce nombre vaut $-1$ si $\eta_\rho=0$ mais n'est
pas connu a priori si $\eta_\rho=1$.
Il est facile de v\'erifier que les donn\'ees $\delta_\rho$
et $\epsilon_\rho$ sont elles uniquement d\'etermin\'ees par
$\psi,\epsilon$.

On pose alors, a priori, $\psi_{cusp},\epsilon_{cusp}$ le
couple relatif \`a un groupe de m\^eme type que $G$ mais de
rang \'eventuellement plus petit, tel que
$Jord(\psi_{cusp}):=\{(\rho,\alpha); \alpha\leq \delta_\rho,
\alpha\equiv \eta_\rho[2]\}$, o\`u $\rho$ parcourt
l'ensemble des repr\'esentations autoduales de $W_F$ tel que
$\delta_\rho\neq 0$ et o\`u, pour tout $(\rho,\alpha)\in
Jord(\psi_{cusp})$, 
$\epsilon_{cusp}(\rho,\alpha):=(-1)^{(\alpha+\eta_\rho)/2+1}
\epsilon_\rho$. On note $G_{cusp}$ le groupe correspondant.
\nl
Le couple $\psi_{cusp},\epsilon_{cusp}$ est uniquement
d\'etermin\'e par $\psi,\epsilon$ et est sa donn\'ee
cuspidal.

\nl \sl
On suppose qu'il existe une repr\'esentation cuspidale
$\pi_{cusp}$ de $G_{cusp}$ avec la propri\'et\'e suivante
pour toute repr\'esentation autoduale $\rho$ de $W_F$
identifi\'ee \`a une repr\'esentation de $\rho$ cuspidale
autoduale de $GL(d_\rho)$ (ce qui d\'efinit $d_\rho$) par la
correspondance de Langlands (\cite{harris},\cite{henniart})
et pour tout entier $f$ tel que $f\equiv \eta_\rho[2]$
l'induite $St(\rho,f)\times \pi_{cusp}$ est irr\'eductible
si et seulement si $(\rho,f)\in Jord(\psi_{cusp})$; ici
$St(\rho,f)$ est la repr\'esentation de Steinberg qui est
plus g\'en\'eralement not\'ee
$<\rho\vert\,\vert^{-(f-1)/2},\cdots,
\rho\vert\,\vert^{(f-1)/2}>$ dans ce papier.\rm
\nl
Avec cette hypoth\`ese, on a montr\'e en \cite{algebra} que
$\pi_{cusp}$ avait les bonnes propri\'et\'es de
r\'eductibilit\'e suppos\'ees en \cite{europe} et
\cite{ams}. On peut donc utiliser les r\'esultats de ces
papiers; on a montr\'e qu'il existe une unique s\'erie
discr\`ete
$\pi(\psi,\epsilon)$ (irr\'eductible) telle que pour tout
$\rho$ comme ci-dessus et pour tout $f$ aussi comme
ci-dessus l'induite
$St(\rho,f)\times \pi(\psi,\epsilon)$ est irr\'eductible si
et seulement si $(\rho,f)\in Jord(\psi)$; de plus pour
$\rho$ comme ci-dessus et pour tout couple $(a>a_-)$
d'entiers cons\'ecutifs dans $Jord_\rho(\psi)$,

$\epsilon(\rho,a)=\epsilon(\rho,a_-)$ $\Leftrightarrow$
$Jac_{(a-1)/2, \cdots, (a_-+1)/2}\pi(\psi,\epsilon)\neq 0$

 et le support
cuspidal partiel de $\pi(\psi,\epsilon)$ est $\pi_{cusp}$ ce
qui veut dire qu'il existe une repr\'esentation $\sigma$
d'un groupe lin\'eaire de rang 1/2 la diff\'erence des rangs
de $G$ et $G_{cusp}$ tel que $\pi(\psi,\epsilon)$ soit un
sous-quotient de l'induite $\sigma\times \pi_{cusp}$.
L'existence est dans \cite{ams} et l'unicit\'e dans
\cite{europe}.

\subsection{Paquets associ\'es \`a des morphismes
\'el\'ementaires} 
Fixons encore $\psi,\epsilon$ comme
ci-dessus d'o\`u
$\psi_{cusp},\epsilon_{cusp}$ et on fait l'hypoth\`ese
cl\'e. Donnons-nous aussi une application de $Jord(\psi)$
dans $\{\pm\}$; cela permet de d\'efinir un morphisme de
$W_F\times SL(2,{\mathbb C})\times SL(2,{\mathbb C})$ dans
$^LG$ de telle sorte que la d\'ecomposition de ce morphisme
en repr\'esentations irr\'eductibles (comme plus haut) est
somme des produits tensoriels $\rho\otimes \sigma_a\otimes
\sigma_b$ o\`u  $(\rho,sup(a,b))$ parcourt $Jord(\psi)$,
$inf(a,b)=1$  et $sup(a,b)=a$ si $\zeta=+$ et $b$ si
$\zeta=-$. On note $\psi_\zeta$ ce morphisme et $\epsilon$
s'identifie encore naturellement \`a un caract\`ere de
$Cent_{^LG}\psi_\zeta$ dans $\{\pm 1\}$. Donc ici $$\psi=\psi_{\zeta}
\circ \Delta.\eqno(*)$$
Dans
\cite{paquet}, on a associ\'e une repr\'esentation
irr\'eductible $\pi(\psi_\zeta,\epsilon)$ \`a un tel couple.
On en a donn\'e 2 d\'efinitions; celle qui donne le plus de
propri\'et\'e est une d\'efinition dans le groupe de
Grothendieck qui g\'en\'eralise la formule d'Aubert
\cite{aubert} et Schneider-Stuhler \cite{SS} pour
l'involution g\'en\'eralisant celle de Zelevinski dans le
cas des groupes lin\'eaires. Pour l'expliquer, nous avons
besoin de la notation suivante: soit $a$ un demi-entier
positif, $\rho$ comme ci-dessus et $\pi$ une repr\'esentation
de
$G$; soit
$P$ un sous-groupe parabolique de $G$ donc de la forme
$\times_{d\in {\cal D}}GL(d)\times G(m)$ o\`u ${\cal D}$ est
un ensemble d'entiers. On note $res_{P,\rho,\leq a}(\pi)$ la
projection de la restriction de $\pi$ le long du radical
unipotent de $P$ sur le support cuspidal pour $\times_{d\in
{\cal D}}Gl(d)$ form\'e de repr\'esentation du type
$\times_{z \in \mathbb{R}}\rho^z$ mais o\`u tous les $z$ qui
apparaissent v\'erifient $\vert z\vert \leq a$. On d\'efinit
de m\^eme $res_{P,\rho,<a}$ en imposant l'in\'egalit\'e
stricte. On pose alors $$inv_{\rho,\leq a}
=\sum_{P}(-1)^{rg_G-rg_P} ind_P^G (res_{P,\rho,\leq a}\pi),$$
c'est une application dans le groupe de Grothendieck des
repr\'esentations lisses de type fini de $G$. On d\'efinit
de fa\c con analogue $inv_{\rho,<a}$ en rempla\c cant
l'in\'egalit\'e large par l'in\'egalit\'e stricte. Et alors:
$$
\pi(\psi_\zeta,\epsilon)=\biggl(\prod_{\begin{array}{l}(\rho,a)\in
Jord(\psi);\\
\zeta (\rho,a)=-
\end{array}} inv_{\rho,<a}\circ inv_{\rho,\leq a}\biggr)
(\pi(\psi,\epsilon)).
$$
 L'ordre dans lequel on effectue ces applications
n'a pas d'importance. Il n'est pas clair du tout sur cette
d\'efinition que $\pi(\psi_\zeta,\epsilon)$ est
irr\'eductible mais c'est un des r\'esultats de
\cite{paquet}. En fait on montre que cette d\'efinition
est la m\^eme que la d\'efinition suivante qui se fait par
induction:

pour toute repr\'esentation $\rho$ fix\'ee, on note
$a_{\rho,\psi,\epsilon}$ l'entier minimum, quand il existe,
v\'erifiant

$(\rho,a_{\rho,\psi,\epsilon})\in Jord(\psi)$,

soit $(\rho,a_{\rho,\psi,\epsilon}-2)\notin Jord(\psi)$,

soit $(\rho,a_{\rho,\psi,\epsilon}-2) \in Jord(\psi)$ (en
acceptant \'eventuellement $a_{\rho,\psi,\epsilon}=2$ et
v\'erifie $$\epsilon(\rho,a_{\rho,\psi,\epsilon})=\epsilon(
\rho,a_{\rho,\psi,\epsilon}-2).$$

On a alors les propri\'et\'es suivantes; dire que
$(\psi,\epsilon)\neq (\psi_{cusp},\epsilon_{cusp})$ est
exactement \'equivalent \`a dire qu'il existe $\rho$ tel que
$a_{\rho,\psi,\epsilon}$ existe. Fixons un tel $\rho$.

\bf propri\'et\'e 1: \rm alors
$\pi(\psi_{\zeta},\epsilon)$ ne d\'epend pas de la valeur de $\zeta$
sur les blocs $(\rho,a) \in Jord(\psi)$ tel que
$a<a_{\rho,\psi,\epsilon}$. On pose
$\zeta_{\rho,\psi,\epsilon}$ la valeur de $\zeta$ sur
$(\rho,a_{\rho,\psi,\epsilon})$;

\bf propri\'et\'e 2: \rm supposons que
$(\rho,a_{\rho,\psi,\epsilon}-2)\notin Jord(\psi)$, alors on
note $\psi',\epsilon'$ le couple qui se d\'eduit de
$\psi,\epsilon$ en changeant simplement le bloc
$(\rho,a_{\rho,\psi,\epsilon})$ en
$(\rho,a_{\rho,\psi,\epsilon}-2)$. Alors l'induite:
$$
\rho\vert\,\vert^{\zeta_{\rho,\psi,\epsilon}
(a_{\rho,\psi,\epsilon}-1)/2} \times \pi(\psi'_{\zeta},\epsilon')$$a
un unique sous-module irr\'eductible qui est pr\'ecis\'ement
$\pi(\psi_{\zeta},\epsilon)$; ceci est aussi vrai si
$a_{\rho,\psi,\epsilon}=2$ avec n\'ecessairement
$\epsilon(\rho,a_{\rho,\psi,\epsilon})=+$;

\bf propri\'et\'e 3: \rm supposons que
$(\rho,a_{\rho,\psi,\epsilon}-2)\in Jord(\psi)$ avec
$a_{\rho,\psi,\epsilon}-2\neq 0$. On note $\psi',\epsilon'$
le couple qui se d\'eduit de $\psi,\epsilon$ en enlevant
les 2 blocs $(\rho,a_{\rho,\psi,\epsilon})$ et $(\rho,
a_{\rho,\psi,\epsilon}-2)$. L'induite:
$$
<\rho\vert\,\vert^{\zeta_{\rho,\psi,\epsilon}(
a_{\rho,\psi,\epsilon}-1)/2}, \cdots, \rho\vert\,\vert^{  
-\zeta_{\rho,\psi,\epsilon}(a_{\rho,\psi,\epsilon}-3)/2}>
\times
\pi(\psi'_{\zeta},\epsilon')$$
a exactement 2 sous-modules irr\'eductibles et
$\pi(\psi_{\zeta},\epsilon)$ est l'un des deux. Le choix est
pr\'ecis\'e en \cite{paquet}, il y a une part
d'arbitraire et nous n'avons pas forc\'ement fait le
meilleur choix. Mais cela n'a aucun importance pour ce que
nous faisons.

\bf propri\'et\'e 4: \rm soit ${\cal E}$ un ensemble de
demi-entiers tels que pour tout $x\in {\cal E}$, $\vert
x\vert < (a_{\rho,\psi,\epsilon}-1)/2$. Alors pour tout
sous-quotient irr\'eductible, $\tilde{\pi}$ de l'induite $\times_{x\in {\cal E}}\rho\vert\,\vert^x\times \pi(\psi_{\zeta},\epsilon)$, il existe un
ensemble ${\cal E}'$ totalement ordonn\'e de demi-entiers
v\'erifiant:
$$
\{{\cal E}\} \cup \{-{\cal E}\}=
\{{\cal E}'\} \cup \{-{\cal E}'\}$$
et une inclusion: $\tilde{\pi}\hookrightarrow \times_{x\in
{\cal E}'}\rho\vert\,\vert^x\times \pi(\psi_{\zeta},\epsilon).$
\nl
On remarque que les propri\'et\'es 1, 2 et 3 permettent de
d\'efinir $\pi(\psi_{\zeta},\epsilon)$ par induction et le travail
\cite{paquet} consiste justement \`a prouver que cette
d\'efinition et celle donn\'ee dans le groupe de
Grothendieck co\"{\i}ncident.

\bf propri\'et\'e 5: \rm soit $x\in {\mathbb R}$; on suppose que $Jac_{x}\pi(\psi_{\zeta},\epsilon)\neq 0$, alors il existe $a\in {\mathbb N}$  tel que $(\rho,a)\in Jord(\psi)$ et $x=\zeta_{(\rho,a)}(a-1)/2$. R\'eciproquement si $(\rho,a)\in Jord(\psi)$ et $(\rho,a-2)\notin Jord(\psi)$ alors $Jac_{\zeta_{(\rho,a)}(a-1)/2}\pi(\psi_{\zeta},\epsilon)\neq 0$. 

\bf propri\'et\'e 6: \rm
 fixons $\rho$ comme
ci-dessus et $a>a_-$ des entiers cons\'ecutifs de
$Jord_\rho(\psi)$; on suppose ici que
$\zeta(\rho,a)=\zeta(\rho,a_-)$ alors $$Jac_{\zeta(a-1)/2,
\cdots, \zeta(a_-+1)/2}\pi(\psi_\zeta,\epsilon)\neq 0
\Leftrightarrow \epsilon(\rho,a)=\epsilon(\rho,a_-).$$

L'analogue de la propri\'et\'e 6 pour les s\'eries discr\`etes est essentiellement la d\'efinition de la param\'etrisation (cf \cite{europe}); en utilisant la formule dans le groupe de Grothendieck d\'efinissant $\pi(\psi_{\zeta}, \epsilon)$ on obtient la propri\'et\'e \'ecrite mais elle est  faible puisqu'il faut
supposer au d\'epart que $\zeta(\rho,a)=\zeta(\rho,a_-)$.
Sans cette hypoth\`ese, on a quand m\^eme une propri\'et\'e
des modules de Jacquet qui distingue le cas
$\epsilon(\rho,a)=\epsilon(\rho,a_-)$ mais c'est plus
compliqu\'e; on l'\'ecrit simplement dans le cas
que nous utiliserons.
On suppose que $Jord_\rho(\psi)$ contient $1,3$ et $5$ et
que $\zeta:=\zeta(\rho,1)=\zeta(\rho,3)\neq \zeta(\rho,5)$ et
$\epsilon(\rho,1)=\epsilon(\rho,3)$. Alors
$\epsilon(\rho,3)=\epsilon(\rho,5)$ si et seulement si
$Jac_{-\zeta 2,-\zeta 1}\pi(\psi_\zeta,\epsilon)\neq 0$.

\subsection{Propri\'et\'e d'induction\label{induction}}
Ici on revient \`a $\psi,\epsilon$ non n\'ecessairement
\'el\'ementaire, c'est-\`a-dire comme dans tout ce travail.
\nl
\bf Proposition: \sl
Soit $(\psi,\epsilon)$ et soit $b$ un entier. Soit aussi
$\zeta=\pm$; on suppose que $(\rho,b)\notin Jord(\psi\circ
\Delta)$ et que $\psi\circ \Delta$ est discret. L'induite $<\rho\vert\,\vert^{\zeta(b-1)/2},
\cdots, \rho\vert\,\vert^{-\zeta(b-1)/2}>\times
\pi(\psi,\epsilon)$ est semi-simple de longueur 2 fois la
longueur de $\pi(\psi,\epsilon)$.\rm
\nl
On la fait pour les s\'eries discr\`etes en \cite{europe}, mais il faut bien tenir compte du fait que le point cl\'e est alors le cas cuspidale o\`u en fait on admet le r\'esultat dans les hypoth\`eses de base. On en d\'eduit le cas de $\psi$ \'el\'ementaire: si $a_{\rho,\psi,\epsilon}>b$, on l'a d\'emontr\'e dans \cite{paquet} et on en d\'eduit le cas g\'en\'eral par r\'ecurrence en \'echangeant socle et induction comme ci-dessous; on ne fait pas les d\'etails. Supposons donc
qu'il existe
$(\rho,a,b)\in Jord(\psi)$ tel que $inf(a,b)\geq 2$. On
suppose pour simplifier que $a\geq b$. On  d\'emontre l'assertion par
r\'ecurrence. On sait (\ref{enonce}) que $\pi(\psi,\epsilon)$
est la somme des sous-modules irr\'eductibles inclus dans
$<\rho\vert\,\vert^{(a-b)/2},\cdots,\rho\vert\,\vert^{-(a+b)/2
+1}>\times \pi(\psi',\epsilon',
(\rho,a,b-2,\epsilon(\rho,a,b)))$ et de repr\'esentation
$\pi(\psi'',\epsilon'')$ pour un
bon choix de morphismes $\psi'',\epsilon''$ auxquels on peut appliquer la r\'ecurrence.

 On pose $\sigma:=
<\rho\vert\,\vert^{\zeta(b-1)/2},
\cdots, \rho\vert\,\vert^{-\zeta(b-1)/2}>$.
On sait que
$\pi(\psi',\epsilon',(\rho,a,b-2,\epsilon(\rho,a,b))$ est
semi-simple; on note $Y'$ l'un de ses sous-modules
irr\'eductibles et $Y$ l'unique sous-module irr\'eductible
inclus dans l'induite $$
<\rho\vert\,\vert^{(a-b)/2},\cdots, \rho\vert\,\vert^{(a+b)/
2-1}>\times Y'.
$$Ce qui reste \`a d\'emontrer est que $\sigma\times Y$ est semi-simple de longueur 2.
On applique le lemme technique, \ref{technique}, \`a
$\sigma\times Y$ avec
$\delta=\sigma$ et $\delta'=
<\rho\vert\,\vert^{(a-b)/2},\cdots, \rho\vert\,\vert^{-(a+b)/
2+1}>$; pour cela on v\'erifie que $\delta\times \delta'$ et
$\delta\times (\delta')^*$ sont irr\'eductibles. Si
$\zeta=+$, cela r\'esulte de Zelevinski car les segments
correspondant ne sont pas li\'es. Si $\zeta=-$, $\delta$
correspond \`a un segment qui n'a pas la m\^eme croissance
que celui qui correspond \`a $\delta'$; il faut donc un
autre argument. On a soit $\beta< (a-b+1)$ soit $\beta
>a+b-1$. Dans le premier cas, pour tout $x\in [(\beta-1)/2,
-(\beta-1)/2]$, le segment r\'eduit \`a $x$ n'est pas li\'e
au segment correspondant \`a $\delta'$ ni \`a celui
correspondant \`a $(\delta')^*$. Dans le deuxi\`eme cas, 
pour tout $y \in [(a-b)/2,-(a+b)/2-1]\cup
[-(a-b)/2,(a+b)/2-1]$, le segment r\'eduit \`a $y$ n'est pas
li\'e \`a $[-(\beta-1)/2,(\beta-1)/2]$. Cela suffit
largement pour avoir l'irr\'eductibilit\'e; en effet on veut
que l'op\'erateur d'entrelacement (normalis\'e \`a la
Langlands) $\delta\times \delta' \rightarrow \delta'\times
\delta$ soit un isomorphisme. Mais on le d\'ecompose en
op\'erateurs \'el\'ementaires qui eux sont n\'ecessairement
des isomorphismes par irr\'eductibilit\'e.

Les autres hypoth\`eses sont satisfaites. Ainsi on sait que
l'ensemble des sous-quotients irr\'eductibles de
$\sigma\times Y$ sont en bijection avec l'ensemble des
sous-quotients irr\'eductibles de
$\sigma\times Y'$. Par r\'ecurrence, on admet que cette
repr\'esentation est semi-simple. Il est facile de voir
qu'elle est alors de longueur au plus 2  et en fait elle
est de longueur au moins 2 car il faut que $\sigma\times
\pi(\psi',\epsilon', (\rho,a,b-2,\epsilon(\rho,a,b))$ double
sa longueur. Ainsi $\sigma\times Y$ a exactement 2
sous-quotients irr\'eductibles. Les sous-quotients, $Y'_i$
pour $i=1,2$ de
$\sigma\times Y'$ \'etant des sous-modules, ils v\'erifient
$Jac_{\zeta(\beta-1)/2, \cdots -\zeta(\beta-1)/2}Y'_i\neq
0$. On note $Y_i$ pour $i=1,2$, le sous-quotient de
$\sigma\times Y$ qui correspond \`a $Y'_i$ par
$Jac_{(a-b)/2, \cdots, -(a+b)/2}$ et n\'ecessairement
$Jac_{\zeta(\beta-1)/2,\cdots,-\zeta(\beta-1)/2}Y_i\neq 0$.
Ainsi $Y_i \hookrightarrow
\rho\vert\,\vert^{\zeta(\beta-1)/2}\times \cdots \times
\rho\vert\,\vert^{-\zeta(\beta-1)/2}\times Y$ et il existe
un sous-quotient $\sigma'$ de $\rho\vert\,\vert^{\zeta(\beta-1)/2}\times \cdots \times
\rho\vert\,\vert^{-\zeta(\beta-1)/2}$ tel que
$Y_i\hookrightarrow \sigma'\times Y$. On v\'erifie que
$\sigma'\simeq \sigma$; pour cela il suffit de montrer que
n\'ecessairement $Jac_{\rho\vert\,\vert^{\zeta(\beta-1)/2},
\cdots, \rho\vert\,\vert^{-\zeta(\beta-1)/2}}\sigma'\neq 0$.
Or s'il n'en est pas ainsi, il existe n\'ecessairement
$x\in ]\zeta(\beta-1)/2,-\zeta(\beta-1)/2]$ tel que
$Jac_{x,\cdots, -\zeta(\beta-1)/2}\sigma'\neq 0$; mais comme
$Y_i$ est un sous-quotient de $\sigma\times Y$, cela
entra\^{\i}ne que $Jac_{x,\cdots, -\zeta(\beta-1)/2} Y\neq
0$. Or ceci 
est exclu par \ref{techniquepoursocle}. D'o\`u le r\'esultat
cherch\'e.
\subsubsection{Irr\'eductibilit\'e\label{irreductibilite}}
On fixe $\psi,\epsilon$ comme dans tout ce travail ainsi
qu'une repr\'esentation cuspidale irr\'eductible $\rho$. On suppose encore ici que $\psi\circ \Delta$ est discret.
\nl
\bf Proposition: \sl Soit $x$ un demi-entier strictement
positif. On suppose que
$(\rho,2x-1)\notin Jord(\psi\circ \Delta)$. Alors la
repr\'esentation induite $\rho\vert\,\vert^x\times
\pi(\psi,\epsilon)$ est semi-simple de m\^eme longueur que
$\pi(\psi,\epsilon)$.\rm
\nl
La proposition est \'equivalente \`a dire que pour $X$ une
sous-repr\'esentation irr\'eductible de $\pi(\psi,\epsilon)$
l'induite $\rho\vert\,\vert^x\times X$ est irr\'eductible.
Ceci se d\'emontre par r\'ecurrence comme ci-dessus.

\section{D\'efinition de $\pi(\psi,\epsilon)$ dans le cas g\'en\'eral.\label{casgeneral}}
Soit ici un morphisme $\psi:W_{F}\times SL(2,{\mathbb C}) \times SL(2,{\mathbb C})$ dans $^LG$ et on ne suppose ici que le fait que $\psi$ est discret c'est-\`a-dire que $Cent_{^LG}\psi$ est un groupe fini. La d\'ecomposition de $\psi$ en repr\'esentation irr\'eductible donne toujours un ensemble $Jord(\psi)$; si on \'ecrit cet ensemble en terme de repr\'esentations cela donne un ensemble de triplet $(\rho,a,b)$ avec une condition de parit\'e reliant $\rho$ et $a+b$. Il est plus commode de poser $A=(a+b)/2-1$ et $B=(a-b)/2$ en notant $\zeta$ le signe de $a-b$ \'etant entendu qu'un nombre nul est positif. Ainsi $Jord(\psi)$ est vu comme un ensemble de quadruplet $(\rho,A,B,\zeta)$ pour $\rho$ une repr\'esentation irr\'eductible autodual de $W_{F}$, $\zeta$ un signe, $A,B$ des demi-entiers tels que $A-B\in {\mathbb Z}_{\geq 0}$ avec la condition que $A$ est entier non demi-entier pr\'ecis\'ement si $\rho$ est de m\^eme type que $^LG$ (le type \'etant  symplectique ou orthogonal) et $\zeta=+$ si $B=0$.

Fixons $\rho$ et d\'efinissons $Jord_{\rho}(\psi)$ comme l'ensemble des quadruplets $(\rho,A,B,\zeta)$ avec ce $\rho$. On ordonne totalement $Jord_{\rho}(\psi)$ par:
$
(\rho,A',B',\zeta')>(\rho,A,B,\zeta)$ si ces 2 quadruplets sont distincts et
\nl
soit $B'>B$; soit $B'=B$ mais $A'>A$; soit $B'=B,A'=A$ mais $\zeta'=+$ d'o\`u $\zeta=-$.
\nl
Soit $G'$ un groupe de m\^eme type que $G$ mais de rang plus grand et soit $\psi'$ un morphisme de m\^eme type que $\psi$ mais relativement \`a $G'$. On dit que $\psi'$ domine $\psi$ si pour tout $\rho$, il existe une bijection entre $Jord_{\rho}(\psi')$ et $Jord_{\rho}(\psi)$, bijection qui pr\'eserve l'ordre, ce qui la d\'etermine totalement et on la note $b_{\psi',\psi}$ et on demande qu'elle v\'erifie, pour tout $(\rho,A',B',\zeta')$ dans $Jord_{\rho}(\psi')$, en posant $(\rho,A,B,\zeta)$ sont image par $b_{\psi',\psi}$, 
$\zeta'=\zeta$, $A'-A=B'-B \geq 0$. 

\nl
Fixons $\psi'$ dominant $\psi$ et $(\rho,A,B,\zeta)\in Jord(\psi)$ avec $(\rho,A',B',\zeta)$ son image inverse par $b_{\psi',\psi}$; on d\'efinit ${\cal E}_{(\rho,A,B,\zeta}$ comme l'ensemble totalement ordonn\'e r\'eunion des segments $\cup_{\ell \in [1,A'-A]} [B'-\ell+1,A'-\ell+1]$, ensemble qui est vide si $A'=A$. On d\'efinit alors $Jac_{(\rho,A',B',\zeta)\mapsto (\rho, A,B,\zeta)}:=Jac_{x\in {\cal E}_{\rho,A,B,\zeta}}$.

Quitte \`a supposer le rang de $G'$ suffisamment grand, on peut ais\'ement construire de tels $\psi'$ dominant $\psi$ et on peut en plus imposer que $\psi'\circ \Delta$ soit discret. Gr\^ace \`a $b_{\psi',\psi}$, on identifie les groupes  $Cent_{^LG}(\psi)$ et $Cent_{^LG'}(\psi')$ ainsi que leurs groupes de caract\`eres. Ainsi pour $\epsilon$ un caract\`ere de $Cent_{^LG}(\psi)$ on sait d\'efinir $\pi(\psi',\epsilon)$. Et on pose, avec les notations pr\'ec\'edentes:
\nl
\bf D\'efinition: \sl $\pi(\psi,\epsilon)_{\psi'}:=\prod_{(\rho,A,B,\zeta)\in Jord(\psi)} Jac_{(\rho,A',B',\zeta)\mapsto (\rho,A,B,\zeta)}) \pi(\psi',\epsilon)$, o\`u dans le produit les $(\rho,A,B\zeta)$ sont rang\'es dans l'ordre croissant.\rm

\nl
\bf Propri\'et\'e: \sl $\pi(\psi,\epsilon)_{\psi'}$ ne d\'epend pas du choix de $\psi'$ dominant $\psi$ et tel que $\psi'\circ \Delta$ est discret.\rm
\nl
Cela r\'esulte de \ref{lecastroue}; en effet fixons $\psi'$ et $\psi''$ dominant $\psi$. On peut construire $\tilde{\psi}$ dominant \`a la fois $\psi'$ et $\psi''$ et a fortiori $\psi$. Il suffit de d\'emontrer que $\pi(\psi,\epsilon)_{\tilde{\psi}}=\pi(\psi,\epsilon)_{\psi'}$ et une \'egalit\'e analogue pour $\psi''$. Par sym\'etrie, il suffit de le d\'emontrer pour $\psi'$. On descend de $\tilde{\psi}$ vers $\psi'$ en appliquant progressivement les $Jac_{(\rho,\tilde{A},\tilde{B},\zeta)\mapsto (\rho,A',B',\zeta)}$ en commen\c{c}ant par les quadruplets les plus petits et donc de proche en proche on  est ramen\'e \`a d\'emontrer $\pi(\psi,\epsilon)_{\tilde{\psi}}=\pi(\psi,\epsilon)_{\psi'}$ quand $\tilde{\psi}$ et $\psi'$ sont tels qu'il existe un unique quadruplet $(\rho,A_{0}',B_{0}',\zeta_{0}') \in Jord(\tilde{\psi}$ tel que $(\rho, A_{0}'-1,B_{0}'-1,\zeta_{0}')\in Jord(\psi')$ et 
$$Jord(\tilde{\psi})-(\rho,A_{0}',B_{0}',\zeta_{0}')=Jord(\psi')-(\rho,A_{0}'-1,B_{0}'-1,\zeta_{0}').$$On note $(\rho, A_{0},B_{0},\zeta_{0})$ l'image de ces quadruplets par les bijections; elles ont n\'ecessairement la m\^eme image. Le premier point est de v\'erifier que 
$$
\prod_{(\rho,A,B,\zeta)\in Jord(\psi)} Jac_{(\rho,A',B',\zeta)\in Jord(\tilde{\psi})\mapsto (\rho,A,B,\zeta)}=$$
$$\prod_{(\rho,A,B,\zeta)\in Jord(\psi)} Jac_{(\rho,A',B',\zeta)\in Jord(\psi') \mapsto (\rho,A,B,\zeta)})Jac_{\zeta'_{0}B'_{0}, \cdots, \zeta'_{0} A'_{0}}.$$
En effet les ensembles 
$(\rho,A',B',\zeta)\mapsto (\rho,A,B,\zeta)$ sont les m\^emes pour $\tilde{\psi}$ et $\psi'$ sauf exactement quand $(\rho,A',B',\zeta)=(\rho,A'_{0},B'_{0},\zeta'_{0})$ o\`u par d\'efinition
$$
Jac_{(\rho,A'_{0},B'_{0},\zeta'_{0})\mapsto (\rho,A_{0},B_{0},\zeta'_{0})}=Jac_{(\rho,A'_{0}-1,B'_{0}-1,\zeta'_{0})\mapsto (\rho,A_{0},B_{0},\zeta'_{0})}Jac_{\zeta'_{0}B'_{0}, \cdots, \zeta'_{0} A'_{0}}.
$$
Il faut ensuite v\'erifier que $Jac_{\zeta'_{0}B'_{0}, \cdots, \zeta'_{0} A'_{0}}$ commute au-dessus des $
Jac_{(\rho,A',B',\zeta)\mapsto (\rho,A,B,\zeta)}$ pour tout $(\rho,A',B',\zeta)\in Jord(\psi')$ strictement plus petit que $(\rho,A'_{0}-1,B'_{0}-1,\zeta'_{0})$. Or, par hypoth\`ese $\psi'\circ \Delta$ est discret donc $B'_{0}-1>B'$, seule possibilit\'e dans la d\'efinition de l'ordre et l'hypoth\`ese $\psi'\circ \Delta $ discret renforce encore en $B'_{0}-1>A'$.  Soit $x$ un \'el\'ement de l'ensemble $(\rho,A',B',\zeta)\mapsto (\rho,A,B,\zeta)$ et pour tout \'el\'ement $y\in [\zeta'_{0}B'_{0}, \zeta'_{0}A'_{0}]$, on a s\^urement $\vert x\vert \leq A' < B'_{0}-1<\vert y\vert -1$. D'o\`u la commutation annonc\'ee (cf \ref{notation}).

Les hypoth\`eses de \ref{lecastroue} pour $\tilde{\psi},\epsilon$ et le quadruplet $(\rho,A'_{0},B'_{0}\zeta'_{0})$ sont pr\'ecis\'ement satisfaites car $\psi'$ est le $\psi_{-}$ de loc.cit. et par hypoth\`ese est tel que $\psi'\circ \Delta$ est discret. Donc on a $$Jac_{\zeta'_{0}B'_{0}, \cdots, \zeta'_{0} A'_{0}}\pi(\tilde{\psi},\epsilon)=\pi(\psi',\epsilon)$$
et
$$
\pi(\psi,\epsilon)_{\tilde{\psi}}=\prod_{(\rho,A,B,\zeta)\in Jord(\psi)} Jac_{(\rho,A',B',\zeta)\in Jord(\tilde{\psi})\mapsto (\rho,A,B,\zeta)} \pi(\tilde{\psi},\epsilon)=
$$
$$
\prod_{(\rho,A,B,\zeta)\in Jord(\psi)} Jac_{(\rho,A',B',\zeta)\in Jord(\psi') \mapsto (\rho,A,B,\zeta)} Jac_{\zeta'_{0}B'_{0}, \cdots, \zeta'_{0} A'_{0}}\pi(\tilde{\psi},\epsilon)=$$
$$
\prod_{(\rho,A,B,\zeta)\in Jord(\psi)} Jac_{(\rho,A',B',\zeta)\in Jord(\psi') \mapsto (\rho,A,B,\zeta)}\pi(\psi',\epsilon)
$$
ce qui est le r\'esultat cherch\'e $\pi(\psi,\epsilon)_{\psi'}$.
\nl
\bf D\'efinition: \sl on pose $\pi(\psi,\epsilon)=\pi(\psi,\epsilon)_{\psi'}$. \rm
\nl
\bf Propri\'et\'e: \sl Supposons que pour au moins un $\psi'$ dominant $\psi$v\'erifiant $\psi'\circ \Delta$ est discret, on ait
$$
\sum_{\epsilon}\epsilon\bigl(1_{W_{F}}\times 1_{SL(2,{\mathbb C}) }\times (\begin{matrix} &-1&0\\ &0 &-1 \end{matrix})\bigr)\pi(\psi',\epsilon)$$ est une distribution stable, il en est alors de m\^eme pour
$$
\sum_{\epsilon}\epsilon\bigl(1_{W_{F}}\times 1_{SL(2,{\mathbb C})} \times (\begin{matrix} &-1&0\\ &0 &-1 \end{matrix})\bigr) \pi(\psi,\epsilon),$$
la somme porte sur tous les caract\`eres $\epsilon$ du groupe $Cent_{^LG}\psi$.\rm
\nl
Cela est cons\'equence du fait que nos $Jac_{\cdots}$ pr\'eserve la stabilit\'e.

Remarquons que pour $G=SO(2n+1)$ si la restriction de $\psi$ \`a $W_{F}$ se factorise par le Frobenius, toutes nos hypoth\`eses sont satisfaites gr\^ace \`a \cite{algebra} qui montre que les constructions de Lusztig des s\'eries discr\`etes ont bien les propri\'et\'es que nous voulons et \`a \cite{inventiones}

\nl
Remarquons aussi que dans certains cas $\pi(\psi,\epsilon)$ peut \^etre 0; par exemple consid\'erons, $G=SO(9,F)$, d\'eploy\'e,  le morphisme $\psi$ tel que $\psi$ est trivial sur $W_{F}$ et ($Rep_{a}$ est par d\'efinition la repr\'esentation irr\'eductible de dimension $a$ de $SL(2,{\mathbb C})$)
$$
\psi=tr_{W_{F}}\otimes Rep_{4}\otimes tr_{SL(2,{\mathbb C})} \oplus
tr_{W_{F}}\otimes Rep_{2}\otimes tr_{SL(2,{\mathbb C})}\oplus
tr_{W_{F}}\otimes tr_{SL(2,{\mathbb C})}\otimes Rep_{2}.
$$
En d'autres termes $\psi$ est \'el\'ementaire et $$Jord(\psi)=\{(tr_{W_{F}},3/2,3/2,+),(tr_{W_{F}},1/2,1/2,+),(tr_{W_{F}},1/2,1/2,-)\}.$$
Notons $\psi'$ le morphisme qui v\'erifient:
$$
\psi=tr_{W_{F}}\otimes Rep_{6}\otimes tr_{SL(2,{\mathbb C})} \oplus
tr_{W_{F}}\otimes Rep_{4}\otimes tr_{SL(2,{\mathbb C})}\oplus
tr_{W_{F}}\otimes tr_{SL(2,{\mathbb C})}\otimes Rep_{2}.
$$On a $Jord(\psi')=\{(tr_{W_{F}},5/2,5/2,+),(tr_{W_{F}},3/2,3/2,+),(tr_{W_{F}},1/2,1/2,-)\}$. Consid\'erons le morphisme de $Jord(\psi)$ dans $\pm 1$ qui correpond dans l'ordre aux signes $(-,+,-)$. On sait que $\pi(\psi',\epsilon)$ est cuspidal, donc, avec nos d\'efinitions $\pi(\psi,\epsilon)=0$. Un calcul que j'ai fait avec Waldspurger prouve que c'est bien ce dont on a besoin pour avoir les propri\'et\'es d'endoscopie. C'est d'ailleurs cet exemple qui a motiv\'e la d\'efinition.

\end{document}